\documentclass[aap,seceqn,citesort,MSNbibl,dvips]{arximspdf}
\usepackage{mathbh}

\doi{10.1214/10-AAP683}
\volume{20}
\issue{6}
\pubyear{2010}
\firstpage{2204}
\lastpage{2260}

\makeatletter

\newtheorem{theorem}{Theorem}[section]
\newproclaim{remark}[theorem]{Remark}
\newtheorem{lemma}[theorem]{Lemma}
\newtheorem{cor}[theorem]{Corollary}
\newproclaim{defn}[theorem]{Definition}
\newproclaim{ass}{Assumption}[section]
\newtheorem{prop}[theorem]{Proposition}

\makeatother

\begin{document}
\begin{frontmatter}

\title{Fluid limits of many-server queues with reneging}
\runtitle{Fluid limits of many-servers queues with reneging}

\begin{aug}
\author[A]{\fnms{Weining} \snm{Kang}\corref{}\ead[label=e1]{wkang@umbc.edu}} and
\author[B]{\fnms{Kavita} \snm{Ramanan}\thanksref{t2}\ead[label=e2]{Kavita\_Ramanan@brown.edu}}
\runauthor{W. Kang and K. Ramanan}
\affiliation{University of Maryland, Baltimore County and
Brown University}
\address[A]{Department of Mathematics\\
\quad and Statistics \\
University of Maryland\\
Baltimore County \\
1000 Hilltop Circle\\
Baltimore, Maryland, 21250\\
USA\\
\printead{e1}} 
\address[B]{
Division of Applied Mathematics\\
Brown University\\
Providence, Rhode Island 02912\\
USA\\
\printead{e2}}
\end{aug}

\thankstext{t2}{Supported in part by the NSF Grants DMS-04-06191,
CMMI-0728064 and CMMI-0928154. Any opinions, findings and conclusions
or recommendations expressed in this material are those of the author(s)
and do not necessarily reflect those of the National Science Foundation.}

\received{\smonth{12} \syear{2008}}
\revised{\smonth{9} \syear{2009}}

%
\begin{abstract}
This work considers a many-server queueing system in which impatient
customers with i.i.d., generally distributed service times and i.i.d.,
generally distributed patience times enter service in the order of
arrival and abandon the queue if the time before possible entry into
service exceeds the patience time. The dynamics of the system is
represented in terms of a pair of measure-valued processes, one that
keeps track of the waiting times of the customers in queue and the
other that keeps track of the amounts of time each customer being
served has been in service. Under mild assumptions, essentially only
requiring that the service and reneging distributions have densities,
as both the arrival rate and the number of servers go to infinity, a
law of large numbers (or fluid) limit is established for this pair of
processes. The limit is shown to be the unique solution of a coupled
pair of deterministic integral equations that admits an explicit
representation. In addition, a fluid limit for the virtual waiting time
process is also established. This paper extends previous work by Kaspi
and Ramanan, which analyzed the model in the absence of reneging.
A~strong motivation for understanding performance in the presence of
reneging arises from models of call centers.
\end{abstract}

%

\end{frontmatter}

\tableofcontents

\section{Introduction} \goodbreak
\label{sec-intro}

\subsection{Background and motivation}
\label{subs-back}

We consider a many-server queueing system in which customers with independent,
identically\vadjust{\goodbreak} distributed (henceforth, i.i.d.) service requirements
chosen from
a general distribution are processed in the order of arrival.
In addition, a customer is assumed to abandon the queue if his/her time
spent waiting in queue reaches his/her patience time. The patience times
of customers are also assumed to be i.i.d. and drawn from a general
distribution.
When there are $N$ servers and the
cumulative customer arrival process is assumed to be a renewal process,
this reduces to the
so-called G/GI/N${}+{}$GI model.

Over the last couple of decades, several applications have spurred the
study of
many-server models with abandonment \cite{bachet81,boxwaal94,gansetal}.
Specifically, in applications to telephone contact
centers and (more generally) customer contact centers,
the effect of customers' impatience has been shown to
have a substantial impact on the performance of the system
\cite{gansetal}.
For example, customer abandonment can stabilize
a system that was formerly unstable.
Under the assumption that the interarrival,
service and abandonment time distributions are (possibly time-varying)
exponential, process-level
fluid and diffusion approximations
were obtained by Mandelbaum, Massey and Reiman \cite{manmasrei}
for the total number in system
in networks of multiserver queues
with abandonments and retrials.

On the other hand, for the case of Poisson arrivals, exponential
service times and
general abandonment distributions (the M/M/N${}+{}$GI queue),
explicit formulae for the steady state distributions of the queue
length and
virtual waiting time were obtained by Baccelli and Hebuterne
\cite{bachet81} (see Sections IV and V.2 therein),
whereas several other steady state performance measures and their
asymptotic approximations, in the limit as the arrival rates
and servers go to infinity, were derived by Mandelbaum and Zeltyn \cite
{manzel05}.
In addition, approximations for performance measures
suggested by these limit theorems were
used by Garnett et al. \cite{garmanrei02} and Mandelbaum and Zeltyn
\cite{manzel08} for the
case of exponential and general abandonment distributions, respectively,
to provide insight into the design of large call centers.

In all the previously mentioned works, the service times
were assumed to be exponential.
However, statistical analysis of real call centers has shown that both
service times and abandon times
are typically not exponentially distributed \cite{brownetal,manzel05},
thus providing strong motivation for considering many-server systems with
general service and abandonment distributions.
A step toward incorporating
more realistic general service distributions was taken in
the insightful paper by Whitt~\cite{whifluid06}, where a
deterministic fluid approximation for a G/GI/N${}+{}$GI queue with general
service and abandonment distributions was proposed.
However, the convergence of the discrete system starting empty
to this fluid approximation was left as a conjecture (see Conjecture
2.1 in
\cite{whifluid06}).
In this work,
we rigorously identify the functional law of large numbers or mean-field limit, as the number of servers goes to
infinity, of a many-server queueing system with general service and abandonment
distributions starting from general initial conditions. In a recent
work, Mandelbaum and Momcilovic
\cite{MM09} have established diffusion approximations for the queue-length and virtual waiting time processes in
a G/GI/N${}+{}$GI queue.

With a view to providing a Markovian representation of the dynamics
with a state space that is independent of the number of servers,
we introduce a pair of measure-valued processes to describe the
evolution of
the system. One measure-valued process keeps track of the waiting
times of customers in queue and the other keeps track of the
amounts of time each customer present in the system has been in service.
Under rather general assumptions (specified in Sections~\ref{subs-modyn}
and~\ref{subs-scaling}),
we establish an asymptotic limit theorem for the scaled (divided by $N$)
pair of measure-valued processes, as the number of servers $N$ and the mean
arrival rate into the system simultaneously go to infinity. In a recent
independent study,
Zhang \cite{Zhang} also considered the fluid limit for the same
G/GI/N${}+{}$GI system by using a measure-valued representation. His
approach is based on tracking the ``residual'' service and patience
times rather than tracking the ``ages'' in system and service as
considered in this work.
As in \cite{kasram07} and \cite{kasram08}, an advantage of the particular measure-valued
representation used here, in terms of ages in system and service,
rather than residual service and residual patience times, is that it
facilitates the
application of martingale techniques, which streamlines the analysis
and also
allows for a more intuitive representation of the dynamics of
the limiting process. In addition, the measure-valued approach also
simultaneously allows for the characterization of asymptotic limits of
several other functionals of interest.
In order to illustrate this point,
we also derive a limit theorem for the virtual waiting time of
a customer, defined to be the time before entry to service of a
(virtual) customer with infinite patience.

This work generalizes the framework of Kaspi and Ramanan
\cite{kasram07}, in which the corresponding model without abandonments
was considered. The presence of two coupled measure-valued processes,
rather than just one as in \cite{kasram07}, makes the analysis here
significantly more involved. In addition, an important step is the
identification of an explicit expression for the cumulative reneging
process. This paper also forms the basis of subsequent work in which we establish,
under suitable conditions, the convergence of the stationary distributions of the
fluid-scaled $N$-server systems to the invariant state of the
fluid limit, as $N$ tends to infinity \cite{kanram08c}.

It is worthwhile to mention that the models discussed above are
relevant when the mean demand
of customers is known (or can be accurately learned from an initial
period of measurements), which is a realistic assumption
in many applications.
In other scenarios, it may be more natural to model the demand
as being doubly stochastic. This approach was adopted
by Harrison and Zeevi \cite{harzee05} (see also \cite{basharzee05}),
who proposed optimal staffing and design of multi-class call centers
with several agent pools in the presence of abandonment
under the assumption that the dominant variability arises from the
randomness in the mean demand, rather than fluctuations around the
mean demand.

\subsection{Outline of the paper}

The outline of the paper is as follows.
We provide a more precise description of the model and the
measure-valued representation of the state, and describe the dynamical
equations governing the evolution of the system
in Section \ref{sec-mode} (the explicit
construction of the state process
is relegated to Appendix~\ref{ap-markov} and the strong Markov
property of the state process is established in Appendix \ref{sec:SMF}).
A key result here is
Theorem \ref{th-prelimit}, which provides a succinct
characterization of the state dynamics.
An analog of this characterization for continuous state
processes leads to
the fluid equations, which are
introduced in Section \ref{subs-fleqs}
(see Definition~\ref{def-fleqns}).
Next, the main results of the paper are summarized in
Section \ref{subs-mainres}. The first (Theorem \ref{thm:1})
is a uniqueness result that states that
(under the assumption that the service
and abandonment distributions have densities and finite first moments)
there exists at most one solution to the fluid equations.
The proof of this result, which is considerably more
involved than in the case without abandonment,
is the subject of Section \ref{sec-uniq}.
The second and main result of the paper (Theorem \ref{thm:2})
states that
under mild additional assumptions (namely,
Assumptions \ref{ass-init}--\ref{ass-h} introduced in
Section \ref{subs-scaling}),
the scaled sequence of state processes converges weakly to
the (unique) solution of the fluid equations, and provides
a fairly explicit representation for the solution.
The proof of this result consists of two main steps.
First, in Section \ref{Sec:relcom},
the sequence of scaled state processes is shown to be tight and
then, in Section \ref{sec:CSL},
it is shown that every subsequential limit is a solution to the
fluid equations.
Both of these results make use of properties of a family of
martingales that are established in Section \ref{subs-prelim}.
Finally, the last result (Theorem \ref{thm:3})
formulates the asymptotic limit theorem for
the virtual waiting time process, which is proved
in Section \ref{subs-prf3}.
To start with, in Section \ref{subs-notat},
we first collect some basic notation and
terminology used throughout the paper.

\subsection{Notation and terminology}
\label{subs-notat}

The following notation will be used throughout the paper.
${\mathbb Z}$ is the set of integers, ${\mathbb N}$ is the set of
strictly positive
integers, ${\mathbb R}$ is set of real numbers,
${\mathbb R}_+$ the set of nonnegative real numbers and ${\mathbb Z}_+$
is the set of
nonnegative integers.
For $a, b \in{\mathbb R}$, $a \vee b$ denotes the maximum of $a$ and $b$,
$a \wedge b$ the minimum of $a$ and $b$ and the short-hand $a^+$ is
used for $a \vee0$.
Given $A \subset{\mathbb R}$ and $a \in{\mathbb R}$, $A - a$ equals
the set
$\{x \in{\mathbb R}\dvtx x + a \in A\}$ and $\mathbh{1}_B$ denotes
the indicator
function of the set $B$
[i.e., $\mathbh{1}_B (x) = 1$ if $x \in B$ and $\mathbh{1}_B(x) = 0$
otherwise].

\subsubsection{Function and measure spaces}
\label{subsub-funmeas}

Given any metric space $E$, $\mathcal{C}_b(E)$ and $\mathcal{C}_c
(E)$ are, respectively,
the space of bounded, continuous functions and
the space of continuous real-valued functions with compact support
defined on
$E$, while $\mathcal{C}^1(E)$ is the space of real-valued,
once continuously differentiable functions on~$E$, and $\mathcal{C}^1_c(E)$
is the subspace of functions in $\mathcal{C}^1(E)$ that have compact support.
The subspace of functions in $\mathcal{C}^1(E)$ that, together with their
first derivatives, are bounded, will be denoted by $\mathcal
{C}^1_b(E)$. For
$H\le\infty$, let
$\mathcal{L}^1[0,H)$ and $\mathcal{L}^1_{\mathrm{loc}}[0,H)$, respectively,
represent the
spaces of
integrable and locally integrable functions on $[0,H)$, where a locally
integrable function $f$ on $[0,H)$
is a measurable function on $[0,H)$ that satisfies
$\int_{[0,a]}f(x)\,dx<\infty$ for all $a<H$. The constant functions
$f \equiv1$ and $f \equiv0$ will be represented by the symbols
${\mathbf{1}}$ and $\mathbf{0}$, respectively. Given any c\`{a}dl\`{a}g,
real-valued function $\varphi$ defined on $[0,\infty)$, we define
$\Vert\varphi\Vert_T \doteq\sup_{s \in[0,T]} |\varphi(s)|$ for
every $T < \infty$,
and let $\Vert\varphi\Vert_\infty\doteq\sup_{s \in[0,\infty)}
|\varphi(s)|$, which could possibly take the value $\infty$.
In addition,
the support of a function $\varphi$ is denoted by $\operatorname
{supp}(\varphi)$.
Given a
nondecreasing function $f$ on $[0,\infty)$, $f^{-1}$ denotes the inverse
function of $f$ in the sense that
%
%
\begin{equation}\label{inverse} f^{-1}(y)=\inf\{x\geq0\dvtx f(x)\geq
y\}.
\end{equation}
For each differentiable function $f$ defined on ${\mathbb R}$, $f'$
denotes the
first derivative of~$f$. For each function $f(t,x)$ defined on
${\mathbb R}
\times{\mathbb R}^n$, $f_t$ denotes the partial derivative of $f$ with respect
to $t$, and
$f_x$ denotes the partial derivative of $f$ with respect to $x$.

The space of Radon measures on a metric space $E$, endowed with the
Borel $\sigma$-algebra,
is denoted by $\mathcal{M}(E)$, while
$\mathcal{M}_F(E)$, $\mathcal{M}_1(E)$ and $\mathcal{M}_{\leq1}
(E)$ are, respectively, the subspaces of finite,
probability and sub-probability measures in $\mathcal{M}(E)$. Also,
given $B <
\infty$,
$\mathcal{M}_{\leq B} (E)\subset\mathcal{M}_F(E)$ denotes the space
of measures $\mu$ in $\mathcal{M}_F(E)$
such that $|\mu(E)| \leq B$.
Recall that a Radon measure is one that assigns finite measure to every
relatively
compact subset of ${\mathbb R}_+$.
The space $\mathcal{M}(E)$ is equipped with the vague topology, that
is, a
sequence of measures
$\{\mu_n\}$ in $\mathcal{M}(E)$ is said to converge to $\mu$ in the vague
topology (denoted
$\mu_n \stackrel{v}{\rightarrow}\mu$) if and only if for every
$\varphi\in\mathcal{C}_c (E)$,
%
%
\begin{equation}
\label{w-limit}
\int_{E} \varphi(x) \mu_n(dx) \rightarrow\int_E \varphi(x) \mu(dx)
\qquad\mbox{as } n \rightarrow\infty.
\end{equation}
By identifying a Radon measure $\mu\in\mathcal{M}(E)$ with the
mapping on
$\mathcal{C}_c (E)$ defined by
\[
\varphi\mapsto\int_{E} \varphi(x) \mu(dx),
\]
one can equivalently define a Radon measure on $E$ as a linear mapping
from $\mathcal{C}_c (E)$ into ${\mathbb R}$ such that for every
compact set $\mathcal{K}
\subset E$, there exists
$L_{\mathcal{K}} < \infty$ such that
\[
\biggl| \int_{E} \varphi(x) \mu(dx) \biggr| \leq L_{\mathcal{K}}
\Vert\varphi\Vert_\infty\qquad\forall
\varphi\in\mathcal{C}_c (E)\mbox{ with } \operatorname
{supp}(\varphi) \subset{\mathcal{K}}.
\]
%
On $\mathcal{M}_F(E)$, we will also consider the
weak topology, that is, a sequence $\{\mu_n\}$ in $\mathcal{M}_F(E)$
is said to
converge weakly to $\mu$ (denoted $\mu_n \stackrel{w}{\rightarrow
}\mu$) if and only
if (\ref{w-limit}) holds for every $\varphi\in\mathcal{C}_b(E)$.
As is
well known, $\mathcal{M}(E)$ and $\mathcal{M}_F(E)$, endowed with the
vague and weak
topologies, respectively,
are Polish spaces.
The symbol $\delta_x$ will be used to denote the measure with unit
mass at the point $x$, and,
by some abuse of notation, we will use $\mathbf{0}$ to denote the
identically zero
Radon measure on $E$. When $E$ is an interval, say $[0,H)$, for
notational conciseness,
we will often write $\mathcal{M}[0,H)$ instead of $\mathcal
{M}([0,H))$. For
any finite measure $\mu$ on $[0,H)$,
we define
%
%
\begin{equation}
\label{def-cdfmu}
F^\mu(x) \doteq\mu[0,x],\qquad x \in[0,H).
\end{equation}
We say a measure $\mu$ is continuous at $x$ if and only if $\mu(\{x\})=0$.

We will mostly be interested in the case when
$E = [0,H)$ and $E = [0,H) \times{\mathbb R}_+$, for some $H \in
(0,\infty]$.
To distinguish these cases, we will usually use $f$ to denote generic functions
on $[0,H)$ and $\varphi$ to denote generic
functions on $[0,H) \times{\mathbb R}_+$. By some abuse of notation,
given $f$ on $[0,H)$, we will sometimes
also treat it as a function on $[0,H) \times{\mathbb R}_+$ that is
constant in
the second
variable. For any
Borel measurable function $f\dvtx[0,H) \rightarrow{\mathbb R}$ that
is integrable
with respect to $\xi\in\mathcal{M}[0,H)$, we often use the
short-hand notation
\[
\langle f, \xi\rangle\doteq\int_{[0,H)} f(x) \xi(dx).
\]
%
Also, for ease of notation, given $\xi\in\mathcal{M}[0,H)$ and an
interval $(a,b)
\subset[0,H)$, we will
use $\xi(a,b)$ and $\xi(a)$ to denote $\xi((a,b))$ and $\xi(\{a\}
)$, respectively.

\subsubsection{Measure-valued stochastic processes}

Given a Polish space $\mathcal{H}$, we denote by $\mathcal
{D}_{\mathcal{H}}[0,T]$
(resp., $\mathcal{D}_{\mathcal{H}}[0,\infty)$) the space of
$\mathcal{H}$-valued,
c\`{a}dl\`{a}g functions on $[0,T]$ (resp., $[0,\infty)$),
and we endow this space with the usual Skorokhod $J_1$-topology
\cite{parbook}. Then $\mathcal{D}_{\mathcal{H}}[0,T]$ and $\mathcal
{D}_{\mathcal{H}}[0,\infty)$ are also Polish
spaces (see \cite{parbook}). In this work, we will be interested
in $\mathcal{H}$-valued stochastic processes, where $\mathcal{H} =
\mathcal{M}_F[0,H)$
for some $H\le\infty$.
These are random elements that are defined on a probability space
$(\Omega, \mathcal{F}, \mathbb{P})$ and take values in $\mathcal
{D}_{\mathcal{H}}[0,\infty)$,
equipped with the Borel $\sigma$-algebra (generated by open sets under the
Skorokhod $J_1$-topology).
A~sequence $\{X_n\}$ of c\`{a}dl\`{a}g, $\mathcal{H}$-valued processes,
with $X_n$ defined on the probability space $(\Omega_n, \mathcal
{F}_n, \mathbb{P}_n)$,
is said to converge in distribution
to a c\`{a}dl\`{a}g $\mathcal{H}$-valued process $X$ defined on
$(\Omega,
\mathcal{F}, \mathbb{P})$ if, for every bounded, continuous functional
$F\dvtx\mathcal{D}_{\mathcal{H}}[0,\infty)\rightarrow{\mathbb R}$,
we have
\[
\lim_{n \rightarrow\infty} \mathbb{E}_n [ F(X_n) ] = \mathbb{E}[
F(X) ],
\]
where $\mathbb{E}_n$ and $\mathbb{E}$ are the expectation operators
with respect to
the probability measures
$\mathbb{P}_n$ and $\mathbb{P}$, respectively.
Convergence in distribution of $X_n$ to $X$ will be denoted by $X_n
\Rightarrow X$. Let $\mathcal{I}_{{\mathbb R}_+}[0,\infty)$ be the
subset of nondecreasing functions
$f \in\mathcal{D}_{{\mathbb R}_+}[0,\infty)$ with $f(0) = 0$.

\section{Description of model and state dynamics}
\label{sec-mode}

In Section \ref{subs-modyn} we describe the basic model and the
model primitives.
In Section \ref{sec:repdyn} we introduce the state descriptor and some
auxiliary processes, and derive some equations that describe the
dynamics of
the state. Finally, in Section \ref{subs-prelimit} (see Theorem
\ref{th-prelimit}), we provide a succinct characterization of the
state dynamics.
This characterization motivates the
form of the fluid equations, which are introduced in Section
\ref{subs-fleqs}.

\subsection{Model description and primitive data}
\label{subs-modyn}

Consider a system with $N$ servers, in which arriving
customers are served in a nonidling, First-Come-First-Serve (FCFS)
manner, that is, a newly
arriving customer immediately enters service if there are any idle
servers or,
if all servers are busy, then the customer joins the back of
the queue, and the customer at the head of the queue (if one is
present) enters
service as soon as a server becomes free.
Our results are not sensitive to the exact mechanism used to assign an
arriving customer to an idle server, as long as the nonidling condition,
that there cannot simultaneously be a positive queue and an idle server,
is satisfied.
It is assumed that customers are impatient, and that a customer reneges
from the queue as soon
as the amount of time he/she has spent in queue reaches his/her
patience time.
Customers do not renege once they have entered service.
The patience times of customers are given by an i.i.d. sequence,
$\{r_i, i\in{\mathbb Z}\}$, with common cumulative distribution
function $G^r$ on $[0,\infty]$,
while the service requirements of customers are given by another
i.i.d.
sequence, $\{v_i, i\in{\mathbb Z}\}$, with common cumulative
distribution function
$G^s$ on $[0,\infty)$. For $i\in{\mathbb N}$, $r_i$ and $v_i$
represent, respectively, the
patience time and the service requirement of the $i$th customer to
enter the
system after time zero, while $\{r_i, i\in-{\mathbb N}\cup\{0\}\}$
and $\{v_i,
i\in-{\mathbb N}\cup\{0\}\}$ represent, respectively, the patience
times and the
service requirements of customers that arrived prior to time zero (if such
customers exist), ordered according to their arrival times
(prior to time zero). We assume that $G^s$ has density $g^s$ and $G^r$,
restricted on $[0,\infty)$, has density $g^r$.
This implies, in particular, that
$G^r(0+) = G^s(0+) = 0$. Let
\begin{eqnarray*}
H^r & \doteq& \sup\{x \in[0,\infty)\dvtx g^r(x) >0 \}, \\
H^s & \doteq
& \sup\{x \in[0,\infty)\dvtx g^s(x) >0 \}
\end{eqnarray*}
denote the right
ends of the supports of $g^r$ and $g^s$, respectively. The superscript
$(N)$ will be used to refer to quantities associated with the system
with $N$ servers.

Let $E^{(N)}$ denote the cumulative arrival process, with $E^{(N)}(t)$
representing the total number of customers that arrive into the system
with $N$ servers in the time interval $[0,t]$. Also, consider the c\`
{a}dl\`{a}g, real-valued process $\alpha_E^{(N)}$ defined by $\alpha
_E^{(N)}(s)=s$ if
$E^{(N)}(s)=0$ and, if $E^{(N)}(s)>0$, then
%
%
\begin{equation}
\label{def-ren}
\alpha_E^{(N)}(s) \doteq s- \sup\bigl\{u<s\dvtx E^{(N)}(u)<E^{(N)}(s) \bigr\},
\end{equation}
which denotes the time elapsed since the last arrival. If $E^{(N)}$ is
a renewal
process, then $\alpha_E^{(N)}$ is simply the backward recurrence time process.
Also, let $\mathcal{E}^{(N)}_0$ be an a.s. ${\mathbb Z}_+$-valued
random variable that
represents the
number of customers that entered the system prior to time zero.
This random variable does not play an important role in the analysis,
but is
used for bookkeeping purposes to keep track of the indices of customers.

The following mild assumptions on $E^{(N)}$ will be imposed throughout,
without explicit
mention:
\begin{itemize}
\item$E^{(N)}$ is a nondecreasing, pure jump process with $E^{(N)}(0)
= 0$
and a.s.,
for $t\in[0, \infty), E^{(N)}(t) < \infty$ and
$E^{(N)}(t)-E^{(N)}(t-) \in\{
0, 1\}$;
\item
the process $\alpha_E^{(N)}$ is Markovian with respect to its own natural
filtration (this holds, e.g., when $E^{(N)}$ is a renewal process);
\item
the cumulative arrival process $E^{(N)}$, the sequence of service requirements
$\{v_j, j \in{\mathbb Z}\}$ and the sequence of patience times
$\{r_j, j \in{\mathbb Z}\}$ are independent.
\end{itemize}
The assumption on the jump size of $E^{(N)}$ is not crucial and is imposed
mainly for convenience. On the other hand, the assumed
independence of the service and patience times is a genuine
restriction. It would be of interest to consider the case of correlated
service and patience times.

\subsection{State descriptor and dynamical equations}
\label{sec:repdyn}

As mentioned in Section \ref{subs-back}, our representation of the
state of the system with $N$ servers involves a pair of measure-valued
processes, the ``potential queue measure'' process, $\eta^{(N)}$,
which keeps track of the waiting times of customers
in queue and the ``age measure'' process, $\nu^{(N)}$,
which encodes the amounts of time
that customers currently receiving service have been in service.
In fact, the potential queue measure process
keeps track not only of the waiting times of customers in
queue, but also of the potential waiting times (equivalently,
times since entry into system) of those customers
who may have already entered service (and possibly departed the system),
but for whom the time since entry into the system has not yet exceeded the
patience time. In order to determine which subset of these customers
is actually in queue, the process $X^{(N)}$, which
represents the total number of
customers in system with $N$ servers (including those in service and
those in queue), is also incorporated into the state descriptor.
Thus the state of the system is represented by the
vector of processes $(\alpha_E^{(N)}, X^{(N)}, \nu^{(N)}, \eta
^{(N)})$, where
$\alpha_E^{(N)}$ determines the cumulative arrival process via (\ref
{def-ren}).
The reason for introducing the process
$\eta^{(N)}$ into the state (rather than working directly with a
restricted measure that only encodes the waiting times
of customers in queue) is that its dynamics is decoupled
from the service dynamics. It is governed
purely by the primitive data $E^{(N)}$ and $G^r$, and is thus
more amenable to analysis
(see Remark \ref{rem-compdyn} for further elaboration of this point).
Indeed, the queue measure process $\eta^{(N)}$ can also be viewed as
describing the ages of customers in an infinite server queue that
has cumulative arrivals $E^{(N)}$ and i.i.d. service requirements
distributed according to $G^r$. Thus the dynamics of the process
$\eta^{(N)}$ is also of independent interest.

Precise mathematical descriptions of $\eta^{(N)}$ and $\nu^{(N)}$
are given in Sections \ref{subsub-queue} and \ref{subsub-service},
respectively. Some auxiliary processes that are useful for describing the
evolution of the state are introduced in Section \ref{subsub-aux}.
Finally, in Section \ref{subsub-filt},
a~filtration $\{\mathcal{F}_t^{(N)}\}$ corresponding to the system
with $N$ servers is introduced, and it is shown that the state
processes and
auxiliary processes are all adapted to this filtration.
In fact, it is shown in Appendix \ref{sec:SMF} that
the state process is a strong Markov process with respect to this filtration.

\subsubsection{Description of queue dynamics}
\label{subsub-queue}

The potential waiting time process $w^{(N)}_j$ of customer $j$ is
(for every realization) defined to be the piecewise linear function
on $[0,\infty)$ that is identically zero till the customer enters the system,
then increases linearly, representing the amount of time elapsed since
entering the system,
and then remains constant (equal to the patience time)
once the time elapsed exceeds the patience time. More precisely,
for $j\in{\mathbb N}$, if $\zeta_j^{(N)}$ is the time at which the $j$th
customer arrives into the system after time $0$, then for $j \in
{\mathbb N}$
$\zeta_j^{(N)} = (E^{(N)})^{-1} (j) \doteq\inf\{ t > 0\dvtx E^{(N)} (t)
= j \}$
and
%
%
\begin{equation}
\label{def-waitjn}
w^{(N)}_j (t) = \cases{
\bigl[t - \zeta_j^{(N)} \bigr] \vee0, &\quad if $t - \zeta
{}^{(N)}_j < r_j$, \vspace*{2pt}\cr
r_j, &\quad otherwise.}
\end{equation}
For $j\in-{\mathbb N}\cup\{0\}$, $w^{(N)}_j$ represents the potential
waiting time
process of the $j$th customer who entered the system before time zero (if
such a customer exists). Observe that the potential waiting time
$w^{(N)}_j(t)$
of a customer at time $t$ equals its actual waiting time (equivalently,
time spent in queue) if and only if the customer has neither entered service
nor reneged by time $t$. For $t\in[0,\infty)$, let $\eta^{(N)}_t$ be the
nonnegative Borel measure on $[0,H^r)$ that has a unit mass at the potential
waiting time of each customer that has entered the system by time $t$ and
whose potential waiting time has not yet reached its patience time. Recall
that $\delta_x$ represents the Dirac mass at $x$. The potential queue
measure $\eta^{(N)}_t$ can be written in the form
%
%
\begin{eqnarray}
\label{def-etan} \eta^{(N)}_t &=& \sum_{j = -\mathcal{E}^{(N)}_0+
1}^{E^{(N)}(t)}
\delta_{w^{(N)}_j(t)} \mathbh{1}_{\{w^{(N)}_j (t) <
r_j\}}\nonumber\\[-8pt]\\[-8pt]
&=&\sum_{j = -\mathcal{E}^{(N)}_0+ 1}^{E^{(N)}(t)} \delta_{w^{(N)}_j(t)}
\mathbh{1}_{ \{{dw^{(N)}_j }/{dt}(t+) >0 \}},\nonumber
\end{eqnarray}
where the last equality holds because at any time $t$, the potential waiting
time process of any customer has a right derivative that is positive if and
only if the customer has entered the system and the customer's potential
waiting time has not yet reached its patience time.

For $t\in[0,\infty)$, let $Q^{(N)}(t)$ be the number of customers
waiting in queue
at time~$t$.
Due to the nonidling condition, the queue length process is then given by
%
%
\begin{equation}
\label{def-qnt}Q^{(N)}(t)=\bigl[X^{(N)}(t)-N\bigr]^+.
\end{equation}
Moreover, since the head-of-the-line customer is the customer in queue
with the
longest waiting time, the quantity
%
%
\begin{equation}\label{def-chi} \chi^{(N)}(t)\doteq
\inf\bigl\{x>0\dvtx\eta^{(N)}_t[0,x]\geq Q^{(N)}(t) \bigr\} =
\bigl(F^{\eta^{(N)}_t} \bigr)^{-1}
\bigl(Q^{(N)}(t)\bigr)
\end{equation}
represents the waiting time of the head-of-the-line customer in the
queue at
time $t$.
Here, recall from (\ref{def-cdfmu}) that $F^{\eta^{(N)}_t}$ is the
c.d.f. of the
measure $\eta^{(N)}_t$ and $(F^{\eta^{(N)}_t})^{-1}$ represents its
inverse, as defined in (\ref{inverse}).
Since this is an FCFS system,
any mass in $\eta^{(N)}_t$ that lies to the right of $\chi^{(N)}(t)$
represents
a customer that has already entered service by time $t$.
Therefore, the queue length process $Q^{(N)}$ admits the following
alternative representation in terms of $\chi^{(N)}$ and $\eta^{(N)}$:
%
%
\begin{eqnarray} \label{qn}
Q^{(N)}(t)&=& \sum_{j=-\mathcal{E}^{(N)}_0+ 1}^{E^{(N)}(t)}
\mathbh{1}_{\{w^{(N)}_j (t)\leq\chi^{(N)}(t), w^{(N)}_j (t)<r_j \}
}
\nonumber\\[-8pt]\\[-8pt]
&=&
\eta^{(N)}_t\bigl[0,\chi^{(N)}(t)\bigr]. \nonumber
\end{eqnarray}

\subsubsection{Description of service dynamics}
\label{subsub-service}

Analogous to the potential waiting process $w_j^{(N)}$,
the age process $a_j^{(N)}$ associated with customer $j$ is (for every
realization)
defined to be the piecewise linear function on $[0,\infty)$ that
equals $0$ till the customer enters service, then increases linearly
while the customer is in service (representing the amount of time
elapsed since entering service) and is then constant (equal to the
total service requirement)
after the customer completes service and departs the system.
For\vspace*{-2pt} $j=-\mathcal{E}^{(N)}_0+1, \ldots,0$, let $a^{(N)}_j(0)$
represent the age of the $j$th
customer in service at time 0 and for $j\in{\mathbb N}$, we set
$a^{(N)}_j(0)=0$.
Due to the First-Come-First-Serve (FCFS) nature of the system,
customers in service at
time $t$ are those that did not renege, that
have been in the system longer than
the head-of-the-line customer at time $t$, but
have not yet completed service and departed.
Therefore, a.s., for $j=-\mathcal{E}^{(N)}_0+1, \ldots,0,\ldots
,E^{(N)}(t)$, $t\geq0$,
%
%
\begin{equation}\label{adif}
\frac{d a^{(N)}_j(t+)}{dt}=\cases{
0, &\quad if $a^{(N)}_j(t)=0, w^{(N)}_j(t)=r
_j$, \vspace*{2pt}\cr
&\quad or $a^{(N)}_j(t)=0, w^{(N)}_j(t)\leq\chi^{(N)}(t)$, \vspace*{2pt}\cr
&\quad or $a^{(N)}_j(t) = v_j$, \vspace*{2pt}\cr
1, &\quad if $a^{(N)}_j(t)=0, \chi^{(N)}(t)< w^{(N)}_j(t) < r_j$, \vspace*{2pt}\cr
&\quad or $0< a^{(N)}_j(t) < v_j$.}
\end{equation}
Note that the condition in the penultimate line of the right-hand side
above represents the scenario in which a
customer enters service precisely at time $t$, which causes
$\chi^{(N)}$ to have a downward jump at time $t$ since the
condition that the arrival process increases only in unit jumps
ensures that there is at most one customer with a given potential
waiting time.

Now, for $t\in[0,\infty)$, let $\nu^{(N)}_t$ be the discrete nonnegative
Borel measure on $[0,H^s)$ that has a unit mass at the age of each of the
customers in service at time~$t$. Then, in a fashion analogous to
(\ref{def-etan}), the age measure $\nu^{(N)}_t$ can be
explicitly represented as
%
%
\begin{equation}
\label{def-nun}
\nu^{(N)}_t = \sum_{j = -\mathcal{E}^{(N)}_0+ 1}^{E^{(N)}(t)} \delta
_{a^{(N)}_j(t)}
\mathbh{1}_{ \{{da^{(N)}_j }/{dt}(t+) >0 \}}.
\end{equation}

\subsubsection{Auxiliary processes}
\label{subsub-aux}

We now introduce certain auxiliary processes that will be useful for
the study of the evolution
of the system.
\begin{itemize}
\item
The cumulative reneging process $R^{(N)}$, where $R^{(N)}(t)$ is the
cumulative number
of customers that have reneged from the system in the time interval $[0,t]$;
\item
the cumulative potential reneging process $S^{(N)}$, where $S^{(N)}(t)$
represents the cumulative number of customers whose potential waiting times
have reached their patience times in the interval $[0,t]$;
\item
the cumulative departure process $D^{(N)}$, where $D^{(N)}(t)$ is the cumulative
number of customers that have departed the system after completion of service
in the interval $[0,t]$;
\item
the process $K^{(N)}$, where $K^{(N)}(t)$ represents the cumulative
number of
customers that
have entered service in the interval $[0,t]$.
\end{itemize}
Now, a customer $j$ completes service (and therefore departs the
system) at
time $s$ if and only if, at time $s$, the left derivative of
$a_j^{(N)}$ is
positive and the right derivative of $a_j^{(N)}$ is zero. Therefore, we
can write
%
%
\begin{equation}
\label{def-depart} D^{(N)}(t) = \sum_{j = -\mathcal{E}^{(N)}_0+
1}^{E^{(N)}(t)}
\sum_{s\in[0,t]} \mathbh{1}_{ \{{da^{(N)}_j }/{dt}(s-) >0,
{da^{(N)}_j }/{dt}(s+)=0 \}}.
\end{equation}
Note that the second sum in (\ref{def-depart}) is well defined since
for each $t\geq0$ and each $j$ between $-\mathcal{E}^{(N)}_0+ 1$ and
$E^{(N)}(t)$, the
piecewise linear structure of\vspace*{1pt}
$a^{(N)}_j$ ensures that the indicator function in the sum is nonzero
for at most one $s\in[0,t]$, that is, there exists at most one $s\in
[0,t]$ such that the customer $j$ completes service at time $s$. A
similar logic shows that the cumulative potential reneging
process $S^{(N)}$ admits the representation
%
%
\begin{equation}
\label{def-cvrp} S^{(N)}(t) = \sum_{j = -\mathcal{E}^{(N)}_0+
1}^{E^{(N)}(t)} \sum
_{s\in[0,t]} \mathbh{1}_{ \{{dw^{(N)}_j}/{dt}(s-) >0,
{dw^{(N)}_j }/{dt}(s+)=0 \}},
\end{equation}
and the cumulative reneging process $R^{(N)}$ admits the representation
%
%
\begin{eqnarray}
\label{def-crp}
&&
R^{(N)}(t) \nonumber\\[-10pt]\\[-10pt]
&&\qquad= \sum_{j = -\mathcal{E}^{(N)}_0+
1}^{E^{(N)}(t)} \sum
_{s\in[0,t]} \mathbh{1}_{ \{w^{(N)}_j(s)\leq\chi^{(N)}(s-),
{dw^{(N)}_j}/{dt}(s-) >0, {dw^{(N)}_j }/{dt}(s+)=0 \}},\nonumber\hspace*{-28pt}
\end{eqnarray}
where the additional restriction $w^{(N)}_j(s)\leq\chi^{(N)}(s-)$ is imposed
so as to only count the reneging of customers actually in queue (and
not the
reneging of all
customers in the potential queue, which is captured by $S^{(N)}$).
Here, one considers the left limit $\chi^{(N)}(s-)$ of $\chi^{(N)}$
at time
$s$ to capture the situation in which $\chi^{(N)}$ jumps down at time $s$
due to the head-of-the-line customer reneging from the queue or
entering service.

Now, $\langle{\mathbf{1}}, \nu^{(N)}_t\rangle= \nu
^{(N)}_t[0,\infty)$ represents
the total number of customers in service at time $t$.
Therefore, mass balances on the total number of customers in the
system, the
number of customers waiting in the ``potential queue'' and the number
of customers in service show that
%
%
\begin{eqnarray}
\label{def-dn}
X^{(N)}(0) + E^{(N)}&=& X^{(N)}+ D^{(N)}+ R^{(N)},
\\
%
%
\label{def-sn}
\bigl\langle{\mathbf{1}}, \eta^{(N)}_0 \bigr\rangle+ E^{(N)}&=& \bigl\langle
{\mathbf{1}}, \eta^{(N)}
\bigr\rangle+ S^{(N)}
\end{eqnarray}
and
%
%
\begin{equation}
\label{def-kn}
\bigl\langle{\mathbf{1}}, \nu^{(N)}_0 \bigr\rangle+ K^{(N)}= \bigl\langle{\mathbf
{1}}, \nu^{(N)}
\bigr\rangle+ D^{(N)}.
\end{equation}
In addition, it is also clear that
%
%
\begin{equation}
\label{def-xn}
X^{(N)}= \bigl\langle{\mathbf{1}}, \nu^{(N)}\bigr\rangle+ Q^{(N)}.
\end{equation}
Combining (\ref{def-dn}), (\ref{def-kn}) and (\ref{def-xn}), we
obtain the following mass balance for the number of customers in queue:
%
%
\begin{equation}
\label{mass-queue}
Q^{(N)}(0) + E^{(N)}= Q^{(N)}+ R^{(N)}+ K^{(N)}.
\end{equation}
Furthermore, the nonidling condition takes the form
\[
N-\bigl\langle{\mathbf{1}}, \nu^{(N)}\bigr\rangle= \bigl[N - X^{(N)}\bigr]^+.
\]
Indeed, note that this ensures that when $X^{(N)}(t) < N$, the number in
the system is equal to the
number in service, and so there is no queue, while if $X^{(N)}(t) > N$,
there is a positive queue and
$\langle{\mathbf{1}}, \nu^{(N)}_t \rangle=N$, indicating that there
are no
idle servers.

An advantage of the measure-valued state representation that we adopt is
that it allows us to simultaneously
study several other functionals of interest. As an example,
we consider the so-called virtual waiting time process, which
is important for applications.
For each $t\geq0$, the virtual waiting time $W^{(N)}(t)$
is defined to be the amount of time a (virtual) customer with infinite
patience would have to wait before entering service if he were to
arrive at time $t$. For a more precise definition of $W^{(N)}$,
let $t\in[0,\infty)$ and
for each $s\in
[0,\infty)$, define
%
%
\begin{eqnarray} \label{dis:Tn}
\mathcal T_t^{(N)}(s)
&\doteq& \sum_{u \in[t,t+s]} \sum
_{j=-\mathcal{E}^{(N)}_0+
1}^{E^{(N)}(t)}
\mathbh{1}_{ \{{dw^{(N)}_j }/{dt}(u-) >0, {dw^{(N)}_j }/{dt}(u
+)=0 \}}\nonumber\\[-8pt]\\[-8pt]
&&\hspace*{79.28pt}{}\times
\mathbh{1}_{\{w^{(N)}_j (u)\leq\chi^{(N)}(u-)\}}. \nonumber
\end{eqnarray}
Observe that $\mathcal T_t^{(N)}(s)$ equals the cumulative number of
customers
who arrived before time $t$ and reneged from the queue (before entering
service) in the time interval $[t,t+s]$. Once again, for each $j$ there
is at most one $u\in[t,t+s]$ for which both indicator functions in the
summation are nonzero, and so
the sum is well defined. The virtual waiting time
$W^{(N)}(t)$ of a customer at time $t$ is the least amount of
time $s$ that elapses
after time $t$ such that the cumulative departure from the system of
customers that arrived prior to time $t$
strictly exceeds the queue length at time $t$.
Observing that this cumulative departure in the interval
$[t,t+s]$ can be
due to either departure from service or
reneging of customers that arrived prior to time $t$,
we can express the virtual waiting time as
%
%
\begin{equation} \label{T}\qquad
W^{(N)}(t)
\doteq \inf\bigl\{s\geq0\dvtx D^{(N)}(t+s)- D^{(N)}(t)+ {\mathcal
T}_t^{(N)}(s)> Q^{(N)}(t)\bigr\}.
\end{equation}
Here, we have used the fact that for all $s$ such that $D^{(N)}(t+s)- D^{(N)}
(t)+ {\mathcal T}_t^{(N)}(s)\leq Q^{(N)}(t)$, every customer that departed
in the time interval $[t,t+s]$ must have arrived prior to time $t$.

\subsubsection{Filtration}
\label{subsub-filt}
The total number of customers
in service at time $t$ is given by $\langle{\mathbf{1}}, \nu^{(N)}_t
\rangle
=\nu^{(N)}_t [0,H^s)$
and is bounded above by $N$. In addition, from (\ref{def-sn}) it
follows that
\[
\bigl\langle{\mathbf{1}}, \eta^{(N)}_t \bigr\rangle= \eta^{(N)}_t[0,H^r)
\leq E^{(N)}(t) + \bigl\langle{\mathbf{1}}, \eta^{(N)}_0 \bigr\rangle\leq
E^{(N)}(t) + \mathcal{E}_0^{(N)},
\]
which is a.s. finite by assumption. Therefore,
for every $t\in[0,\infty)$, a.s., $\nu^{(N)}_t \in\mathcal{M}_F
[0,H^s)$ and
$\eta^{(N)}_t \in\mathcal{M}_F [0,H^r)$.
Hence, the state descriptor
$(\alpha_E^{(N)},X^{(N)},\break\nu^{(N)},\eta^{(N)})$ takes values in
${\mathbb R}_+ \times{\mathbb Z}_+
\times\mathcal{M}_F
[0,H^s) \times\mathcal{M}_F [0,H^r)$.
For purely technical purposes we will find it
convenient to also introduce the additional ``station process''
$s^{(N)} \doteq(s_j^{(N)}, j \in{\mathbb Z})$, defined on the same
probability space $(\Omega,\mathcal F,\mathbb{P})$. For each $t \in
[0,\infty)$,
if customer $j$ has already
entered service by time $t$, then
$s_j^{(N)} (t)$ is equal to the index $i \in\{1, \ldots, N\}$ of the
station at which
customer $j$ receives/received service and $s_j^{(N)} (t) \doteq
0$ otherwise.
For $t\in[0,\infty)$, let $\tilde{\mathcal{F}}_t^{(N)}$ be the
$\sigma$-algebra generated by
\begin{eqnarray*}
&&\bigl\{\mathcal{E}^{(N)}_0,X^{(N)}(0),\alpha_E^{(N)}(s), w^{(N)}_j(s),
a^{(N)}_j(s), s_j^{(N)},\\
&&\hspace*{45pt} j \in\{-\mathcal{E}^{(N)}_0+1, \ldots, 0\}
\cup{\mathbb N}, s \in
[0,t] \bigr\},
\end{eqnarray*}
and let $\{\mathcal{F}_t^{(N)}\}$ denote the associated
right-continuous filtration, completed with respect to $\mathbb{P}$. In
Appendix \ref{ap-markov}, an explicit construction of
the state descriptor and auxiliary
processes is provided, which shows in particular
that the state descriptor $(\alpha_E^{(N)},X^{(N)},\nu^{(N)},\eta
^{(N)})$ and
auxiliary processes are c\`{a}dl\`{a}g.
Moreover, in Lemma~\ref{app:adapt}, it is proved that the state process
$V^{(N)}\doteq(\alpha_E^{(N)},X^{(N)},\nu^{(N)},\eta^{(N)})$ and
the processes
$E^{(N)}$, $Q^{(N)}$, $S^{(N)}$, $R^{(N)}$, $D^{(N)}$ and $K^{(N)}$ are
all $\mathcal
{F}_t^{(N)}$-adapted, and in Lemma \ref{lem:Mark}, it is shown that
$(V^{(N)}, \mathcal{F}_t^{(N)})$ is a strong Markov process.

\subsection{A succinct characterization of the dynamics}
\label{subs-prelimit}
The main result of this section is Theorem \ref{th-prelimit}, which provides
equations that more succinctly characterize the dynamics of the state
$(\alpha_E^{(N)},X^{(N)},\nu^{(N)},\eta^{(N)})$ described in Section
\ref{sec:repdyn}.
First, we introduce some notation that is required to state the result.

For any measurable function $\varphi$ on $[0,H^s)\times{\mathbb
R}_+$, consider
the process $D^{(N)}_\varphi$ that takes
values in ${\mathbb R}$, and is given by
%
%
\begin{equation}
\label{def-baren}
D^{(N)}_\varphi(t) \doteq\sum_{s \in[0,t]} \sum_{j=-\mathcal
{E}^{(N)}_0+ 1}^{E^{(N)}(t)}
\mathbh{1}_{ \{{da^{(N)}_j }/{dt}(s-) >0, {da^{(N)}_j
}/{dt}(s+)=0 \}}
\varphi\bigl(a^{(N)}_j (s),s\bigr)\hspace*{-32pt}
\end{equation}
for $t \in[0,\infty)$.
It follows immediately from (\ref{def-baren}) and the right continuity
of the filtration $\{\mathcal{F}_t^{(N)}\}$ that
$D^{(N)}_\varphi$ is $\{\mathcal{F}_t^{(N)}\}$-adapted.
Also, comparing (\ref{def-baren}) with (\ref{def-depart}), it is
clear that when $\varphi$ is the constant function ${\mathbf{1}}$,
$D^{(N)}_{\mathbf{1}}$
is exactly the cumulative departure process $D^{(N)}$, that is,
%
%
\begin{equation}D^{(N)}_{\mathbf{1}}
=D^{(N)}. \label{equivD}
\end{equation}
%
In an exactly analogous fashion, for any measurable function $\psi$ on
$[0,H^r)\times{\mathbb R}_+$, consider the
process $S^{(N)}_\psi$ that takes
values in ${\mathbb R}$, and is given by
%
%
\begin{equation}
\label{def-barwn}
S^{(N)}_\psi(t) \doteq\sum_{s \in[0,t]} \sum_{j=-\mathcal
{E}^{(N)}_0+ 1}^{E^{(N)}(t)}
\mathbh{1}_{ \{{dw^{(N)}_j }/{dt}(s-) >0, {dw^{(N)}_j
}/{dt}(s+)=0 \}}
\psi\bigl(w^{(N)}_j (s),s\bigr).\hspace*{-37pt}
\end{equation}
It follows immediately from (\ref{def-barwn}) and the right continuity
of the filtration $\{\mathcal{F}_t^{(N)}\}$ that for $t\in[0,\infty)$,
$S^{(N)}_\psi$ is $\{\mathcal{F}_t^{(N)}\}$-adapted.
Moreover, $S^{(N)}_{\mathbf{1}}$ is clearly equal to the cumulative potential
reneging process $S^{(N)}$, that is,
%
%
\begin{equation}S^{(N)}_{\mathbf{1}}=S^{(N)}. \label{equivS}
\end{equation}
In addition, using
(\ref{def-dn}), (\ref{def-xn}) and the nonnegativity of $Q^{(N)}$, $R^{(N)}$
and $\langle{\mathbf{1}}, \nu^{(N)}\rangle$,
it follows that for any $t \in[0,\infty)$ and bounded,
measurable $\varphi$,
%
%
\begin{equation}
\label{bd-1}
\mathbb{E}\bigl[ \bigl|D^{(N)}_\varphi(t) \bigr| \bigr]
\leq\Vert\varphi\Vert_\infty\mathbb{E}\bigl[X^{(N)}(0) + E^{(N)}(t)
\bigr]<\infty
\end{equation}
and likewise, for each $t \in
[0,\infty)$ and bounded measurable $\psi$, (\ref{def-sn}) shows that
%
%
\begin{equation}
\label{bd-2}
\mathbb{E}\bigl[ \bigl|S^{(N)}_\psi(t) \bigr| \bigr] \leq
\Vert\psi\Vert_\infty\mathbb{E}\bigl[\bigl\langle{\mathbf{1}}, \eta
^{(N)}_0\bigr\rangle+
E^{(N)}
(t) \bigr] <
\infty.
\end{equation}

Next, comparing (\ref{def-crp}) with (\ref{def-barwn}), it is clear
that the cumulative reneging process $R^{(N)}$ satisfies
%
%
\begin{equation} \label{RQ}
R^{(N)}(t)=S^{(N)}_{\theta^{(N)}}(t),\qquad t\geq0,
\end{equation}
where $\theta^{(N)}$ is given by
%
%
\begin{equation}\label{ps}\theta^{(N)}(x,s)=\mathbh{1}_{[x,\infty
)}\bigl(\chi
^{(N)}(s-)\bigr),\qquad
x\in{\mathbb R}, s\geq0.
\end{equation}

We now state the main result of this section.
For $s, r \in[0,\infty)$, recall that
$\langle\varphi(\cdot+ r,s), \nu^{(N)}_s \rangle$ is used as a
short form for
$\int_{[0, H^s)} \varphi(x+ r,s) \nu^{(N)}_s (dx)$, and likewise
for~$\eta^{(N)}$.
\begin{theorem}
\label{th-prelimit}
The processes $(E^{(N)}, X^{(N)},
\nu^{(N)},\eta^{(N)})$ a.s. satisfy the following coupled set of
equations: for
$\varphi\in\mathcal{C}^1_c([0,H^s)\times{\mathbb R}_+)$ and $t \in
[0,\infty)$,
%
%
\begin{eqnarray}
\label{eqn-prelimit1}
\bigl\langle\varphi(\cdot,t), \nu^{(N)}_{t} \bigr\rangle
& = & \bigl\langle\varphi(\cdot, 0), \nu^{(N)}_{0} \bigr\rangle+ \int
_{0}^t \bigl\langle\varphi_x(\cdot,s) + \varphi_s(\cdot,s), \nu^{(N)}_s
\bigr\rangle \,ds \nonumber\\[-8pt]\\[-8pt]
& &{} - D^{(N)}_\varphi(t) + \int_{[0,t]} \varphi(0,s)\, dK^{(N)}(s)
,\nonumber
\end{eqnarray}
for
$\psi\in\mathcal{C}^1_c([0,H^r)\times{\mathbb R}_+)$ and $t \in
[0,\infty)$,
%
%
\begin{eqnarray}\quad
\label{eqn-prelimit3} \bigl\langle\psi(\cdot, t), \eta^{(N)}_{t}
\bigr\rangle
& = & \bigl\langle\psi(\cdot, 0), \eta^{(N)}_{0} \bigr\rangle+ \int
_{0}^t \bigl\langle\psi_x(\cdot,s) + \psi_s(\cdot,s), \eta^{(N)}_s
\bigr\rangle \,ds \nonumber\\[-8pt]\\[-8pt]
& &{} - S^{(N)}_\psi(t) + \int_{[0,t]} \psi(0,s) \,dE^{(N)}(s)
,\nonumber\\
\label{eqn-prelimit2}
X^{(N)}(t) & = & X^{(N)}(0) + E^{(N)}(t) - D^{(N)}_{\mathbf{1}}(t)- R^{(N)}
(t), \\
\label{comp-prelimit}
N - \bigl\langle{\mathbf{1}}, \nu^{(N)}_t \bigr\rangle &=& \bigl[N - X^{(N)}(t) \bigr]^+,
\end{eqnarray}
where $K^{(N)}$ satisfies (\ref{def-kn}), $R^{(N)}$ satisfies
(\ref{RQ}) and $D^{(N)}_\varphi$ and $S^{(N)}_\psi$ are the processes
defined in (\ref{def-baren}) and (\ref{def-barwn}), respectively.
\end{theorem}
\begin{remark}
\label{rem-compdyn}
In the service dynamics, customer arrivals into
service are governed by the process $K^{(N)}$, the random duration in
service is determined by the distribution $G^s$ and departures are
represented by $D^{(N)}$. As captured by
(\ref{eqn-prelimit1}) and (\ref{eqn-prelimit3}), the dynamics of the
\textit{potential queue} is exactly analogous, with the customer
arrivals now governed by the process $E^{(N)}$, the random duration of stay
in the potential queue determined by $G^r$ and potential departures due
to reneging represented by $S^{(N)}$. Moreover, \textit{given
$K^{(N)}$}, the
dynamics of $\nu^{(N)}$ are exactly the same as in the case without
abandonment, which was well studied in~\cite{kasram07}. However, in the
presence of reneging, there is a significantly more complicated
coupling between $\nu^{(N)}$ and $K^{(N)}$, as captured by
(\ref{eqn-prelimit2}) and~(\ref{comp-prelimit}). In particular, this
involves the cumulative reneging process $R^{(N)}$, which has no analogy
with any quantity in the system without abandonments. Instead, as shown
in the sequel [see Lemma \ref{cor:1}, (\ref{rep-rcomp2}) and
Proposition \ref{prop:3}], we will exploit the representation
(\ref{RQ}) of $R^{(N)}$ in terms of the ``known'' quantity $S^{(N)}$
in order
to characterize the limit of the scaled sequence of reneging processes.
\end{remark}
\begin{pf*}{Proof of Theorem \ref{th-prelimit}}
The proof of (\ref{eqn-prelimit1}) can be carried out in exactly the
same way as the proof of (5.2) in Theorem 5.1 of \cite{kasram07}, since
the definition of $\nu^{(N)}$ in \cite{kasram07} is equivalent to the
definition given in (\ref{def-nun}) here since $da^{(N)}_j(t+)/dt = 0$
for all $j > K^{(N)}(t)$ in \cite{kasram07}. For the reasons mentioned in
Remark \ref{rem-compdyn}, the proof of (\ref{eqn-prelimit3}) is also
analogous except that the condition that each $\nu^{(N)}_t$ has total mass
no greater than $N$ is replaced by the argument below, which shows that
each $\eta^{(N)}_t$ has finite mass. We know that for $k = 0, \ldots,
\lfloor nt \rfloor$,
\[
\bigl\langle{\mathbf{1}},\eta^{(N)}_{({k+1})/{n}}\bigr\rangle \leq
E^{(N)}\biggl(\frac{k+1}{n} \biggr)+ \bigl\langle{\mathbf{1}},\eta^{(N)}_0
\bigr\rangle\leq
E^{(N)}(t+1)+\bigl\langle
{\mathbf{1}}, \eta^{(N)}_0\bigr\rangle.
\]
Thus, by taking the supremum over $k = 0, \ldots, \lfloor nt \rfloor$,
we have a.s.
%
%
\begin{equation} \label{renegnf}
\sup_{k = 0, \ldots, \lfloor nt
\rfloor} \bigl\langle{\mathbf{1}},\eta^{(N)}_{({k+1})/{n}} \bigr\rangle\leq
E^{(N)}(t+1)+\mathcal{E}^{(N)}_0<\infty.
\end{equation}
Equation (\ref{eqn-prelimit2})
follows from (\ref{def-dn}), (\ref{equivD}) and (\ref{RQ}), while
(\ref{comp-prelimit}) is just the nonidling condition
formulated in Section \ref{subsub-aux}.
\end{pf*}

\section{Main results}
\label{sec-flconj}
In this section we summarize our main results.
First, in Section \ref{subs-scaling}, we introduce the
fluid-scaled quantities and state our basic assumptions.
Then, in Section \ref{subs-fleqs}, we introduce the so-called
fluid equations, which provide a continuous analog of the characterization
of the discrete model given in
Theorem~\ref{th-prelimit}.
In Section \ref{subs-mainres} we present our main theorems. In
particular, we show that the fluid equations uniquely characterize the
strong law of large numbers or ``fluid''
limit of the many-server system, as the number of
servers goes to infinity.

\subsection{Fluid scaling and basic assumptions}
\label{subs-scaling}

Consider the following scaled versions of the basic processes
described in Section \ref{sec-mode}. For each $N \in{\mathbb N}$, the scaled
version of the state descriptor
$(\overline{\alpha}{}^{(N)}_E,\overline{X}{}^{(N)},\overline{\nu
}^{(N)},\overline{\eta}^{(N)})$ is given by
%
%
\begin{eqnarray} \label{fl-scaling}
\overline{\alpha}_E^{(N)}(t)& \doteq&\alpha_E^{(N)}(t),\qquad \overline
{X}{}^{(N)}(t)
\doteq\frac{X^{(N)}(t)}{N}, \\ \overline{\nu}{}^{(N)}_t (B) &\doteq&
\frac{\nu^{(N)}_t (B)}{N}, \qquad\overline{\eta}^{(N)}_t (B) \doteq
\frac{\eta^{(N)}_t (B)}{N},
\end{eqnarray}
for $t \in[0,\infty)$ and any Borel subset $B$ of ${\mathbb R}_+$.
Analogously, define
%
%
\begin{equation}
\label{fl-scaling2} 
\overline{I}^{(N)}\doteq\frac{I^{(N)}}{N} \qquad\mbox{for } I=E,D,K,Q,R,S,
\mathcal{T}_t.
\end{equation}

Recall that $\mathcal{I}_{{\mathbb R}_+}[0,\infty)$ is the subset of
nondecreasing functions
$f \in\mathcal{D}_{{\mathbb R}_+}[0,\infty)$ with $f(0) = 0$,
$H^s=\sup\{x \in[0,\infty)\dvtx g^s(x) > 0\}$ and $H^r=\sup\{x \in
[0,\infty)\dvtx\break g^r(x) > 0\}$. Define
%
%
\begin{equation}\qquad
\label{def-newspace}
\mathcal{S}_0\doteq\left\{\matrix{
(e,x,\nu, \eta) \in\mathcal{I}_{{\mathbb R}_+}[0,\infty)\times
{\mathbb R}_+ \times
\mathcal{M}_F[0,H^s)\times\mathcal{M}_F[0,H^r)\mbox{:}\cr
1 - \langle{\mathbf{1}}, \nu\rangle= [1-x]^+}\right\}.
\end{equation}
$\mathcal{S}_0$ serves as the space of possible input data for the
fluid equations.
Our goal is to identify
the limit in distribution of the quantities
$(\overline{X}{}^{(N)}, \overline{\nu}{}^{(N)},\overline{\eta
}^{(N)})$, as \mbox{$N \rightarrow\infty$}. To this end,
we impose some natural assumptions on the sequence of
initial conditions $(\overline{E}{}^{(N)}, \overline{X}{}^{(N)}(0),
\overline{\nu}{}^{(N)}_0,\overline{\eta}{}^{(N)}_0)$.
\begin{ass}[(Initial conditions)]
\label{ass-init}
There exists an $\mathcal{S}_0$-valued random variable $(\overline
{E}, \overline{X}(0),
\overline{\nu}_0,\overline{\eta}_0)$
such that,
as $N \rightarrow\infty$, the following limits hold:
\begin{enumerate}
\item
$\overline{E}{}^{(N)}\rightarrow\overline{E}$ in $\mathcal
{D}_{{\mathbb R}_+}[0,\infty)$ $\mathbb{P}$-a.s., and
$\mathbb{E}[\overline{E}{}^{(N)}(t) ] \rightarrow\mathbb{E}[\overline
{E}(t) ]<\infty$ for
every $t \in[0,\infty)$;
%
\item
$\overline{X}{}^{(N)}(0) \rightarrow\overline{X}(0)$ in ${\mathbb
R}_+$ $\mathbb{P}$-a.s.;
\item
$\overline{\nu}{}^{(N)}_0 \stackrel{w}{\rightarrow}\overline{\nu
}_0$ in $\mathcal{M}_F[0,H^s)$;
\item
$\overline{\eta}^{(N)}_0 \stackrel{w}{\rightarrow}\overline{\eta
}_0$ in $\mathcal{M}_F[0,H^r)$, and $\mathbb{E}[\langle
1, \overline{\eta}^{(N)}_0 \rangle] \rightarrow\mathbb{E}[\langle
1, \overline{\eta}_0 \rangle] < \infty$.
\end{enumerate}
\end{ass}
\begin{remark}
\label{rem-skorep}
If the limits in (1) and (2) of Assumption \ref{ass-init} hold only in
distribution rather than
almost surely, then using the Skorokhod representation theorem in the
standard way,
it can be shown that all the stochastic process convergence results in the
paper continue to hold. Also, (1) and (4) of Assumption \ref{ass-init} and
(\ref{comp-prelimit}) imply that, for every $T\in[0,\infty)$,
%
%
\begin{equation}\label{init-bd}\qquad
\sup_{t\in[0,T]}\sup_N\mathbb
{E}\bigl[\overline{X}{}^{(N)}
(0)+\overline{E}{}^{(N)}
(t) \bigr]\leq\mathbb{E}\bigl[1+\bigl\langle1, \overline{\eta}^{(N)}_0 \bigr\rangle
+\overline{E}{}^{(N)}(T)
\bigr]<\infty.
\end{equation}
\end{remark}

The next assumption imposes some regularity conditions on $\overline
{\eta}_0$ and
$\overline{E}$.
\begin{ass} \label{ass-jump}
For each $t\geq0$, if $\overline{\eta}_0(\{t\})>0$, then $\overline
{\eta}
_0(t,t+\varepsilon)>0$ for every $\varepsilon>0$ and if $\overline
{E}(t)-\overline{E}
(t-)>0$, then $\overline{E}(t-)-\overline{E}(t-\varepsilon)>0$ for
every $\varepsilon>0$.
\end{ass}
\begin{remark} Assumption \ref{ass-jump} is trivially satisfied if
$\overline{\eta}_0$ and $\overline{E}$ are continuous, that is,
$\overline{\eta}_0(\{t\})=0$
for all $t\geq0$ and the function $\overline{E}$ is continuous.
\end{remark}

In order to state our last assumption,
define the hazard rate functions of $G^r$ and $G^s$ in the usual manner
%
%
\begin{eqnarray}
\label{def-h}
h^r(x) &\doteq&\frac{g^r(x)}{1 - G^r(x)},\qquad x\in[0,H^r),\\
h^s(x) &\doteq&\frac{g^s(x)}{1 - G^s(x)},\qquad x \in[0,H^s).
\end{eqnarray}
It is easy to verify that $h^r$ and $h^s$ are locally integrable on
$[0,H^r)$ and $[0,H^s)$, respectively.
\begin{ass} \label{ass-h}
There exists $L^s<H^s$ such that $h^s$ is either bounded or
lower-semicontinuous on $(L^s,H^s)$, and, likewise, there exists
$L^r<H^r$ such that $h^r$ is either bounded or lower-semicontinuous on
$(L^r,H^r)$.
\end{ass}

\subsection{Fluid equations}
\label{subs-fleqs}

We now introduce the so-called fluid
equations and provide some intuition as to why the
limit of any sequence $(\overline{X}{}^{(N)},\overline{\nu
}^{(N)},\overline{\eta}^{(N)})$ should be expected
to be a solution to these equations.
In Section \ref{sec:CSL}, we provide a rigorous proof of this fact.
\begin{defn}[(Fluid equations)]
\label{def-fleqns}
The c\`{a}dl\`{a}g function $(\overline{X},
\overline{\nu},\overline{\eta})$ defined on $[0,\infty)$ and
taking values in ${\mathbb R}_+
\times
\mathcal{M}_F[0,H^s)\times\mathcal{M}_F[0,H^r)$ is said to solve the
\textit{fluid equations} associated with
$(\overline{E}, \overline{X}(0), \overline{\nu}_0,\overline{\eta
}_0) \in\mathcal{S}_0$ and the hazard
rate functions $h^r$ and $h^s$
if and only if
for every $t \in[0,\infty)$,
%
%
\begin{equation}
\label{cond-radon}
\int_0^t \langle h^r, \overline{\eta}_s \rangle \,ds < \infty,\qquad \int_0^t
\langle h^s, \overline{\nu}_s \rangle \,ds < \infty,
\end{equation}
and the following relations are satisfied: for every $\varphi\in
\mathcal{C}^1_c([0,H^s)\times{\mathbb R}_+)$,
%
%
\begin{eqnarray}
\label{eq-ftmeas}\quad
\langle\varphi(\cdot,t), \overline{\nu}_t \rangle& = & \langle
\varphi(\cdot,0),
\overline{\nu}_0 \rangle+
\int_0^t \langle\varphi_s(\cdot,s), \overline{\nu}_s \rangle \,ds +
\int_0^t \langle\varphi_x(\cdot,s), \overline{\nu}_s \rangle \,ds
\nonumber\\[-8pt]\\[-8pt]
& &{}
- \int_0^t \langle h^s(\cdot) \varphi(\cdot,s), \overline{\nu}_s
\rangle \,ds+
\int_0^t \varphi(0,s) \,d\overline{K}(s),\nonumber
\end{eqnarray}
where
%
%
\begin{equation}
\label{eq-fk}
\overline{K}(t) = \langle{\mathbf{1}}, \overline{\nu}_t \rangle-
\langle{\mathbf{1}},
\overline{\nu}_0 \rangle+
\int_0^t \langle h^s, \overline{\nu}_s \rangle \,ds;
\end{equation}
for every $\psi\in\mathcal{C}^1_c([0,H^r)\times{\mathbb R}_+)$
%
%
\begin{eqnarray}
\label{eq-ftreneg}\qquad\quad
\langle\psi(\cdot,t), \overline{\eta}_t \rangle& = & \langle\psi
(\cdot
,0), \overline{\eta}_0 \rangle+
\int_0^t \langle\psi_s(\cdot,s), \overline{\eta}_s \rangle \,ds +
\int_0^t \langle\psi_x(\cdot,s), \overline{\eta}_s \rangle \,ds
\nonumber\\[-8pt]\\[-8pt]
& &{}
- \int_0^t \langle h^r(\cdot) \psi(\cdot,s), \overline{\eta}_s
\rangle \,ds+ \int_0^t \psi(0,s)\, d\overline{E}(s);\nonumber
%
%
\\
\label{fq}
\overline{Q}(t)&=& \overline{X}(t) - \langle{\mathbf{1}}, \overline
{\nu}_t \rangle; \\
\label{fqfreneg} \overline{Q}(t)&\leq&\langle{\mathbf{1}},
\overline{\eta}_t \rangle; \\
\label{fr}
\overline{R}(t) &=& \int_0^t \biggl(\int_0^{\overline
{Q}(s)}h^r((F^{\overline\eta_s}
)^{-1}(y))\,dy \biggr) \,ds,
\end{eqnarray}
where we recall that $F^{\overline\eta_t}(x)=\overline{\eta}_t[0,x];$
%
%
\begin{equation}
\label{eq-fx}
\overline{X}(t) = \overline{X}(0) + \overline{E}(t) - \int_0^t
\langle h^s, \overline{\nu}_s
\rangle \,ds - \overline{R}(t)
\end{equation}
and
%
%
\begin{equation}
\label{eq-fnonidling}
1 - \langle{\mathbf{1}}, \overline{\nu}_t \rangle= [1 - \overline
{X}(t)]^+.
\end{equation}
\end{defn}

It immediately follows from (\ref{fq}) and (\ref{eq-fnonidling}) that
for each $t\in[0,\infty)$,
%
%
\begin{equation}
\label{fqfx}\overline{Q}(t)=[\overline{X}(t)-1]^+.
\end{equation}
For future use, we also observe that (\ref{eq-fk}), (\ref{fq})
and (\ref{eq-fx}), when combined, show that for every $t\in[0,\infty)$,
%
%
\begin{equation}\label{qt-conserve} \overline{Q}(0)+\overline
{E}(t)=\overline{Q}(t)+\overline{K}
(t)+\overline{R}(t).
\end{equation}

We now provide an informal, intuitive explanation for the form of the fluid
equations. Equations (\ref{eq-fk}), (\ref{fq}) and (\ref{eq-fx}) are simply
mass conservation equations, that are fluid analogs of (\ref{def-kn}),
(\ref{def-xn}) and (\ref{eqn-prelimit2}), respectively, while (\ref
{fqfreneg})
expresses a bound, whose analog clearly holds in the pre-limit, as can
be seen
from (\ref{qn}). The relation (\ref{eq-fnonidling}) is simply the
fluid analog
of the nonidling condition (\ref{comp-prelimit}). Equations (\ref{eq-ftmeas})
and (\ref{eq-ftreneg}), which govern the evolution of the fluid age measure
$\overline{\nu}$ and queue measure $\overline{\eta}$, respectively,
are natural
analogs of the
pre-limit equations (\ref{eqn-prelimit1}) and (\ref{eqn-prelimit3}),
respectively. It is worthwhile to comment further on~the fourth terms
on the
right-hand sides of (\ref{eq-ftmeas}) and (\ref{eq-ftreneg}), which
characterize the fluid departure rate and potential reneging rate,
respectively, as integrals
of the corresponding hazard rate with respect to the age and queue
measures. Note that $\overline{\nu}_s(dx)$ represents the amount of
mass (limiting
fraction of customers) whose age lies in the range $[x,x+dx)$ at time
$s$, and
$h^s(x)$ represents the fraction of mass with age $x$ (i.e., with
service time
no less than $x$) that would depart from the system while having age in
$[x,x+dx)$. Hence, it is natural to expect $\langle h^s,\overline{\nu
}_s \rangle$
to represent the
departure rate of mass from the fluid system at time $s$. This was rigorously
proved in the case without abandonments in \cite{kasram07} (see Proposition
5.17 therein). By exploiting the exact analogy between $(\overline{\nu
},\overline{K},
\overline{D})$ and
$(\overline{\eta},\overline{E}, \overline{S})$ (see Remark \ref
{rem-compdyn}), it is clear
that the
potential reneging rate at time $s$ can be similarly represented as
$\langle
h^r,\overline{\eta}_s \rangle$. Thus the fluid potential reneging
process $\overline{S}
$, defined by
%
%
\begin{equation}
\label{def-fs}
\overline{S}(t) \doteq\int_0^t \langle h^r, \overline{\eta}_s
\rangle \,ds,\qquad t
\in[0,\infty),
\end{equation}
represents the cumulative amount of potential reneging from the fluid system
in the interval $[0,t]$. Due to the FCFS nature of the system, the
fluid queue at
time $s$ contains all the mass in $\overline{\eta}$ that is to the
left of
$(F^{\overline\eta_s})^{-1}(\overline{Q}(s))$, where recall
$F^{\overline\eta_s}$ is the
c.d.f. of $\overline{\eta}_s$.
Moreover, roughly speaking, given any $y\in[0,\overline{Q}(s)]$,
there is a
mass of
$dy$ customers in the queue whose waiting time at $s$ is
$(F^{\overline\eta_s})^{-1}(y)$ and the mean reneging rate of
customers with this
waiting time is $h^r((F^{\overline\eta_s})^{-1}(y))$. Thus the total
actual reneging
that has occurred in the interval $[0,t]$, is represented by the
integral, as specified in (\ref{fr}).

We close the section with a simple result on the action of time-shifts
on solutions to the fluid equations. For this, we need the following
notation: for any $t\in[0,\infty)$,
\begin{eqnarray*}
\overline{E}^{[t]}&\doteq&\overline{E}(t+\cdot)-\overline{E}(t),\qquad
\overline{K}{}^{[t]}\doteq\overline{K}
(t+\cdot)-\overline{K}(t),\\
\overline{X}{}^{[t]}&\doteq&\overline{X}(t+\cdot),\qquad
\overline{\nu}^{[t]}\doteq\overline{\nu}_{t+\cdot},
\\
\overline{R}{}^{[t]}&\doteq&\overline{R}(t+\cdot)-\overline{R}(t),\qquad
\overline{\eta}^{[t]}\doteq
\overline{\eta}_{t+\cdot},\qquad \overline{Q}{}^{[t]}\doteq\overline
{Q}(t+\cdot).
\end{eqnarray*}
\begin{lemma} \label{lem:shift}
Suppose the c\`{a}dl\`{a}g function $(\overline{X},
\overline{\nu},\overline{\eta})$ defined on $[0,\infty)$ and
taking values in ${\mathbb R}_+
\times\mathcal{M}_F[0,H^s)\times\mathcal{M}_F[0,H^r)$ solves the
fluid equations associated with
$(\overline{E}, \overline{X}(0), \overline{\nu}_0,\overline{\eta
}_0) \in\mathcal{S}_0$, then $(\overline{X}{}^{[t]},
\overline{\nu}^{[t]},\overline{\eta}^{[t]})$ solves the fluid
equations associated with
$(\overline{E}^{[t]}, \overline{X}(t), \overline{\nu}_t,\overline
{\eta}_t) \in\mathcal{S}_0$, where
$\overline{K}{}^{[t]}, \overline{R}{}^{[t]}, \overline{Q}{}^{[t]}$ are the
corresponding processes
that satisfy (\ref{eq-fk}), (\ref{fr}), (\ref{fq}) with $\overline
{\nu}
^{[t]}$, $\overline{\eta}^{[t]}$ and $\overline{X}{}^{[t]}$ in place
of $\overline{\nu}$,
$\overline{\eta}$ and $\overline{X}$.
\end{lemma}

The proof of the lemma just involves a rewriting of the fluid
equations, and is thus omitted.

\subsection{Summary of main results}
\label{subs-mainres}

Our first result establishes uniqueness of solutions to the fluid equations.
\begin{theorem} \label{thm:1}
Given any $(\overline{E}, \overline{X}(0), \overline{\nu}_0,
\overline{\eta}_0) \in\mathcal{S}_0$, there
exists at
most one solution $(\overline{X},\overline{\nu},\overline{\eta})$
to the associated fluid equations
(\ref{cond-radon})--(\ref{eq-fnonidling}). Moreover, if $\overline
{\nu}$ and
$\overline{\eta}$
satisfy (\ref{cond-radon}), then $(\overline{X},\overline{\nu
},\overline{\eta})$ is a
solution to the
fluid equations if and only if for every $f \in\mathcal
{C}_b({\mathbb R}_+)$,
%
%
\begin{eqnarray} \label{eq-freneg2}
\int_{[0,H^r)} f (x) \overline{\eta}_t (dx) & = & \int_{[0,H^r)} f(x+t)
\frac{1 - G^r(x+t)}{1 - G^r(x)} \overline{\eta}_0 (dx)
\nonumber\\[-8pt]\\[-8pt]
& &{}
+ \int_{[0,t]} f(t-s) \bigl(1
- G^r(t-s)\bigr)\, d \overline{E}(s), \nonumber\\
\label{eq-fmeas2}
\int_{[0,H^s)} f (x) \overline{\nu}_t (dx) & = & \int_{[0,H^s)} f(x+t)
\frac{1 - G^s(x+t)}{1 - G^s(x)} \overline{\nu}_0 (dx) \nonumber\\[-8pt]\\[-8pt]
& &{}
+ \int_{[0,t]} f(t-s) \bigl(1
- G^s(t-s)\bigr)\, d \overline{K}(s), \nonumber
\end{eqnarray}
where
%
%
\begin{eqnarray}\label{dis:fk}
\overline{K}(t) &=& [ \overline
{X}(0)-1]^+ - [\overline{X}(t) -
1]^+ + \overline{E}
(t)\nonumber\\[-8pt]\\[-8pt]
&&{} - \int_0^t \biggl( \int_0^{[\overline{X}(s) - 1]^+} h^r ( (
F^{\overline{\eta}_s} )^{-1} (y)
) \,dy \biggr) \,ds\nonumber
\end{eqnarray}
and for all $t\in[0,\infty)$, $\overline{X}$ satisfies $[\overline
{X}(t)-1]^+\leq
\langle{\mathbf{1}},\overline{\eta}_t \rangle$, the nonidling
condition (\ref
{eq-fnonidling}) and
%
%
\begin{eqnarray}\label{dis:fx}
\overline{X}(t) &=& \overline{X}(0) +
\overline{E}(t) - \int_0^t
\langle h^s,
\overline{\nu}_s \rangle \,ds\nonumber\\[-8pt]\\[-8pt]
&&{} -
\int_0^t \biggl( \int_0^{[\overline{X}(s) - 1]^+} h^r ( (
F^{\overline{\eta}_s} )^{-1} (y)
) \,dy \biggr) \,ds.\nonumber
\end{eqnarray}
Moreover, $\overline{K}$ also satisfies
%
%
\begin{eqnarray}
\label{eq-fk2}\qquad
\overline{K}(t) &=& \langle{\mathbf{1}}, \overline{\nu}_{t-s}
\rangle- \langle{\mathbf{1}},
\overline{\nu}_0
\rangle+\int_{[0,H^s )}
\frac{G^s(x+t-s)- G^s(x)}{1 - G^s(x)} \overline{\nu}_0 (dx)
\nonumber\\
&&{} + \int_0^t \biggl( \langle{\mathbf{1}}, \overline{\nu}_{t-s} \rangle-
\langle
\mathbf{1}, \overline{\nu}_0 \rangle\\
&&\hspace*{31.4pt}{} +\int_{[0,H^s )}
\frac{G^s(x+t-s)- G^s(x)}{1 - G^s(x)} \overline{\nu}_0 (dx) \biggr)
u^s(s) \,ds,\nonumber
\end{eqnarray}
where $u^s$ is the density of the renewal function $U^s$ associated with
$G^s$
($u^s$ exists since $G^s$ is assumed to have a density).
\end{theorem}


Next, we state the main result of the paper, which shows that, under
fairly general conditions, a solution to the fluid equations exists and
is the functional law of large numbers limit, as $N\rightarrow\infty
$, of the $N$-server system with abandonment.
\begin{theorem} \label{thm:2}
Suppose that Assumptions \ref{ass-init}--\ref{ass-h} hold, and let
$(\overline{E}, \overline{X}(0)$,\break $\overline{\nu}_0,\overline{\eta
}_0)\in\mathcal{S}_0$ be the limiting
initial condition. Then there exists a unique solution $(\overline
{X},\overline{\nu}
,\overline{\eta})$ to the associated fluid equations, and the
sequence $(\overline{X}{}^{(N)}
,\overline{\nu}{}^{(N)},\break\overline{\eta}^{(N)})$ converges weakly, as
$N\rightarrow\infty$, to
$(\overline{X},\overline{\nu},\overline{\eta})$.
\end{theorem}

Theorem \ref{thm:2} follows from Theorem \ref{th-tight}, which establishes
tightness of the sequence $\{\overline{X}{}^{(N)},\overline{\nu
}^{(N)},\overline{\eta}^{(N)}\}$, Theorem \ref{thm:FE},
which shows that any subsequential limit of the sequence
$\{\overline{X}{}^{(N)},\overline{\nu}{}^{(N)},\overline{\eta}^{(N)}\}
$ satisfies the fluid equations, and the uniqueness
of solutions to the fluid equations stated in Theorem \ref{thm:1}.
\begin{cor}
Suppose that Assumptions \ref{ass-init}--\ref{ass-h} hold. Given any
$(\overline{E}$, $\overline{X}(0), \overline{\nu}_0,\overline{\eta
}_0) \in\mathcal{S}_0$, let $(\overline{X}
,\overline{\nu},\overline{\eta})$ be the unique solution to the
associated fluid
equations (\ref{cond-radon})--(\ref{eq-fnonidling}) specified in
Theorem \ref{thm:1}. If the function $\overline{E}$ is absolutely continuous
and $\overline{\nu}_0$ and $\overline{\eta}_0$ are absolutely
continuous measures,
then the function $\overline{X}$ is also absolutely continuous and for every
$t\in[0,\infty)$, the measures $\overline{\nu}_t$ and $\overline
{\eta}_t$ are also
absolutely continuous.
\end{cor}
\begin{pf}
Since $\overline{E}$ is absolutely continuous, (\ref{dis:fx}) allows
us to
deduce that $\overline{X}$ is absolutely continuous. In turn, (\ref{dis:fk})
shows that $\overline{K}$ is also absolutely continuous. Then the
argument used
in proving Lemma 5.18 of \cite{kasram07} can be adapted, together with
(\ref{eq-freneg2}) and (\ref{eq-fmeas2}), to show that $\overline
{\nu}_t$ and
$\overline{\eta}_t$ are absolutely continuous for every $t\in
[0,\infty)$. This
proves the corollary.
\end{pf}

We now state the fluid limit result for the virtual waiting time process
$W^{(N)}$. This result is of particular interest in the context of call
centers.
Note that in the fluid system, for any $u > t$
the total mass of customers in queue at time
$u$ that arrived before time $t$ equals $\overline{Q}(u) -
\overline{\eta}_{u}[0,u-t]$, and the ages of these (fluid) customers
lie in the interval $(u-t, \overline\chi(u)]$, where
%
%
\begin{equation}\overline\chi(u) \doteq(F^{\overline{\eta
}_{u}})^{-1}(\overline{Q}
(u)).
\end{equation}
Observe that this definition is analogous to the definition of
$\chi^{(N)}$ given in (\ref{def-chi}). Therefore, by the same logic
that was used to justify the
expression (\ref{fr}) for $\overline{R}$ in Definition \ref{def-fleqns},
it is natural to conjecture that, for each
$t \in[0,\infty)$, the fluid limit of the sequence $\{\overline
{\mathcal
T}{}^{(N)}_t\}$ equals $\overline
{\mathcal T}_t$, where for $s\in[0,\infty)$,
%
%
\begin{eqnarray}
\label{dis:T}
\overline
{\mathcal T}_t(s) & \doteq&
\int_t^{t+s} \biggl( \int_{\overline{\eta}_{u}[0,u-t]}^{\overline{Q}(u)}
h^r( (F^{\overline\eta_{u}})^{-1}(y)) \,dy \biggr) \,du \nonumber\\[-8pt]\\[-8pt]
& = & \int_0^s \biggl( \int_{\overline{\eta}_{t+u}[0,u]}^{\overline{Q}(t
+u)}h^r( (F^{\overline\eta_{t+u}})^{-1}(y)) \,dy \biggr) \,du.
\nonumber
\end{eqnarray}
Also, define
%
%
\begin{equation}\label{Wbar} \overline W(t) \doteq\inf\biggl\{s\geq
0\dvtx\int
_t^{t+s} \langle h^s, \overline{\nu}_u \rangle \,du + \overline
{\mathcal T}_t(s)
\geq
\overline{Q}(t) \biggr\}.
\end{equation}

We will say a function $f\in\mathcal{D}[0,\infty)$ is uniformly strictly
increasing if it is absolutely continuous and there exists $a>0$ such
that the derivative of $f$ is bigger than and equal to $a$ for a.e.
$t\in[0,\infty)$. Note that for
any such
function, $f^{-1}(f(t))=t$ and $f^{-1}$ is continuous and strictly
increasing on $[0,\infty)$. We now characterize the fluid limit of the
(scaled) virtual waiting time in the system.
\begin{theorem} \label{thm:3}
Suppose that the conditions of Theorem \ref{thm:2} hold and that the
function $\int_0^{\cdot} \langle h^s, \overline{\nu}_u \rangle \,du$
is uniformly
strictly increasing. For each $t \geq0$, if $\overline{Q}$ is
continuous at
$t$, then $\overline{\mathcal
T}{}^{(N)}_t
\Rightarrow\overline{\mathcal T}_t$ and
$W^{(N)}(t) \Rightarrow\overline W(t)$, as
$N\rightarrow\infty$.
\end{theorem}

\section{Uniqueness of solutions to the fluid equations}
\label{sec-uniq}

In Section \ref{subs-cont}, we show that if $(\overline{X}, \overline
{\nu}, \overline{\eta})$
solve the fluid equations associated with a given initial condition
$(\overline{E},\overline{X}(0),\overline{\nu}_{0},\overline{\eta
}_{0})\in\mathcal{S}_0$,
then $\overline{\nu}$ (resp., $\overline{\eta}$) can be written
explicitly in terms
of the auxiliary fluid process $\overline{K}$ (resp., cumulative
arrival process $\overline{E}$). In
Section \ref{sec-flsol}, these representations are
used, along with the nonidling condition and the remaining fluid equations,
to show that there is at most one solution to the fluid
equations for a given initial condition.

\subsection{Integral equations for $(\overline{\nu},\overline{K})$
and $(\overline{\eta},\overline{E})$}
\label{subs-cont}

We begin by recalling Theorem~4.1 and Remark 4.3 of \cite{kasram07},
which we state here as Proposition \ref{th-pde}. This proposition
identifies an implicit relation that must be satisfied by the processes
$(\overline{\nu},\overline{K})$ and $(\overline{\eta},\overline
{E})$ that solve (\ref{eq-ftmeas}) and
(\ref{eq-ftreneg}), respectively.
\begin{prop}[\cite{kasram07}]
\label{th-pde} Let $G$ be the cumulative distribution function of a
probability distribution with density $g$ and hazard rate function
$h=g/(1-G)$, let $H\doteq\sup\{x\in[0,\infty)\dvtx g(x)>0\}$. Suppose
$\overline\pi\in\mathcal{D}_{\mathcal{M}_F[0,H)}[0,\infty)$ has
the property
that for every $m\in[0,H)$ and
$T\in[0,\infty)$, there exists $C(m,T)<\infty$ such that
%
%
\begin{equation}
\label{dis:hcontr} \int_0^\infty\langle\varphi(\cdot,s)h(\cdot
),\overline
\pi_s \rangle\,ds<C(m,T)\Vert\varphi\Vert_\infty
\end{equation}
for every $\varphi\in\mathcal{C}_c({\mathbb R}^2)$ with
$\operatorname{supp}(\varphi)\subset
[0,m]\times[0,T]$. Then given
any $\overline\pi_0 \in
\mathcal{M}_F[0,H)$ and $\overline{Z}\in\mathcal{I}_{{\mathbb
R}_+}[0,\infty)$, $\overline\pi$ satisfies the
integral equation
%
%
\begin{eqnarray}\label{eq-pde}\qquad
\langle
\varphi(\cdot,t), \overline\pi_t \rangle& = & \langle\varphi
(\cdot,0),
\overline\pi_0
\rangle+ \int_0^t \langle\varphi_s(\cdot,s), \overline\pi_s
\rangle \,ds +
\int_0^t \langle\varphi_x(\cdot,s), \overline\pi_s \rangle \,ds
\nonumber\\[-8pt]\\[-8pt]
& &{} - \int_0^t \langle\varphi(\cdot,s)h(\cdot), \overline\pi_s
\rangle \,ds + \int_{[0,t]} \varphi(0,s) \,d \overline{Z}(s)\nonumber
\end{eqnarray}
for
every $\varphi\in\mathcal{C}_c((-\infty,H)\times{\mathbb R})$ and
$t \in
[0,\infty)$, if and only if $\overline\pi$ satisfies
%
%
\begin{eqnarray}\label{eq2-fmeas}
\int_{[0,M)} f (x) \overline\pi_t (dx) &=& \int_{[0,M )} f(x+t)
\frac{1 - G(x+t)}{1 - G(x)} \overline\pi_0 (dx)\nonumber\\[-8pt]\\[-8pt]
&&{} + \int_{[0,t]}
f(t-s) \bigl(1
- G(t-s)\bigr) \,d \overline{Z}(s),\nonumber
\end{eqnarray}
for every $f \in\mathcal{C}_b({\mathbb R}_+)$ and $t\in
(0,\infty)$. Moreover, for every $f \in\mathcal{C}^1_b({\mathbb
R}_+)$ and $t\in(0,\infty
)$,
%
%
\begin{eqnarray}
\label{eq-ibp}
&&\int_0^t f(t-s) \bigl(1 - G(t-s)\bigr) \,d \overline{Z}(s) \nonumber\\
&&\qquad=
f(0) \overline{Z}(t) + \int_{[0,t]} f^\prime(t-s) \bigl(1 -
G(t-s)\bigr)\overline{Z}(s)
\,ds \\
&&\qquad\quad{}- \int_{[0,t]} f(t-s)
g(t-s)\overline{Z}(s) \,ds.\nonumber
\end{eqnarray}
\end{prop}

Fluid equations (\ref{cond-radon})--(\ref{eq-ftreneg}) show that
(\ref{dis:hcontr}) and (\ref{eq-pde}) are satisfied with
$(h,\overline
\pi,\overline Z)$ replaced by $(h^s,\overline{\nu},\overline{K})$
and $(h^r,\overline{\eta}
,\overline{E})$,
respectively. Therefore, the next result
follows from Proposition \ref{th-pde}.
\begin{cor} \label{cor:nueta} Processes $(\overline{\eta},\overline
{E})$ and $(\overline{\nu}
,\overline{K})$ satisfy (\ref{eq-freneg2}) and (\ref{eq-fmeas2}) for every
bounded Borel measurable function $f$ and $t\in[0,\infty)$.
Moreover, $\overline{K}$ satisfies the renewal equation
%
%
\begin{eqnarray}\label{dis:fkl}
\overline{K}(t)&=&\langle{\mathbf{1}}, \overline{\nu}_t \rangle-
\langle{\mathbf{1}}, \overline{\nu}
_0 \rangle
+\int_{[0,H^s )}
\frac{G^s(x+t)- G^s(x)}{1 - G^s(x)} \overline{\nu}_0 (dx)\nonumber\\[-8pt]\\[-8pt]
&&{} + \int_0^t
g^s(t-s)\overline{K}(s) \,ds\nonumber
\end{eqnarray}
for each $t\geq0$ and admits the representation
\begin{eqnarray*}
\overline{K}(t) &=& \int_{[0,t]} (\langle{\mathbf{1}}, \overline
{\nu}_{t-s} \rangle
- \langle
{\mathbf{1}}, \overline{\nu}_0 \rangle) \,dU^s(s) \\
& &{} +\int
_{[0,t]} \biggl(\int_{[0,H^s )}
\frac{G^s(x+t-s)- G^s(x)}{1 - G^s(x)} \overline{\nu}_0 (dx) \biggr) \,dU^s(s),
\end{eqnarray*}
where $dU^s$ is the renewal measure associated with the distribution $G^s$.
\end{cor}
\begin{remark}
Strictly speaking, in \cite{kasram07} the cumulative distribution
function $G$ was
assumed to be absolutely continuous and supported on $[0,\infty)$. However,
the proofs given there only use the local integrability
of the hazard rate function $h$ on $[0,H)$ and so continue
to hold for $G^r$ here, which may possibly have a positive mass at
$\infty$.
In fact, in the case that $G^r$ has a positive mass at $\infty$ the
hazard rate function $h^r$ is globally integrable on $[0,H^r)$.
\end{remark}

\subsection{Uniqueness of solutions}
\label{sec-flsol}

Let $(\overline{X},
\overline{\nu},\overline{\eta})$ be a solution to the fluid
equations associated with
$(\overline{E},
\overline{X}(0), \overline{\nu}_0,\overline{\eta}_0)$. Recall the
definitions of $\overline{Q}$ and
$\overline{R}$ that
are given in (\ref{fq}) and (\ref{fr}).
As an immediate consequence of (\ref{fr}), we have the following elementary
property.
\begin{lemma}
\label{fq0}
For any $0 \leq a \leq b < \infty$, if $\overline{Q}(t)=0$ for all
$t\in
[a,b]$, then $\overline{R}(b)-\overline{R}(a)=0$.
\end{lemma}

Next, we establish the intuitive result that the process $\overline
{K}$ that
represents the cumulative entry of ``fluid'' into service is nondecreasing.
\begin{lemma}
The function $\overline{K}$ is nondecreasing.
\end{lemma}
\begin{pf}
Fix $t\in[0,\infty)$ and $0\leq s<t$. If $\overline{X}(t)\geq1$,
then $\langle
{\mathbf{1}}, \overline{\nu}_t \rangle=1\geq\langle{\mathbf
{1}},\overline{\nu}_s \rangle
$ by
(\ref{eq-fnonidling}). Hence, by (\ref{eq-fk}), it follows that
%
%
\begin{equation}
\label{eq-kinc1} \overline{K}(t)-\overline{K}(s)=\langle{\mathbf
{1}}, \overline{\nu}_t \rangle-
\langle
{\mathbf{1}},
\overline{\nu}_s \rangle+ \int_s^t \langle h^s, \overline{\nu}_l
\rangle \,dl \geq0.
\end{equation}
If
$\overline{X}(t)< 1$, we consider two cases.

\textit{Case} 1. $\overline{X}(v)<1$ for all $v\in(s,t]$. In this case,
by (\ref{fq}) and (\ref{eq-fnonidling}), $\overline{Q}(v)=0$ for all
$v\in
(s,t]$. Hence, by Lemma \ref{fq0} and the right continuity of
$\overline{R}$,
$\overline{R}(t)-\overline{R}(s)=0$. By (\ref{qt-conserve}), it then
follows that
\begin{eqnarray*}
\overline{K}(t)-\overline{K}(s)&=& \overline{K}(t)-\overline{K}(s)
+ \overline{R}(t)-\overline{R}(s) +\overline{Q}(t)-\overline{Q}(s)\\
&=& \overline{E}(t)-\overline{E}(s)\\ &\geq& 0.
\end{eqnarray*}

\textit{Case} 2. There exists $v\in(s,t]$ such that $\overline{X}
(v)\geq
1$. Define $l\doteq\sup\{v\leq t\dvtx\overline{X}(v)\geq1\}$.
Then, clearly
$l\in(s,t]$
and $\overline{X}(l-)\geq1$. Now, (\ref{fr}) implies that $\overline
{R}$ is
continuous and
hence, by (\ref{eq-fx}), $\overline{X}(v)-\overline{X}(v-)\geq0$
for every $v\in
(0,\infty)$. Therefore, $\overline{X}(l)\geq1=\langle{\mathbf
{1}},\overline{\nu}_l
\rangle$,
with the latter equality being a consequence of the nonidling
condition (\ref{eq-fnonidling}). Due to the case assumption $\overline
{X}(t)<1$,
we must have $l<t$. Then (\ref{eq-kinc1}), with $t$
replaced by $l$, shows that $\overline{K}(l)-\overline{K}(s)\geq0$.
On the other hand, since
$\overline{X}(v)<1$ for all $v\in(l,t]$, the argument in case 1 above
shows that
$\overline{K}(t)-\overline{K}(l)\geq0$. Thus, in this case too, we
have $\overline{K}(t)-\overline{K}
(s)\geq0$.
\end{pf}

We now state the main result of this section.
\begin{theorem}\label{th-unique2}
For $i=1,2$, let $(\overline{X}{}^i,\overline{\nu}^i,\overline{\eta
}^i)$ be a solution to the
fluid equations associated with $(\overline{E},\overline
{X}(0),\overline{\nu}_0,\overline{\eta}_0)
\in\mathcal{S}_0$. Then $\overline{X}{}^1=\overline{X}{}^2, \overline
{\nu}^1=\overline{\nu}^2$ and
$\overline{\eta}^1=\overline{\eta}{}^2$.
\end{theorem}
\begin{pf}
For each $i=1,2$, let $\overline{Q}{}^i, \overline{K}{}^i, \overline
{D}{}^i, \overline{R}{}^i$ be the processes
associated with the solution $(\overline{X}{}^i,\overline{\nu
}^i,\overline{\eta}^i)$ to the fluid
equations for $(\overline{E},\overline{X}(0),\overline{\nu
}_0,\overline{\eta}_0)\in\mathcal{S}_0$. It
follows directly from Corollary \ref{cor:nueta} that
$\overline{\eta}^1=\overline{\eta}^2$. Let $\triangle A$ denote
$A^2-A^1$ for $A=\overline{Q},
\overline{K}, \overline{D}, \overline{R}$ and $\overline{\nu}$.
For each $t\geq0$, let $\triangle
\overline{\nu}_t$ be the measure that satisfies $\triangle
\overline{\nu}_t(\Xi)=\overline{\nu}_t^2(\Xi)-\overline{\nu
}_t^1(\Xi)$ for every
measurable set
$\Xi\subset[0,\infty)$. Choose $\delta>0$ and define
\[
\tau=\tau_\delta\doteq\inf\{t\geq0\dvtx\triangle\overline
{K}(t)\vee
\triangle
\overline{K}(t-)\geq\delta\}.
\]
We shall argue by contradiction to show that
$\tau=\infty$. Suppose that $\tau<\infty$.

We first claim that for each $t\in[0,\tau]$,
%
%
\begin{equation}\label{claim1}
\triangle\overline{K}(t)< \delta\qquad\mbox{if } \langle{\mathbf{1}},
\overline{\nu}^1_t
\rangle
=1.
\end{equation}
To see why this is true, suppose that $\langle{\mathbf{1}}, \overline
{\nu}
^1_t \rangle=1$ for some $t\in[0,\tau]$. Since $\langle{\mathbf{1}},
\overline{\nu}
^2_t \rangle\leq1$, we have $\langle{\mathbf{1}}, \triangle
\overline{\nu}_t
\rangle\leq
0$. When combined with (\ref{dis:fkl}) and the identity $\triangle
\overline{\nu}_0=0$, this shows that
%
%
\begin{equation} \label{dis:Kest}\qquad
\triangle\overline{K}(t)= \langle{\mathbf{1}}, \triangle\overline
{\nu}_t \rangle+ \int_0^t
g^s(t-s)\triangle\overline{K}(s) \,ds \leq\int_0^t g^s(t-s)\triangle
\overline{K}(s) \,ds.
\end{equation}
If $G^s(t)>0$ then, along with the fact that $\triangle\overline
{K}(s)<\delta
$ for all
$s\in[0,t)$, this implies $\triangle\overline{K}(t) <\delta
G^s(t)\leq\delta
$. On
the other hand, if $G^s(t)=0$, it must be that $g^s(s)=0$ for a.e.
$s\in
[0,t]$ and so (\ref{dis:Kest}) implies that $\triangle\overline
{K}(t)=0\leq
\delta$. Thus (\ref{claim1}) follows in either case. In addition, the
right-continuity of $\overline{K}{}^1$ and $\overline{K}{}^2$ implies
that $\triangle\overline{K}
(\tau)\geq\delta$. When combined with (\ref{claim1}), (\ref{fq})
and (\ref{eq-fnonidling}), this shows that
%
%
\begin{equation}\label{claim2} \overline{X}{}^1(\tau) = \langle
{\mathbf{1}},
\overline{\nu}^1_\tau\rangle
<1 \quad\mbox{and}\quad \overline{Q}{}^1(\tau)=0.
\end{equation}
Now, define
\[
r\doteq\sup\{t<\tau\dvtx\overline{Q}{}^2(t)<\overline{Q}{}^1(t) \}
\vee0.
\]
Then for every $t\in[r,\tau]$, $\overline{Q}{}^2(t)\geq\overline{Q}{}^1(t)$.
If $r=0$, then $\triangle\overline{K}(r)=\triangle\overline
{K}(0)=0<\delta$. On the
other hand, if $r>0$, there exists a sequence of $\{t_n\}_{n=1}^\infty
$ such that $t_n<r$ and $t_n\rightarrow r$ as $n\rightarrow\infty$
and $0\leq\overline{Q}{}^2(t_n)<\overline{Q}{}^1(t_n)$ for each $n\in
{\mathbb N}$. Since $\overline{Q}{}^1$
and $\overline{Q}{}^2$ are c\`{a}dl\`{a}g, this implies that
%
%
\begin{equation}
\label{uniq-qineq}\overline{Q}{}^2(r-)\leq\overline{Q}{}^1(r-),
\end{equation}
and, due to (\ref{fq}) and
(\ref{eq-fnonidling}), it also follows that $\overline{X}{}^1(t_n)>
\langle\mathbf{1}, \overline{\nu}^1_{t_n} \rangle=1$ for every
$n\in{\mathbb N}$. When combined with
(\ref{dis:Kest}), this shows that for $n\in{\mathbb N}$,
\[
\triangle\overline{K}(t_n)\leq\int_0^{t_n} g^s(t_n-s)\triangle
\overline{K}(s)
\,ds=\int_0^{t_n} g^s(s)\triangle\overline{K}(t_n-s) \,ds.
\]
Since $\overline{K}{}^1$ and $\overline{K}{}^2$ are c\`{a}dl\`{a}g, this
implies that
\[
\triangle\overline{K}(r-)\leq\int_0^{r} g^s(s)\triangle\overline
{K}\bigl((r-s)-\bigr) \,ds.
\]
Using the fact that $\triangle\overline{K}((r-s)-)<\delta$ for all
$s\in
(0,r)$, it is easy to see [once again, as in the analysis of (\ref
{dis:Kest}), by considering the cases $G^s(r)>0$ and $G^s(r)=0$
separately] that this implies
%
%
\begin{equation}\label{uniq-kineq} \triangle\overline{K}
(r-)<\delta.
\end{equation}

On the other hand, since (\ref{qt-conserve}) is satisfied with
$(\overline{K},\overline{R},\overline{Q})$ replaced by $(\overline{K}{}^i,\overline
{R}{}^i$, $\overline{Q}{}^i)$ for $i=1,2$, and
$\triangle\overline{Q}(0) + \triangle\overline{E}(t)=0$ for each
$t\geq0$, it
follows that
\[
\triangle\overline{K}(\tau)+\triangle\overline{R}(\tau)+\triangle
\overline{Q}(\tau
)=\triangle\overline{K}(r-)+\triangle\overline{R}(r-)+\triangle
\overline{Q}(r-)=0.
\]
Hence,
\[
\triangle\overline{K}(\tau) - \triangle\overline{K}(r-) =
-\bigl(\triangle\overline{R}(\tau) -
\triangle\overline{R}(r-) \bigr) -\triangle\overline{Q}(\tau)+\triangle
\overline{Q}(r-).
\]
Since
$-\triangle\overline{Q}(\tau)=\overline{Q}{}^1(\tau)-\overline
{Q}^2(\tau)=-\overline{Q}{}^2(\tau) \leq
0$ due to
(\ref{claim2}) and $\triangle\overline{Q}(r-)\leq0$ by (\ref{uniq-qineq}),
we obtain
%
%
\begin{equation}\label{KR}
\triangle\overline{K}(\tau) - \triangle\overline{K}(r-) \leq
-\bigl(\triangle
\overline{R}(\tau) - \triangle\overline{R}(r-) \bigr).
\end{equation}
We now show that the right-hand side of the above display is
nonpositive. For each $t\geq0$, by (\ref{fr}), we see that
\begin{eqnarray*}
\triangle\overline{R}(t)&=& \overline{R}{}^2(t)-\overline{R}{}^1(t) \\
&=& \int_{0}^{t}
\biggl(\int_0^{\overline{Q}{}^2(s)} h^r((\overline F^{\overline\eta
^2_s})^{-1}(y))\,dy \biggr) \,ds\\
&&{} -\int_{0}^{t} \biggl(\int_0^{\overline{Q}{}^1(s)}
h^r((\overline F^{\overline\eta^1_s})^{-1}(y))\,dy \biggr) \,ds.
\end{eqnarray*}
Since $\overline{\eta}^1=\overline{\eta}^2$, it follows that
$\overline F^{\overline
\eta^1_{\cdot}}=\overline F^{\overline\eta^2_{\cdot}}$. Together
with the continuity of $\overline{R}{}^1$ and~$\overline{R}{}^2$, this
yields the equation
%
%
\begin{eqnarray}
\label{R12}
&& \triangle\overline{R}(\tau)-\triangle\overline
{R}(r-)\nonumber\\
&&\qquad=
\triangle\overline{R}(\tau)-\triangle\overline{R}(r)\nonumber\\[-8pt]\\[-8pt]
&&\qquad= \int_{r}^{\tau
} \biggl(\int_0^{\overline{Q}{}^2(s)} h^r((\overline F^{\overline\eta
^1_s})^{-1}(y))\,dy \biggr) \,ds\nonumber\\
&&\qquad\quad{} -\int_{r}^{\tau} \biggl(\int_0^{\overline{Q}
^1(s)} h^r((\overline F^{\overline\eta^1_s})^{-1}(y))\,dy \biggr) \,ds.
\nonumber
\end{eqnarray}
However, by the definition of $r$, for each $t\in[r,\tau]$,
$\overline{Q}
^2(t)\geq\overline{Q}{}^1(t)$, and so $\triangle\overline{R}(\tau
)-\triangle\overline{R}
(r-)\geq0$. Together with (\ref{KR}) and (\ref{uniq-kineq}), this implies
\[
\triangle\overline{K}(\tau)\leq\triangle\overline{K}(r-) <\delta.
\]
Essentially the same argument can be used to also show that $\triangle
\overline{K}(\tau-)\leq\triangle\overline{K}(r-)<\delta$. Hence
$\triangle\overline{K}
(\tau)\vee\triangle\overline{K}(\tau-)<\delta$, which contradicts the
definition of $\tau$.

Thus we have proved that $\tau=\infty$ and $\overline
{K}^2(t)-\overline{K}{}^1(t)\leq
\delta$ for each $\delta>0$ and $t\geq0$. By letting $\delta
\rightarrow0$, we have $\overline{K}{}^2(t)\leq\overline{K}{}^1(t)$ for
all $t\geq0$. An
exactly analogous argument yields the reverse inequality $\overline
{K}^1(t)\leq
\overline{K}{}^2(t)$ for each $t\geq0$, and so it must be that
$\overline{K}{}^2= \overline{K}{}^1$.
By Corollary \ref{cor:nueta}, it follows that $\overline{\nu
}^1=\overline{\nu}^2$.
Also, by (\ref{qt-conserve}), we obtain
%
%
\begin{equation}\label{frq}\overline{R}{}^1+\overline{Q}
^1=\overline{R}{}^2+\overline{Q}{}^2.
\end{equation}
We now show that, in fact $\overline{Q}{}^1=\overline{Q}{}^2$ and
$\overline{R}{}^1=\overline{R}{}^2$. If there exists $t\in(0,\infty)$
such that $\overline{Q}
^1(t)>\overline{Q}{}^2(t)$, let
\[
s\doteq\sup\{v<t\dvtx\overline{Q}{}^1(v)\leq\overline{Q}{}^2(v)\}\vee0.
\]
Then $\overline{Q}{}^1(s-)\leq\overline{Q}{}^2(s-)$ and $\overline
{Q}^1(v)>\overline{Q}{}^2(v)$ for each
$v\in(s,t]$. Due to the fact that $\overline\eta^1=\overline\eta
^2$, we have
\begin{eqnarray*}
\overline{R}{}^1(t)-\overline{R}{}^1(s)&=&\int_{s}^{t} \biggl(\int
_0^{\overline{Q}{}^1(v)}
h^r((\overline F^{\overline\eta^1_v})^{-1}(y))\,dy \biggr) \,dv
\\
&\geq&
\int_{s}^{t} \biggl(\int_0^{\overline{Q}{}^2(v)} h^r((\overline
F^{\overline
\eta^2_v})^{-1}(y))\,dy \biggr) \,dv\\
&=& \overline{R}{}^2(t)-\overline{R}{}^2(s).
\end{eqnarray*}
From (\ref{frq}) and the continuity of $\overline{R}{}^i, i=1,2$, we deduce
that $\overline{Q}{}^1(t)-\overline{Q}{}^1(s-)\leq\overline
{Q}{}^2(t)-\overline{Q}{}^2(s-)$. Combining this
with the inequality $\overline{Q}{}^1(s-)\leq\overline{Q}{}^2(s-)$
proved above, we obtain
$\overline{Q}{}^1(t)\leq\overline{Q}{}^2(t)$, which leads to a
contradiction. Hence $\overline{Q}
^1(v)\leq\overline{Q}{}^2(v)$ for all $v\in(0,\infty)$. By symmetry,
we can
also argue that $\overline{Q}{}^1(v)\geq\overline{Q}{}^2(v)$ for all
$v\in(0,\infty)$.
This shows $\overline{Q}{}^1= \overline{Q}{}^2$ and, hence, $\overline
{R}{}^1= \overline{R}{}^2$. Finally, by
(\ref{fq}), we have $\overline{X}{}^1=\overline{X}{}^2$.
\end{pf}
\begin{pf*}{Proof of Theorem \ref{thm:1}}
The first statement in Theorem \ref{thm:1} follows from Theorem
\ref{th-unique2}. The second statement follows directly from Corollary
\ref{cor:nueta} and the fluid equations (\ref{fq}), (\ref{fr}),
(\ref{eq-fx}) and (\ref{fqfx}). The alternative representation
(\ref{eq-fk2}) for $\overline{K}$ is a direct consequence of the renewal
equation (\ref{dis:fkl}) and the fact that the first three terms on the
right-hand side of (\ref{dis:fkl}) are bounded by one.
\end{pf*}
\begin{remark}
\label{rem-compen0}
For future use, we observe here that
the result of Lemma 5.16 in \cite{kasram07} (and the analog with
$\overline{\nu}$ replaced by $\overline{\eta}$), which was
obtained for the model without abandonments,
is also valid in the present context.
This is because equations (\ref{eq-freneg2}) and (\ref{eq-fmeas2}) of
Theorem \ref{thm:1} and Corollary 4.14 of \cite{kasram07} show that,
in the terminology of \cite{kasram07}, $\{\overline{\eta}_s\}$
(resp., $\{\overline{\nu}_s\}$) satisfies
the simplified age equation associated with a certain Radon measure
$\xi(\overline{\eta}_0$, $\overline{E})$ and $h^r$
[resp.,
$\xi(\overline{\nu}_0$, $\overline{K})$ and $h^s$]. Therefore, by
Proposition 4.15
of \cite{kasram07}, it follows that the result of Lemma 5.16 of
\cite{kasram07} is also valid in the present context.
\end{remark}

\section{A family of martingales}
\label{subs-prelim}

In Section \ref{subs-comp}, we identify the compensators (with respect
to the
filtration $\mathcal{F}_t^{(N)}$) of the cumulative
departure, potential reneging and (actual) reneging processes. Then,
in Section \ref{subs-altcomp}, we establish a more convenient representation
for the compensator of the reneging process.

\subsection{Compensators}
\label{subs-comp}

For any bounded measurable function $\varphi$ on $[0,H^s)\times
{\mathbb R}
_+$, consider
the sequence $\{A^{(N)}_{\varphi,\nu}\}$ of processes given by
%
%
\begin{equation}
\label{def-dcompns}
A^{(N)}_{\varphi,\nu} (t) \doteq\int_0^t \biggl(
\int_{[0, H^s)} \varphi(x,s) h^s(x) \nu^{(N)}_s (dx) \biggr)
\,ds,\qquad
t \in[0,\infty).
\end{equation}
Likewise, for any bounded measurable function $\varphi$ on
$[0,H^r)\times{\mathbb R}_+$
and $N \in{\mathbb N}$, let
%
%
\begin{equation}
\label{def-dcompnr}
A^{(N)}_{\varphi,\eta} (t) \doteq\int_0^t \biggl(
\int_{[0, H^r)} \varphi(x,s) h^r(x) \eta^{(N)}_s (dx) \biggr)
\,ds,\qquad
t \in[0,\infty).
\end{equation}
In Proposition \ref{cor-compensatormeasn}, we show that
$A^{(N)}_{\varphi,\nu}$ (resp., $A^{(N)}_{\varphi,\eta}$) is
the $\mathcal{F}_t^{(N)}$-compensator
of the associated ``$\varphi$-weighted'' cumulative departure process
$D^{(N)}_\varphi$
(resp., $S^{(N)}_\varphi$).
A~similar result was established
in \cite{kasram07} for the model without abandonments. However,
the filtration $\{\mathcal{F}_t^{(N)}\}$ considered here
is larger than the one considered in \cite{kasram07}, and
so Proposition \ref{cor-compensatormeasn} does not directly
follow from the results in \cite{kasram07}.
\begin{prop}
\label{cor-compensatormeasn} The following properties hold:
\begin{enumerate}
\item For every bounded measurable function $\varphi$ on
$[0,H^s)\times
{\mathbb R}_+$ such that the function $s \mapsto\varphi
(a^{(N)}_j(s),s)$ is left
continuous on $[0,\infty)$ for each $j$, the process
$M^{(N)}_{\varphi,\nu}$ defined by
%
%
\begin{equation}\label{def-martnmeasn}
M^{(N)}_{\varphi,\nu} \doteq D^{(N)}_\varphi- A^{(N)}_{\varphi
,\nu}
\end{equation}
is a local $\mathcal{F}_t^{(N)}$-martingale.
Moreover, for every $N \in{\mathbb N}$, $t \in[0,\infty)$ and
$m \in[0,H^s)$,
%
%
\begin{equation}
\label{bd-3}
\bigl|A^{(N)}_{\varphi,\nu} (t)\bigr| \leq\Vert\varphi\Vert_{\infty}
\bigl( X^{(N)}(0) + E^{(N)}(t) \bigr) \biggl( \int_0^m h^s(x) \,dx
\biggr) < \infty
\end{equation}
for every $\varphi\in\mathcal{C}_{c} ([0,H^s) \times
{\mathbb R}_+)$ with $\operatorname{supp}(\varphi) \subset[0,m]
\times{\mathbb R}_+$.
In addition, the quadratic variation process $\langle\overline
{M}{}^{(N)}_{\varphi
,\nu} \rangle$ of the scaled process $\overline{M}{}^{(N)}_{\varphi
,\nu}\doteq
M^{(N)}_{\varphi,\nu}/N$ satisfies
%
%
\begin{equation}\label{lim-qv}
\lim_{N\rightarrow\infty}\mathbb{E}\bigl[ \bigl\langle\overline
{M}{}^{(N)}_{\varphi,\nu} \bigr\rangle
(t) \bigr] =
0;\qquad \overline{M}{}^{(N)}_{\varphi,\nu}\Rightarrow\mathbf{0}\qquad\mbox{as
} N \rightarrow\infty.
\end{equation}
%
%
\item
Furthermore, properties (\ref{def-martnmeasn})--(\ref{lim-qv}) also hold
with $D, a_j, \nu, H^s$ and $h^s$, respectively, replaced
by $S, w_j, \eta, H^r$ and $h^r$.
\end{enumerate}
\end{prop}
\begin{pf}
In Lemma 5.4 and Corollary 5.5 of \cite{kasram07}, it was shown that
$A^{(N)}_{\varphi,\nu}$ is the compensator of $D^{(N)}_\varphi$
with respect
to a
certain filtration. The filtration $\{\mathcal{F}_t^{(N)}\}$ that we
consider here is larger than the filtration used in \cite{kasram07}
since it also includes the $\sigma$-algebra generated by the potential
waiting times $\{\eta^{(N)}_j(s), s\leq t, j=-\mathcal
{E}^{(N)}_0+1,\ldots,E^{(N)}(t)\}$.
Thus the results of \cite{kasram07} do not directly apply here.
Nevertheless, as we prove below, the result continues to hold due to
the assumed independence of the patience and service times.

We first claim that for every $\mathcal{F}_t^{(N)}$-stopping time
$\Upsilon$,
%
%
\begin{eqnarray}\label{eqn-KR}
&&
\mathbb{E}\bigl[\mathbh{1}_{\{\theta_n^k\leq
{j}/{2^m}<\Upsilon, \zeta_n^k>{j}/{2^m}\}}\mathbh{1}_{\{\zeta
_n^k\leq
({j+1})/{2^m}\}}|\mathcal{F}^{(N)}_{{j}/{2^m}}
\bigr]\nonumber\\[-8pt]\\[-8pt]
&&\qquad
=\mathbh{1}_{\{\theta_n^k\leq{j}/{2^m}<\Upsilon, \zeta
_n^k>{j}/{2^m}\}}\int_{j/2^m}^{(j+1)/2^m}\frac{g^s(u-\theta
_n^k)}{1-G^s({j}/{2^m}-\theta_n^k)}\,du,\nonumber
\end{eqnarray}
where $\theta_n^k$ (resp., $\zeta_n^k$) is the time at which the $n$th
customer to be served at station $k$ starts (resp., completes)
service. Then $\zeta_n^k-\theta_n^k$ is the service time of the $n$th
customer to be served at station $k$, which has cumulative distribution
function $G^s$. In order to show the equality in (\ref{eqn-KR}), it
suffices to show that for every bounded $\mathcal{F}^{(N)}_{
{j}/{2^m}}$-adapted random variable $H$,
%
%
\begin{eqnarray}\label{eqn-KR2}\quad
&&
\mathbb{E}\bigl[H\mathbh{1}_{\{\theta_n^k\leq
{j}/{2^m}<\Upsilon, \zeta_n^k>{j}/{2^m}\}}\mathbh{1}_{\{\zeta
_n^k\leq
({j+1})/{2^m}\}} \bigr] \nonumber\\[-8pt]\\[-8pt]
&&\qquad
= \mathbb{E}\biggl[H\mathbh{1}_{\{\theta
_n^k\leq{j}/{2^m}<\Upsilon, \zeta_n^k>{j}/{2^m}\}}\int
_{j/2^m}^{(j+1)/2^m}\frac{g^s(u-\theta_n^k)}{1-G^s(
{j}/{2^m}-\theta_n^k)}\,du \biggr].\nonumber
\end{eqnarray}
For $j\in{\mathbb N}$, $m\in{\mathbb N}$, define $\mathcal{G}_
{j/2^m}^{(N)}$
be the
$\sigma$-algebra to be generated by the events $\{(\theta_n^k\leq
x)\cap(\theta_n^k\leq\frac{j}{2^m},\zeta_n^k>\frac{j}{2^m}),
x\geq0\}$. In particular, $\mathcal{G}_{j/2^m}^{(N)}$ contains the
information of the ages of all customers in service at time $\frac
{j}{2^m}$. Recall that the patience times and the service times of
customers are assumed to be independent. Therefore, given $\mathcal
{G}_{j/2^m}^{(N)}$, $\zeta_n^k-\theta_n^k$ and $\mathcal
{F}^{(N)}_{{j/2^m}}$ are conditionally\vspace*{1pt} independent. Hence, it
follows from the left-hand side of (\ref{eqn-KR2}) that
\begin{eqnarray*}
&&\hspace*{-4pt}
\mathbb{E}\bigl[H\mathbh{1}_{\{\theta_n^k\leq
{j}/{2^m}<\Upsilon,
\zeta_n^k>{j}/{2^m}\}}\mathbh{1}_{\{\zeta_n^k\leq(
{j+1})/{2^m}\}}
\bigr] \\
&&\hspace*{-4pt}\qquad=
\mathbb{E}\bigl[\mathbb{E}\bigl[H\mathbh{1}_{\{{j}/{2^m}<\Upsilon\}
}\mathbh{1}_{\{
\theta_n^k\leq
{j}/{2^m},
\zeta_n^k>{j}/{2^m}\}}\mathbh{1}_{\{\zeta_n^k-\theta_n^k\leq
({j+1})/{2^m}-\theta_n^k\}}|\mathcal{G}_{j/2^m}^{(N)}
\bigr] \bigr] \\
&&\hspace*{-4pt}\qquad=
\mathbb{E}\bigl[\mathbb{E}\bigl[H\mathbh{1}_{\{{j}/{2^m}<\Upsilon\}
}|\mathcal
{G}_{j/2^m}^{(N)} \bigr]\\
&&\qquad\hspace*{15pt}{}\times\mathbb{E}\bigl[\mathbh{1}_{\{\theta_n^k\leq
{j}/{2^m},
\zeta_n^k>{j}/{2^m}\}}\mathbh{1}_{\{\zeta_n^k-\theta_n^k\leq
({j+1})/{2^m}-\theta_n^k\}}|\mathcal{G}_{j/2^m}^{(N)}
\bigr] \bigr]
\end{eqnarray*}
and
\begin{eqnarray*}
&&\mathbb{E}\bigl[\mathbh{1}_{\{\theta_n^k\leq{j}/{2^m},
\zeta_n^k>{j}/{2^m}\}}\mathbh{1}_{\{\zeta_n^k-\theta_n^k\leq
({j+1})/{2^m}-\theta_n^k\}}|\mathcal{G}_{j/2^m}^{(N)}
\bigr] \\
&&\qquad=\mathbh{1}_{\{\theta_n^k\leq{j/2^m},
\zeta_n^k>{j/2^m}\}}\int_{j/2^m}^{(j+1)/2^m}\frac
{g^s(u-\theta_n^k)}{1-G^s({j}/{2^m}-\theta_n^k)}\,du.
\end{eqnarray*}
Therefore,
\begin{eqnarray*}
&&\mathbb{E}\bigl[\mathbb{E}\bigl[H\mathbh{1}_{\{{j}/{2^m}<\Upsilon
\}}|\mathcal{G}_{j/2^m}^{(N)} \bigr]\mathbb{E}\bigl[\mathbh{1}_{\{
\theta
_n^k\leq{j}/{2^m}, \zeta_n^k>{j}/{2^m}\}}\mathbh{1}_{\{
\zeta
_n^k-\theta_n^k\leq({j+1})/{2^m}-\theta_n^k\}}|\mathcal{G}_
{j/2^m}^{(N)} \bigr] \bigr]
\\
&&\qquad= \mathbb{E}\biggl[\mathbb{E}
\bigl[H\mathbh{1}_{\{{j}/{2^m}<\Upsilon\}}|\mathcal{G}_
{j/2^m}^{(N)} \bigr]\mathbh{1}_{\{\theta_n^k\leq{j}/{2^m}, \zeta
_n^k>{j}/{2^m}\}}\\
&&\qquad\quad\hspace*{28.6pt}{}\times\int_{j/2^m}^{(j+1)/2^m}\frac{g^s(u-\theta
_n^k)}{1-G^s({j}/{2^m}-\theta_n^k)}\,du \biggr]
\\
&&\qquad= \mathbb{E}
\biggl[\mathbb{E}\biggl[H\mathbh{1}_{\{{j}/{2^m}<\Upsilon\}}\mathbh{1}_{\{
\theta
_n^k\leq{j}/{2^m}, \zeta_n^k>{j}/{2^m}\}}\\
&&\qquad\quad\hspace*{9.6pt}{}\times\int
_{j/2^m}^{(j+1)/2^m}\frac{g^s(u-\theta_n^k)}{1-G^s(
{j}/{2^m}-\theta_n^k)}\,du\Big|\mathcal{G}_{j/2^m}^{(N)}
\biggr] \biggr]
\\
&&\qquad= \mathbb{E}\biggl[H\mathbh{1}_{\{{j}/{2^m}<\Upsilon\}}\mathbh
{1}_{\{
\theta_n^k\leq{j}/{2^m}, \zeta_n^k>{j}/{2^m}\}}\int
_{j/2^m}^{(j+1)/2^m}\frac{g^s(u-\theta_n^k)}{1-G^s(
{j}/{2^m}-\theta_n^k)}\,du \biggr].
\end{eqnarray*}
This shows that (\ref{eqn-KR2}), and therefore (\ref{eqn-KR}), holds.

If $\varphi$ is bounded, measurable and such that the function
$s\mapsto
\varphi(a^{(N)}_j(s),s)$ is left
continuous for each $j$, then the process $\{\varphi(a^{(N)}_j(s),s),
s\geq0\}$ is $\mathcal{F}_t^{(N)}$-predictable. Therefore, it follows
from the standard theory
(cf. Theorem 3.18 of \cite{jacshibook}) that $M^{(N)}_{\varphi,\nu}$
is a local
$\mathcal{F}_t^{(N)}$-martingale. Inequality (\ref{bd-3}) can be established
exactly as in Proposition 5.7 of \cite{kasram07} and assertions
(\ref{lim-qv}) can be proved using the same argument as in Lemma 5.9 of
\cite{kasram07}, thus establishing property (1).
Due to the analogy between the service dynamics
and the potential queue dynamics (see Remark \ref{rem-compdyn}),
property (2) is a direct consequence of
property (1).
\end{pf}
\begin{remark}
\label{rem-compen}
It is easy to see that Lemmas 5.6 and 5.8 of \cite{kasram07}
continue to be valid in the presence of abandonments. Indeed, the proofs
of Lemmas 5.6 and 5.8 of \cite{kasram07} only depend on Assumption 1
and Corollary
5.5 therein (since, as shown in Lemma 5.12 of \cite{kasram07},
the additional conditions (5.32) and (5.33) of Lemma 5.8 of \cite
{kasram07} can be derived from
Assumption 1), which correspond to Assumption \ref{ass-init} and Proposition
\ref{cor-compensatormeasn} of this paper. In addition, due to the parallels
between the dynamics of $\nu^{(N)}$ and $\eta^{(N)}$ (see Remark \ref
{rem-compdyn}),
the analogs of the
results in Lemmas 5.6 and 5.8, with $D^{(N)}, \nu^{(N)}, G^s$ and $H^s$,
respectively, replaced by
$S^{(N)}, \eta^{(N)}, G^r$ and $H^r$, also hold. In this case, even though
$\eta^{(N)}_0$ is (unlike $\nu^{(N)}_0$) not necessarily a
sub-probability measure,
the verification of the conditions analogous to (5.32) and
(5.33) of Lemma 5.8 in \cite{kasram07} can still be carried out in
the same manner since
Assumption \ref{ass-init} implies that the sequence
$\{\langle{\mathbf{1}}, \eta^{(N)}_0\rangle\}$ is tight. Moreover,
even though
$G^r$ is
allowed to have a mass at $\infty$, the proofs of
Lemmas 5.6 and 5.8 are still valid, with the renewal function $U^s$ now
replaced by the function
$U^r(\cdot) = \int_0^\cdot\sum_{n=1}^\infty(g^r)^{*n} (s) \,ds$, where
$(g^r)^{*n}$ is the $n$th convolution of $g^r$ on $[0,\infty)$.
\end{remark}

Now, note from (\ref{RQ}) that $R^{(N)}=S^{(N)}_{\theta^{(N)}}$, where
$\theta^{(N)}$ is defined by (\ref{ps}). In view of the fact that
$A^{(N)}_{\varphi,\eta}$ is the compensator for
$S^{(N)}_{\varphi}$, it is natural to conjecture that the compensator of
$R^{(N)}$ is equal to $A^{(N)}_{\theta^{(N)},\eta}$, where
%
%
\begin{eqnarray}
\label{rep-rcomp1}
A^{(N)}_{\theta^{(N)},\eta} (t) \doteq\int_0^t \biggl(
\int_{[0, H^r)} \mathbh{1}_{[0,\chi^{(N)}(s-)]}(x) h^r(x) \eta^{(N)}_s
(dx) \biggr) \,ds,\nonumber\\[-8pt]\\[-8pt]
\eqntext{ t \in[0,\infty).}
\end{eqnarray}
However, this is not immediate from Proposition \ref
{cor-compensatormeasn}(2) since $\theta^{(N)}(w^{(N)}_j(\cdot),\cdot
)$ is not left continuous for any $j$. Instead, we approximate $\theta
^{(N)}$ by a sequence $\{\theta^{(N)}_m\}_{N\in{\mathbb N}}$ defined by
%
%
\begin{equation}
\label{def-thetam}\theta^{(N)}_m(x,s)\doteq\mathbh{1}_{(x-
{1}/{m},\infty)}\bigl(\chi^{(N)}(s-)\bigr),
\end{equation}
which is shown to be left
continuous in Lemma \ref{psim}. Then in Lemma \ref{cor:1}, we use an
approximation argument to show that $A^{(N)}_{\theta^{(N)},\eta}$
is indeed the compensator of $R^{(N)}$.
\begin{lemma} \label{psim}
For each $m\geq1$, $x\in{\mathbb R}$ and $s\in{\mathbb R}_+$, the
sequence $\{\theta
^{(N)}_m\}_{N\in{\mathbb N}}$ defined by (\ref{def-thetam}) satisfies the
following two properties:
\begin{enumerate}
\item For every $N\in{\mathbb N}, x\in{\mathbb R}, s\in{\mathbb R}$,
$\theta
^{(N)}_m(x,s)$ is nonincreasing in $m$ and converges, as $m\rightarrow
\infty$, to $\theta^{(N)}(x,s)$ for every sample point in $\Omega$.
\item For each $N,m\in{\mathbb R}$, $j\in{\mathbb Z}$, the process
$\theta
^{(N)}_m(w^{(N)}_j(\cdot),\cdot)$ has left continuous paths on
$(0,\infty)$.
\end{enumerate}
\end{lemma}
\begin{pf}
The first property is immediate from the definition of
$\theta^{(N)}_m$. For the second property, fix $N,m\in{\mathbb N}, s>0,
j\in
{\mathbb Z}$ and $\omega\in\Omega$. To ease the notation, we shall suppress
$\omega$ from the notation. Let $\{s_n\}$ be a sequence in $(0,\infty)$
such that $s_n\uparrow s$ as $n\rightarrow\infty$. We now consider two
mutually exclusive cases.

\textit{Case} 1. $\theta^{(N)}_m(w^{(N)}_j(s),s)=1$. Then
$w^{(N)}_j(s)< \chi^{(N)}(s-)+1/m$. Since $w^{(N)}_j$ is
nondecreasing, $w^{(N)}_j(s_n)\leq w^{(N)}_j(s)$ and since the process
$\{\chi^{(N)}(s-), s\geq0\}$ is left continuous, we have, for all
$n$ large enough, $w^{(N)}_j(s_n)<\chi^{(N)}(s_n-)+1/m$. Hence,
$\theta^{(N)}_m(w^{(N)}_j(s_n),s_n)=1$ for all $n\in{\mathbb N}$.
Thus, in
this case, $\theta^{(N)}_m(w^{(N)}_j(\cdot),\cdot)$ is left
continuous at $s$.

\textit{Case} 2. $\theta^{(N)}_m(w^{(N)}_j(s),s)=0$. Then
$w^{(N)}_j(s)\geq\chi^{(N)}(s-)+1/m$. It follows from Lemma \ref
{lem-chi} that for all sufficiently large $n$, $\chi^{(N)}(s-)-\chi
^{(N)}(s_n-)=s-s_n>0$. Since (\ref{def-waitjn}) implies
$w_j^{(N)}(s)-w_j^{(N)}(s_n)\leq s-s_n$ for all $n\in{\mathbb N}$, this
implies $w^{(N)}_j(s_n)\geq\chi^{(N)}(s_n-)+1/m$ for all $n$ large
enough. Hence, $\theta^{(N)}_m(w^{(N)}_j(s_n),s_n)=0$ and $\theta
^{(N)}_m(w^{(N)}_j(\cdot),\cdot)$ is again left continuous at $s$.
\end{pf}
\begin{lemma} \label{cor:1}
For every $N\in{\mathbb N}$,
the process $M^{(N)}_{\theta^{(N)},\eta}$ defined by
%
%
\begin{equation}
M^{(N)}_{\theta^{(N)},\eta} \doteq R^{(N)}- A^{(N)}_{\theta
^{(N)},\eta}
\end{equation}
is a local $\mathcal{F}_t^{(N)}$-martingale.
In addition, as $N
\rightarrow\infty$,
%
%
\begin{equation}\label{lim-psi}
\mathbb{E}\bigl[ \bigl\langle\overline
{M}{}^{(N)}_{\theta^{(N)},\eta}
\bigr\rangle(t) \bigr] \rightarrow 0,\qquad \overline{M}{}^{(N)}_{\psi,\eta} \Rightarrow
\mathbf{0}
\quad\mbox{and}\quad \overline{M}{}^{(N)}_{\theta^{(N)},\eta} \Rightarrow
\mathbf{0}.
\end{equation}
\end{lemma}
\begin{pf}
Fix $N\in{\mathbb N}$, and let $A^{(N)}_{\theta^{(N)}_m,\eta}$,
$m\in{\mathbb N}
$, be defined in the obvious way
%
%
\begin{equation}
\label{def-dcompnpsim}
A^{(N)}_{\theta^{(N)}_m,\eta} (t) \doteq\int_0^t \biggl(
\int_{[0, H^r)} \theta_m^{(N)}(x,s) h^r(x) \eta^{(N)}_s (dx)
\biggr) \,ds.
\end{equation}
By Proposition \ref{cor-compensatormeasn}(2) and Lemma \ref{psim}(2), the
process $A^{(N)}_{\theta^{(N)}_m,\eta}$ is the $\mathcal
{F}_t^{(N)}$-compen\-sator of
the process $S^{(N)}_{\theta_m^{(N)}}$, and the process
$M^{(N)}_{\theta^{(N)}_m,\eta}$ defined by
%
%
\begin{equation}\label{def-martnpsim}
M^{(N)}_{\theta^{(N)}_m,\eta} \doteq S^{(N)}_{\theta_m^{(N)}}-
A^{(N)}_{\theta
^{(N)}_m,\eta}
\end{equation}
is a local $\mathcal{F}_t^{(N)}$-martingale. Now, by Lemma \ref{psim}(1),
$\theta^{(N)}_m\rightarrow\theta^{(N)}$ pointwise on ${\mathbb R}^2_+$,
$|\theta^{(N)}_m(x,s)-\theta^{(N)}(x,s)|\leq1$ for all $(x,s)\in
{\mathbb R}^2_+$,
and $\mathbb{E}[S^{(N)}_{\mathbf{1}}(t) ]<\infty$,
$\mathbb{E}[A^{(N)}_{{\mathbf{1}},\eta}(t) ]<\infty$ for all
$t\in
(0,\infty)$. Hence, an application of the dominated convergence
theorem shows
that for all $t\in(0,\infty)$, as $m\rightarrow\infty$,
\[
\mathbb{E}\Bigl[\sup_{0\leq s\leq t} \bigl|A^{(N)}_{\theta^{(N)}_m,\eta}(s)-
A^{(N)}_{\theta^{(N)},\eta}(s) \bigr| \Bigr]\rightarrow0
\]
and
\[
\mathbb{E}\Bigl[\sup_{0\leq s\leq t} \bigl|S^{(N)}_{\theta_m^{(N)}}(s)-
S^{(N)}_{\theta^{(N)}}(s) \bigr|
\Bigr]\rightarrow
0,
\]
and hence\vspace*{-1pt} $M^{(N)}_{\theta^{(N)}_m,\eta}$ converges in distribution to
$M^{(N)}_{\theta^{(N)},\eta}$. Since $ |S^{(N)}_{\theta
_m^{(N)}}(t)-\break S^{(N)}_{\theta_m^{(N)}}
(t-) |\leq1$
for all $t\in[0,\infty)$ and $m\in{\mathbb N}$, we conclude that
$M^{(N)}_{\theta^{(N)},\eta}$ is a local $\mathcal
{F}_t^{(N)}$-martingale by
Corollary 1.19 of Chapter IX of \cite{jacshibook}.
Given that $M^{(N)}_{\theta^{(N)}, \eta}$ is a martingale,
the proof of the limits (\ref{lim-psi}) is exactly analogous to
the proof of~(\ref{lim-qv}), as carried out in
Lemma 5.9 of \cite{kasram07}.
\end{pf}

\subsection{An alternative representation for the compensator of $R^{(N)}$}
\label{subs-altcomp}

We now derive an alternative, more convenient representation for
$A^{(N)}_{\theta^{(N)},\eta}$, or more generally, for processes of
the form $A^{(N)}_{\theta^{(N)},\eta}$, but with $h^r$ replaced by
an arbitrary measurable function $h$.
In what follows, recall that $F^{\eta_t^{(N)}}(x)=\eta^{(N)}_t[0,x]$
and its inverse
$(F^{\eta_t^{(N)}})^{-1}$ is as defined in (\ref{inverse}).
\begin{prop} \label{lem:uni5}
For each $N\in{\mathbb N}$, $t\geq0$ and measurable function $h$ on
$[0,H^r)$,
%
%
\begin{eqnarray}\label{lem:ext1}
&&\int_{[0,H^r)}\mathbh{1}_{[0,\chi
^{(N)}(t-)]}(x)h(x)\eta^{(N)}_t(dx)
\nonumber\\[-8pt]\\[-8pt]
&&\qquad=\int_0^{Q^{(N)}(t)+\iota^{(N)}(t)}h((F^{\eta_t^{(N)}})^{-1}(y))\,dy,\nonumber
\end{eqnarray}
where
%
%
\begin{equation}
\label{iota}\iota^{(N)}(t)\doteq\cases{
0, &\quad if $\bigl(\chi^{(N)}(t-)-\chi^{(N)}(t)\bigr)\bigl(K^{(N)}
(t)-K^{(N)}(t-)\bigr)=0$,\vspace*{2pt}\cr
1, &\quad if
$\bigl(\chi^{(N)}(t-)-\chi^{(N)}(t)\bigr)\bigl(K^{(N)}(t)-K^{(N)}(t-)\bigr)>0$.}\hspace*{-30pt}
\end{equation}
%
\end{prop}
\begin{pf}
Fix $N\in{\mathbb N}$, $t\geq0$ and a measurable function $h$ on $[0,H^r)$.
By the representation (\ref{def-etan}) for $\eta^{(N)}_t$, we have
%
%
\begin{eqnarray}
\label{dis:2}
&&\int_{[0,H^r)}\mathbh{1}_{[0,\chi
^{(N)}(t-)]}(x)h(x)\eta^{(N)}_t(dx) \nonumber\\[-8pt]\\[-8pt]
&&\qquad= \sum_{j=-\mathcal{E}^{(N)}_0+
1}^{E^{(N)}(t)}
h \bigl(w^{(N)}_j (t) \bigr) \mathbh{1}_{\{w^{(N)}_j (t)\leq\chi
^{(N)}(t-)\}} \mathbh{1}_{\{ w^{(N)}_j (t) < r_j \}}.\nonumber
\end{eqnarray}
Moreover, by (\ref{qn}),
\[
Q^{(N)}(t)=\eta^{(N)}_t\bigl[0,\chi^{(N)}(t)\bigr]=\sum_{j=-\mathcal
{E}^{(N)}_0+ 1}^{E^{(N)}(t)}
\mathbh{1}_{\{w^{(N)}_j (t)\leq\chi^{(N)}(t)\}} \mathbh{1}_{\{
w^{(N)}_j (t) <
r_j \}}.
\]
Thus $Q^{(N)}(t)$ is the total number of customers who have arrived to
the system and have not reneged by $t$ and whose potential waiting
times at $t$ are less than or equal to $\chi^{(N)}(t)$. If we arrange
those customers in increasing order of their potential waiting times at
$t$, then for $i=1,2,\ldots, Q^{(N)}(t)$, $(F^{\eta
_t^{(N)}})^{-1}(i)$ is exactly
the potential waiting time at $t$ of the $i$th customer from the back
of the queue.

Suppose that $(\chi^{(N)}(t-)-\chi^{(N)}(t))(K^{(N)}(t)-K^{(N)}(t-))=0$.
This implies that either $\chi^{(N)}(t-)=\chi^{(N)}(t)$ holds or both
$\chi^{(N)}(t-)>\chi^{(N)}(t)$ and $K^{(N)}(t)=K^{(N)}(t-)$ hold. The
latter condition indicates that the head-of-the-line customer right
before time $t$ reneged at time $t$. In this case,
the right-hand side of (\ref{dis:2}) admits the alternative representation
\[
\int_0^{Q^{(N)}(t)}h((F^{\eta_t^{(N)}})^{-1}(y))\,dy.
\]

On the other hand, suppose that
$(\chi^{(N)}(t-)-\chi^{(N)}(t))(K^{(N)}(t)-K^{(N)}(t-))>0$. In this
case, the
head-of-the-line customer right before time $t$ departs for service at time
$t$ and this customer is counted in the right-hand side of (\ref
{dis:2}) but
not in $Q^{(N)}(t)$. Since $E^{(N)}(t)-E^{(N)}(t-)\leq1$, there is exactly
one such
customer, that is, $K^{(N)}(t)-K^{(N)}(t-)=1$. Hence the right-hand
side of
(\ref{dis:2})
can be rewritten as
\[
\int_0^{Q^{(N)}(t)+1}h((F^{\eta_t^{(N)}})^{-1}(y))\,dy.
\]
\upqed\end{pf}

As an immediate consequence of (\ref{rep-rcomp1}), Lemma \ref{cor:1},
and Proposition \ref{lem:uni5}, we obtain the following alternative
representation for the compensator $A^{(N)}_{\theta^{(N)},\eta}$
of $R^{(N)}$:
%
%
\begin{eqnarray}\label{rep-rcomp2}
A^{(N)}_{\theta^{(N)},\eta} (t) &\doteq& \int_0^t \biggl(\int
_0^{Q^{(N)}(t)+\iota^{(N)}(t)}h^r((F^{\eta_s^{(N)}})^{-1}(y))\,dy
\biggr)\,ds,\nonumber\\
\eqntext{t\in[0,\infty),}
\end{eqnarray}
where $\iota^{(N)}$ is given by (\ref{iota}).

\section{Tightness of pre-limit sequence}
\label{Sec:relcom}

The main objective of this section is to show that, under
suitable assumptions, the sequence of scaled state processes
$\{(\overline{X}{}^{(N)},\overline{\nu}{}^{(N)},\overline{\eta
}^{(N)})\}$ and the sequences of auxiliary
processes are tight.
Specifically, from (\ref{bd-1}) and (\ref{bd-3}) it is clear that for every
$t$, the linear functionals $\overline{D}{}^{(N)}_{\cdot}(t)\dvtx
\varphi\mapsto\overline{D}{}^{(N)}
_{\varphi}(t)$ and
$\overline{A}{}^{(N)}_{\cdot, \nu}(t)\dvtx\varphi\mapsto\overline
{A}^{(N)}_{\varphi, \nu}(t)$
are finite Radon measures on $[0,H^s)\times{\mathbb R}_+$. Likewise, from
(\ref{bd-2}) and the fact that
(\ref{bd-3}) holds with $\nu, h^s$, respectively, replaced by $\eta,
h^r$ by property (2) of Proposition \ref{cor-compensatormeasn}, it
follows that the linear functionals $\overline{S}{}^{(N)}_{\cdot
}(t)\dvtx\psi\mapsto
\overline{S}{}^{(N)}_{\psi}(t)$ and
$\overline{A}{}^{(N)}_{\cdot, \eta}(t)\dvtx\psi\mapsto\overline
{A}^{(N)}_{\psi, \eta}(t)$
define finite Radon measures on $[0,H^r)\times{\mathbb R}_+$. Thus $\{
\overline{D}{}^{(N)}
_{\cdot}(t)\dvtx t\in[0,\infty)\}$ and $\{\overline{A}{}^{(N)}_{\cdot
, \nu
}(t)\dvtx t\in[0,\infty)\}$
can be viewed as $\mathcal{M}_F([0,H^s)\times{\mathbb R}_+)$-valued
c\`{a}dl\`{a}g processes, and $\{\overline{S}{}^{(N)}
_{\cdot}(t)\dvtx t\in[0,\infty)\}$ and $\{\overline{A}{}^{(N)}_{\cdot
, \eta
}(t)\dvtx t\in[0,\infty)\}$
can be viewed as $\mathcal{M}_F([0,H^r)\times{\mathbb R}_+)$-valued
c\`{a}dl\`{a}g processes.
Now, for $N\in{\mathbb N}$, let
%
%
\begin{eqnarray}\label{Y}\overline Y{}^{(N)}
& \doteq &
\bigl(\overline{X}{}^{(N)}(0),\overline{E}{}^{(N)},\overline
{X}{}^{(N)},\overline{R}{}^{(N)}, \overline{\nu}{}^{(N)}_0,\nonumber\\[-8pt]\\[-8pt]
& &\hspace*{4.8pt} \overline{\nu}{}^{(N)},
\overline{\eta}^{(N)}_0,\overline{\eta}^{(N)},
\overline{A}{}^{(N)}_{\cdot,\nu}, \overline{D}{}^{(N)}_{\cdot
},\overline{A}{}^{(N)}_{\cdot,\eta},
\overline{S}{}^{(N)}_{\cdot} \bigr). \nonumber
\end{eqnarray}
Then each $\overline Y{}^{(N)}$ is a $\mathcal Y$-valued process, where
$\mathcal
Y$ is the space
\begin{eqnarray*} \mathcal{Y}&\doteq&{\mathbb R}_+ \times(\mathcal
{D}_{{\mathbb R}
_+}[0,\infty))^3 \times\mathcal{M}_F[0,H^s)
\times\mathcal D_{\mathcal{M}_F[0,H^s)}[0,\infty) \times\mathcal
{M}_F[0,H^r)\\
&&{} \times
\mathcal
D_{\mathcal{M}_F[0,H^r)}[0,\infty) \times\bigl(\mathcal{D}_{\mathcal
{M}_F([0,H^s)\times{\mathbb R}_+)}[0,\infty)\bigr)^2 \\
&&{}\times\bigl(\mathcal
{D}_{\mathcal{M}_F([0,H^r){\mathbb R}_+)}[0,\infty)\bigr)^2
\end{eqnarray*}
equipped with the product metric. Clearly, $\mathcal Y$ is a Polish space.
Now we state
the main result of this section.
\begin{theorem} \label{th-tight} Suppose Assumption \ref{ass-init} is
satisfied. Then the sequence $\{\overline Y{}^{(N)}\}$ defined in (\ref
{Y}) is
relatively compact in the Polish space $\mathcal{Y}$, and is therefore tight.
\end{theorem}

The relative compactness of $\{\overline{Y}{}^{(N)}\}$ follows
from Assumption \ref{ass-init} and Lemmas \ref{lem:rc}, \ref{lem:fnueta},
\ref{lem-tight1} and \ref{lem-tight2} below.
Since $\mathcal{Y}$ is a Polish space, tightness is then a direct
consequence of Prohorov's theorem.

We start by recalling Kurtz's criteria (see Theorem 3.8.6 of \cite
{ethkurbook} for details)
for the relative compactness of a sequence $\{\overline{F}{}^{(N)}\}$ of
processes in $\mathcal{D}_{{\mathbb R}_+}[0,\infty)$.

\begin{prop}[(Kurtz's criteria)]\label{Kurtz}
The sequence of processes $\{\overline{\mathcal Z}{}^{(N)}\}$ is
relatively compact if and only if the following two properties hold:
\begin{enumerate}[K2.]
\item[K1.] For every rational $t\geq0$,
\[
\lim_{R\rightarrow\infty}\sup_N\mathbb{P}\bigl(\overline{\mathcal
Z}{}^{(N)}(t)>R\bigr)=0.
\]
\item[K2.] For each $t>0$, there exists $\beta>0$ such that
%
%
\begin{equation}
\label{newk2}
\lim_{\delta\rightarrow
0}\sup_N\mathbb{E}\bigl[ \bigl|\overline{\mathcal Z}{}^{(N)}(t+\delta
)-\overline{\mathcal Z}{}^{(N)}(t) \bigr|^\beta\bigr]=0.
\end{equation}
\end{enumerate}
\end{prop}
\begin{lemma} \label{lem:rc} Suppose Assumption \ref{ass-init} holds.
Then the sequences $\{\overline{X}{}^{(N)}\}$,
$\{\overline{K}{}^{(N)}\}$, $ \{\overline{R}{}^{(N)}\}$,
$\{\langle{\mathbf{1}},\overline{\nu}{}^{(N)}\rangle\}, \{\langle
{\mathbf{1}},\overline{\eta}^{(N)}
\rangle\}$, the
sequences $\{\overline D{}^{(N)}_\varphi\}, \{\overline
{A}{}^{(N)}_{\varphi,\nu}\}$, for every
$\varphi\in\mathcal{C}_b([0,H^s)\times{\mathbb R}_+)$ and the
sequences $\{S^{(N)}_\psi\}, \{\overline{A}{}^{(N)}_{\psi
,\eta}\}$, for every $\psi\in\mathcal{C}_b([0, H^r)\times{\mathbb
R}_+)$, are relatively compact.
\end{lemma}
\begin{pf}
Fix $T\in(0,\infty)$. It follows from Proposition \ref
{cor-compensatormeasn}(1), (\ref{bd-1}) and (\ref{init-bd}) that for
$\varphi\in\mathcal{C}_b([0,H^s)\times{\mathbb R}_+)$,
\[
\sup_{N}\mathbb{E}\bigl[\overline{A}{}^{(N)}_{\varphi,\nu}(T) \bigr]= \sup
_{N}\mathbb{E}
\bigl[\overline D{}^{(N)}_\varphi(T)\bigr] \leq\Vert\varphi\Vert_\infty\sup
_N\mathbb{E}\bigl[\overline{X}{}^{(N)}(0)+\overline{E}{}^{(N)}(T)\bigr]
<\infty.
\]
Similarly, by Proposition \ref{cor-compensatormeasn}(2), (\ref{bd-2})
and (\ref{init-bd}), we have for every $\psi\in\mathcal
{C}_b([0,\break H^r)\times{\mathbb R}_+)$,
\[
\sup_{N}\mathbb{E}\bigl[\overline{A}{}^{(N)}_{\psi,\eta}(T)\bigr]= \sup
_{N}\mathbb{E}\bigl[\overline S^{(N)}_\psi(T)\bigr]
\leq
\Vert\psi\Vert_\infty\sup_N\mathbb{E}\bigl[\overline
{X}^{(N)}(0)+\overline{E}{}^{(N)}(T)\bigr] <\infty,
\]
which verifies condition K1 for $\mathcal Z = A^{(N)}_{\varphi,\nu},
D^{(N)}_{\varphi},
\varphi\in\mathcal{C}_b([0,H^s) \times{\mathbb R}_+)$ and
$\mathcal Z =
A^{(N)}
_{\psi,\eta},
S^{(N)}_{\psi}, \psi\in\mathcal{C}_b([0,H^r) \times{\mathbb R}_+)$.
The same argument that was used to prove
Lemma 5.8(2) in \cite{kasram07} can then be used to
show that (\ref{newk2}) is also satisfied by the same collection of
$\mathcal Z$ (see Remark \ref{rem-compen}).
The fact that $\overline{R}{}^{(N)}$ and its increments are dominated,
respectively,
by $\overline{S}{}^{(N)}$ and its increments shows that the sequence $\{
\overline{R}{}^{(N)}\}$ also
satisfies conditions K1 and K2, and is thus relatively compact.
Since $\overline{D}{}^{(N)}=\overline D{}^{(N)}_{\mathbf{1}}$ and
$\overline{S}{}^{(N)}=\overline S{}^{(N)}_{\mathbf{1}}$, it follows that
the sequences
$\{\overline{D}{}^{(N)}\}$ and $\{\overline{S}{}^{(N)}\}$ are also
relatively compact. By Assumption
\ref{ass-init}, the sequences $\{\overline{E}{}^{(N)}\}$ and $\{
\overline{X}{}^{(N)}(0)\}$ are relatively
compact.

Since for every $t\geq0$, $\langle{\mathbf{1}},\overline{\nu
}{}^{(N)}_t\rangle\leq\overline{X}{}^{(N)}
(t)\leq\overline{X}{}^{(N)}(0)+\overline{E}{}^{(N)}(t)$ by (\ref
{comp-prelimit}) and (\ref
{def-dn}), it follows from Markov's inequality that
$\langle\mathbf{1},\overline{\nu}{}^{(N)}_t\rangle$ and $\overline
{X}{}^{(N)}$ satisfy K1 of Proposition \ref{Kurtz}.
In addition, (\ref{def-dn}) also shows that
\begin{eqnarray*}
\bigl|\overline{X}{}^{(N)}(t)-\overline{X}{}^{(N)}(s) \bigr| &\leq& \bigl|\overline
{E}{}^{(N)}(t)-\overline{E}{}^{(N)}(s)
\bigr|+ \bigl|\overline{D}{}^{(N)}(t)-\overline{D}{}^{(N)}(s) \bigr|\\
& &{} + \bigl|\overline
{R}{}^{(N)}(t)-\overline{R}{}^{(N)}(s) \bigr|,
\end{eqnarray*}
and by (\ref{comp-prelimit}) and the Lipschitz continuity of the
function $x\mapsto[1-x]^+$ with Lipschitz constant $1$, we have
\[
\bigl|\bigl\langle{\mathbf{1}},\overline{\nu}^{(N)}_t\bigr\rangle-\bigl\langle
{\mathbf{1}},\overline{\nu}{}^{(N)}
_s\bigr\rangle
\bigr|= \bigl|\bigl[1-\overline{X}{}^{(N)}(t)\bigr]^+-\bigl[1-\overline{X}{}^{(N)}(s)\bigr]^+ \bigr|\leq
\bigl|\overline{X}{}^{(N)}(t)-\overline{X}{}^{(N)}
(s) \bigr|.
\]
When combined with the properties of $\overline{E}{}^{(N)}$, $\overline
{D}{}^{(N)}$ and $\overline{R}{}^{(N)}$ established
above, this shows that $\{\overline{X}{}^{(N)}\}$ and $\{\langle
{\mathbf{1}},\overline{\nu}{}^{(N)}
\rangle\}$
satisfy K2
of Proposition \ref{Kurtz} and, are relatively compact. In turn, by
(\ref{mass-queue}), the relative compactness of $\{\overline
{D}{}^{(N)}\}$ and $\{
\langle
{\mathbf{1}}, \overline{\nu}{}^{(N)}\rangle\}$ implies that of $\{
\overline{K}{}^{(N)}\}$. Moreover,
due to
(\ref{def-sn}), for every $s, t \in[0,\infty)$, we have that
%
%
\begin{eqnarray}\quad \bigl|\bigl\langle{\mathbf{1}},\overline{\eta
}^{(N)}_t\bigr\rangle-\bigl\langle
\mathbf{1},\overline{\eta}^{(N)}_s\bigr\rangle\bigr|& \leq& \bigl|\overline
{E}{}^{(N)}(t)-\overline{E}{}^{(N)}(s)
\bigr|+ \bigl|\overline{S}{}^{(N)}(t)-\overline{S}{}^{(N)}(s) \bigr|, \\
\label{dis:etaxe}
\bigl\langle{\mathbf{1}},\overline{\eta}{}^{(N)}
_t\bigr\rangle
&\leq& \bigl\langle{\mathbf{1}},\overline{\eta}^{(N)}_0\bigr\rangle
+\overline{E}{}^{(N)}(t).
\end{eqnarray}
Thus $\langle{\mathbf{1}},\overline{\eta}^{(N)}\rangle$ is also
relatively compact,
and the
proof is complete.
\end{pf}
\begin{lemma} \label{lem:fnueta}
Suppose Assumption \ref{ass-init} holds.
For every $f\in\mathcal{C}^1_c({\mathbb R}_+)$, the sequences $\{
\langle f,\overline{\nu}{}^{(N)}\rangle\}$ and
$\{\langle f,\overline{\eta}^{(N)}\rangle\}$ of $\mathcal
{D}_{{\mathbb R}}[0,\infty)$-valued random variables are
relatively compact.
\end{lemma}
\begin{pf}
Fix $t\in[0,\infty)$. By (\ref{eqn-prelimit1}) and (\ref
{eqn-prelimit3}), for every $f \in\mathcal{C}_c^1({\mathbb R}_+)$, we have
\[
\bigl\langle f, \overline{\nu}{}^{(N)}_t \bigr\rangle- \bigl\langle f, \overline{\nu
}^{(N)}_0\bigr\rangle
= \int_0^t \bigl\langle f^\prime, \overline{\nu}{}^{(N)}_s \bigr\rangle \,ds -
\overline{D}{}^{(N)}_f (t) + f(0) \overline{K}{}^{(N)}(t)
\]
and
\[
\bigl\langle f,\overline{\eta}^{(N)}_t\bigr\rangle-\bigl\langle f,\overline{\eta
}^{(N)}_0\bigr\rangle= \int_0^t\bigl\langle
f',\overline{\eta}^{(N)}_s\bigr\rangle \,ds -\overline
{S}{}^{(N)}_{f}(t)+f(0)\overline{E}{}^{(N)}(t).
\]
Since $\{\overline{D}_f^{(N)}\}$, $\{\overline{K}{}^{(N)}\}$,
$\{\overline{S}_f^{(N)}\}$ and
$\{ \overline{E}{}^{(N)}\}$ are relatively
compact due to Lem\-ma~\ref{lem:rc} and property 1 of Assumption \ref{ass-init},
it suffices to show that the sequences
$\{\int_0^\cdot\langle f^\prime, \overline{\nu}{}^{(N)}_s \rangle \,ds
\}$ and
$\{\int_0^\cdot\langle f^\prime, \overline{\eta}_s \rangle \,ds \}$
are tight.
It follows from (\ref{dis:etaxe}) that for $\delta\in(0,1)$,
\begin{eqnarray*}
\biggl|\int_t^{t+\delta}\bigl\langle f',\overline{\eta}^{(N)}_s\bigr\rangle \,ds \biggr|
&\leq&
\Vert f'\Vert_\infty\int_t^{t+\delta}\bigl|\bigl\langle{\mathbf
{1}},\overline{\eta}^{(N)}_s\bigr\rangle\bigr| \,ds
\\
&\leq&
\Vert f'\Vert_\infty\delta\bigl(\bigl\langle{\mathbf{1}},\overline
{\eta}^{(N)}_0\bigr\rangle
+\overline{E}{}^{(N)}
(t+1) \bigr).
\end{eqnarray*}
Hence, we have
%
%
\begin{equation} \label{dis:preeta}
\mathbb{E}\biggl[ \biggl|\int_t^{t+\delta}\bigl\langle f',\overline{\eta
}^{(N)}_s\bigr\rangle \,ds
\biggr| \biggr]\leq\Vert f'\Vert_\infty\delta\sup_N\mathbb{E}\bigl[\bigl\langle
{\mathbf
{1}},\overline{\eta}^{(N)}
_0\bigr\rangle+\overline{E}{}^{(N)}(t+1)\bigr].
\end{equation}
For each $t\in[0,\infty)$, by (\ref{def-etan}) and Assumption \ref
{ass-init}, it follows that
%
%
\begin{equation}\label{bound-eta} \sup_N\mathbb{E}
\bigl[\bigl\langle{\mathbf{1}},\overline{\eta}^{(N)}_t\bigr\rangle\bigr] \leq\sup
_N\mathbb{E}\bigl[\bigl\langle
\mathbf{1},\overline{\eta}^{(N)}_0\bigr\rangle+\overline{E}{}^{(N)}(t) \bigr]
<\infty.
\end{equation}
Therefore, taking the limit, as $\delta\rightarrow0$, in (\ref{dis:preeta})
and using the last inequality in (\ref{bound-eta}), we have
\[
\lim_{\delta\rightarrow0} \sup_N\mathbb{E}\biggl[ \biggl|\int_t^{t+\delta
}\bigl\langle f',\overline{\eta}^{(N)}_s\bigr\rangle \,ds \biggr| \biggr]=0.
\]
Similarly, since $\langle{\mathbf{1}}, \overline{\nu}{}^{(N)}_s
\rangle\leq1$ for every
$s\in
[0,\infty)$ and $N\in{\mathbb N}$,
\[
\lim_{\delta\rightarrow0} \sup_N\mathbb{E}\biggl[ \biggl|\int_t^{t+\delta
}\bigl\langle f',\overline{\nu}{}^{(N)}_s\bigr\rangle \,ds \biggr| \biggr]\leq\lim_{\delta
\rightarrow0}\Vert f'\Vert_\infty\delta=0.
\]
Moreover, by (\ref{bound-eta}), we also have, for every $t\in
[0,\infty)$,
\[
\sup_N \mathbb{E}\biggl[ \biggl|\int_0^{t}\bigl\langle f',\overline{\eta}^{(N)}
_s\bigr\rangle \,ds \biggr| \biggr] \leq\Vert f'\Vert_\infty t \sup_N\mathbb{E}
\bigl[\bigl\langle{\mathbf{1}},\overline{\eta}^{(N)}_0\bigr\rangle+\overline
{E}{}^{(N)}(t) \bigr]<\infty.
\]
Similarly, we have
\[
\sup_N\mathbb{E}\biggl[ \biggl|\int_0^{t}\bigl\langle f',\overline{\nu}{}^{(N)}
_s\bigr\rangle \,ds \biggr| \biggr] \leq\sup_N\mathbb{E}\biggl[\int_0^{t}\bigl|\bigl\langle
f',\overline{\nu}{}^{(N)}_s\bigr\rangle\bigr| \,ds \biggr] \leq\Vert f'\Vert_\infty
t<\infty.
\]
This implies that $ \{\int_0^\cdot\langle f',\overline{\eta
}^{(N)}_s\rangle
\,ds \}$ and $ \{\int_0^\cdot\langle f',\overline{\nu}{}^{(N)}_s\rangle \,ds
\}$ both satisfy criteria K1 and K2 of Proposition \ref{Kurtz} and
hence are relatively compact. This completes the proof of the lemma.
\end{pf}

Next, we show that $\{\overline{\nu}{}^{(N)}\}$ and $\{\overline{\eta
}^{(N)}\}$ are tight, and
hence are relatively compact with respect to the topology on $\mathcal
D_{\mathcal{M}_F[0,H^s)}[0,\infty)$ and $\mathcal D_{\mathcal
{M}_F[0,H^r)}[0,\infty)$, respectively.
Since, as mentioned in Section \ref{subsub-funmeas}, $\mathcal
{M}_F[0,H^s)$ and
$\mathcal{M}_F[0,H^r)$, equipped with the topology of weak
convergence, are Polish
spaces, we can apply Jakubowski's criteria to establish the tightness
of $\{\overline{\nu}{}^{(N)}\}$ and $\{\overline{\eta}^{(N)}\}$. For
convenience, we recall
Jakubowski's criteria.
%
%
\begin{prop}[(Jakubowski)]\label{prop:JC} A sequence $\{\overline\pi
^{(N)}\}$ of $\mathcal D_{\mathcal{M}_F[0,H)}[0,\infty)$-valued
random elements
defined on $(\Omega,\mathcal F,\mathbb{P})$ is tight if and only if the
following two conditions hold:
\begin{enumerate}[J2.]
\item[J1.] For each $T>0$ and $0<\delta<1$, there are compact subsets
$\tilde C_{T,\delta}$ of $\mathcal{M}_F[0,H)$ such that
\[
\liminf_{N\rightarrow\infty} \mathbb{P}\bigl(\overline{\nu}{}^{(N)}_t\in
\tilde
C_{T,\delta} \mbox{ for all }t\in[0,T] \bigr)>1-\delta.
\]
\item[J2.] There exists a family $\mathbb{F}$ of real continuous
functions $F$
on $\mathcal{M}_F[0,H)$ that separates points in $\mathcal{M}_F[0,H)$
and is closed under
addition, and $\{\overline\pi^{(N)}\}$ is $\mathbb{F}$-weakly tight, that
is, for every $F\in\mathbb{F}$, the sequence $\{F(\overline\pi
^{(N)}),s\in
[0,\infty)\}$ is tight in $\mathcal D_{{\mathbb R}}[0,\infty)$.
\end{enumerate}
\end{prop}
\begin{lemma} \label{lem-tight1}
Suppose Assumption \ref{ass-init} holds. The sequences $\{\overline
{\nu}^{(N)}\}$ and
$\{\overline{\eta}^{(N)}\}$ are relatively compact.
\end{lemma}
\begin{pf}
By Remark 5.11 of \cite{kasram07} and Lemma \ref{lem:fnueta}, it
follows that $\{\overline{\nu}{}^{(N)}\}$ and $\{\overline{\eta
}^{(N)}\}$ satisfy Jakubowski's J2
criterion. Therefore, it suffices to show that they also satisfy Jakubowski's
J1 criterion. By (2) and (3) of Assumption \ref{ass-init}, for almost every
$\omega\in\Omega$, $\sup_N \overline{\nu}{}^{(N)}_0(\omega
)[0,H^s)<\infty$. By
Lemma A 7.5
of \cite{Kal}, for every $\varepsilon>0$, there exists
$k(\omega,\varepsilon)<\infty$ such that $\sup_N
\overline{\nu}{}^{(N)}_0(\omega)(k(\omega,\varepsilon
),\break H^s)<\varepsilon$.
The argument for tightness of $\{\overline{\nu}{}^{(N)}\}$ (in the
absence of reneging)
presented in Lemma 5.12 of \cite{kasram07} can be directly applied to show
that $\{\overline{\nu}{}^{(N)}\}$ satisfies Jakubowski's J1 criterion,
and hence
$\{\overline{\nu}{}^{(N)}\}$ is tight in the presence of reneging as well.
Similarly, due to (2) and (4) of Assumption \ref{ass-init}, for almost every
$\omega\in\Omega$, $\sup_N \overline{\eta}^{(N)}_0(\omega
)[0,H^r)<\infty$.
Once again, by
Lemma A 7.5 of \cite{Kal}, we infer that for every $\varepsilon>0$, there
exists $l(\omega,\varepsilon)<\infty$ such that $\sup_N
\overline{\eta}^{(N)}_0(\omega)(l(\omega,\varepsilon
),H^r)<\varepsilon$. Since
$\{\langle{\mathbf{1}},
\overline{\eta}^{(N)}\rangle\}$ is tight by Lemma \ref{lem:fnueta}, the argument
for tightness of
$\{\overline{\nu}{}^{(N)}\}$ presented in Lemma 5.12 of \cite
{kasram07} can also be adapted
to show that the sequence $\{\overline{\eta}^{(N)}\}$ satisfies
Jakubowski's J1
criterion, and is therefore tight. We omit the details.
\end{pf}

We end this section by
establishing the relative compactness of the measure-valued
processes associated with the cumulative departure and reneging
functionals and their compensators.
\begin{lemma}\label{lem-tight2}
Suppose Assumption \ref{ass-init} holds. Then the sequences $\{
\overline{D}{}^{(N)}
_{\cdot}\}$ and $\{\overline{A}{}^{(N)}_{\cdot,\nu}\}$ are
relatively compact
in $\mathcal{D}_{\mathcal{M}_F([0,H^s)\times{\mathbb R}_+)}[0,\infty
)$. Similarly, the sequences $\{\overline{S}{}^{(N)}_{\cdot}\}$ and $\{
\overline{A}{}^{(N)}_{\cdot,\eta}\}$ are relatively compact in
$\mathcal{D}_{\mathcal{M}_F([0,H^r)\times{\mathbb R}_+)}[0,\infty)$.
\end{lemma}
\begin{pf}
This can be proved by combining Lemma \ref{lem:rc} and Proposition
\ref{cor-compensatormeasn}
with the argument that was used in Lemma 5.13 of \cite{kasram07} to
establish the tightness of the sequences
$\{\overline\mathcal{Q}{}^{(N)}\}$ and $\{\overline\mathcal{A}{}^{(N)}\}$
therein. Since the adaptation of
the argument in \cite{kasram07} is
straightforward, we omit the details.
\end{pf}

\section{Strong law of large numbers limits}
\label{sec:CSL}

\subsection{Characterization of subsequential limits}
\label{subs-csl}

The focus of this section is the following theorem which, in particular,
establishes existence of a solution to the fluid equations.
\begin{theorem}\label{thm:FE}
Suppose that Assumptions \ref{ass-init}--\ref{ass-h} hold. Let
$(\overline{X}
,\overline{\nu},\overline{\eta})$ be the limit of any subsequence
of $\{\overline{X}{}^{(N)},\overline{\nu}{}^{(N)}
,\overline{\eta}^{(N)}\}$. Then $(\overline{X},\overline{\nu
},\overline{\eta})$ solves the fluid equations.
\end{theorem}

The rest of the section is devoted to the proof of this theorem. Let
$(\overline{E},
\overline{X}(0)$, $\overline{\nu}_0,\overline{\eta}_0)$ be the
$\mathcal{S}_0$-valued random
variable\vspace*{1pt} that
satisfies Assumption~\ref{ass-init}, and let $\{\overline Y{}^{(N)}\}
_{N\in{\mathbb N}}$
be the sequence of processes defined in (\ref{Y}). Then, by Assumption
\ref{ass-init}, Theorem \ref{th-tight} and the limits $\overline{M}{}^{(N)}
_{\cdot,\nu}=
\overline{D}{}^{(N)}_{\cdot} - \overline{A}{}^{(N)}_{\cdot,\nu
}\Rightarrow0$ and
$\overline{M}{}^{(N)}_{\cdot,\eta}= \overline{S}{}^{(N)}_{\cdot} -
\overline{A}{}^{(N)}_{\cdot,\eta
}\Rightarrow0$
established in Proposition \ref{cor-compensatormeasn},
there exist processes
$\overline{X}\in\mathcal{D}_{{\mathbb R}_+}[0,\infty), \overline
{R}\in\mathcal{D}_{{\mathbb R}
_+}[0,\infty), \overline{\nu}\in\mathcal
D_{\mathcal{M}_F[0,H^s)}[0,\infty)$, $\overline{\eta}\in\break \mathcal
D_{\mathcal{M}_F[0,H^r)}[0,\infty)$,
$\overline{A}_{\cdot,\nu}\in\mathcal{D}_{\mathcal
{M}_F([0,H^s)\times{\mathbb R}_+)}[0,\infty)$, $\overline{D}_{\cdot
}\in$ $\mathcal{D}_{\mathcal{M}_F([0,H^s)\times{\mathbb
R}_+)}[0,\infty)$,
$\overline{A}_{\cdot,\eta}\in\mathcal{D}_{\mathcal
{M}_F([0,H^r)\times{\mathbb R}_+)}[0,\infty)$, $\overline{S}_{\cdot
}\in\mathcal{D}_{\mathcal{M}_F([0,H^r)\times{\mathbb
R}_+)}[0,\infty)$
such that $\overline Y{}^{(N)}$ converges weakly (along a suitable
subsequence) to
\[
\overline Y\doteq(\overline{X}(0),\overline{E},\overline
{X},\overline{R}, \overline{\nu}_0, \overline{\nu},
\overline{\eta}_0,\overline{\eta},\overline{A}_{\cdot,\nu},
\overline{A}_{\cdot,\nu
},\overline{A}_{\cdot,\eta}, \overline{A}_{\cdot,\eta} ) \in
\mathcal{Y}.
\]
Denoting this subsequence again by $\overline Y{}^{(N)}$ and invoking the
Skorokhod representation theorem, with a slight abuse of notation, we
can assume that, $\mathbb{P}$
a.s., $\overline Y{}^{(N)}\rightarrow\overline Y$ as $N\rightarrow
\infty$.
Without loss of generality, we may further assume that the above
convergence holds everywhere.

We now identify some properties of the limit that will be used to prove
Theorem \ref{thm:FE}. From Proposition \ref{cor-compensatormeasn}(1),
it follows that, as $N\rightarrow\infty$, $(\overline Y{}^{(N)},\break
\overline{D}{}^{(N)}_{\cdot
}) \rightarrow(\overline Y,\overline{A}_{\cdot,\nu})$. Together
with (\ref
{def-dn}), this implies that
%
%
\begin{equation}\label{limX}
\overline{X}= \overline{X}(0)+\overline{E}-\overline{A}_{{\mathbf
{1}},\nu}-\overline{R}.
\end{equation}
Moreover, we claim that
%
%
\begin{equation}\label{dis:13}
\overline{A}_{\varphi,\nu}=\int_0^\cdot\langle\varphi
h^s,\overline{\nu}_s\rangle \,ds.
\end{equation}
This
corresponds to relation (5.48) established in Proposition 5.17 of \cite
{kasram07} for the
model without abandonments. However, essentially the same argument can
be used
here as well. Specifically, the proof of (5.48) in \cite{kasram07}
relies on
Lemmas 5.8(1) and~5.16 of \cite{kasram07},
which continue to be valid in the presence of abandonments due to Remarks
\ref{rem-compen} and
\ref{rem-compen0}.
On substituting (\ref{dis:13}) into (\ref{limX}), we see that the
fluid equation (\ref{eq-fx}) is satisfied.

Next, in Proposition \ref{prop:3}, we establish representation (\ref{fr})
for $\overline{R}$ given in the fluid equations. The proof of this result
relies on the alternative representation for the compensator $A^{(N)}
_{\theta^{(N)},\eta}$ of $R^{(N)}$ given in
(\ref{rep-rcomp2}).
\begin{prop}
\label{prop:3}
For every $T\in[0,\infty)$, as $N\rightarrow\infty$,
%
%
\begin{equation}\label{rep-Aconv}
\mathbb{E}\biggl[\sup_{t\in[0,T]} \biggl|\overline{A}{}^{(N)}_{\theta
^{(N)},\eta}(t) - \int_0^t \biggl(\int_0^{\overline
{Q}(s)}h^r((F^{\overline\eta_s}
)^{-1}(y))\,dy \biggr)\,ds \biggr| \biggr]\rightarrow0.
\end{equation}
Moreover, almost surely,
%
%
\begin{equation}\label{rep-R}\overline{R}(t)=\int_0^t \biggl(\int
_0^{\overline{Q}(s)}h^r((F^{\overline\eta_s})^{-1}(y))\,dy \biggr) \,ds,\qquad
t\in [0,\infty).
\end{equation}
\end{prop}

The proof of Proposition \ref{prop:3} is given near the end of this
section and relies on the following preliminary
observations.
Let $\tilde R(t)$ be defined by the right-hand side of (\ref{rep-R})
for $t\in[0,\infty)$. We first show how (\ref{rep-R}) can be deduced
from (\ref{rep-Aconv}). From (\ref{rep-Aconv}), it follows that
$\overline{A}{}^{(N)}_{\theta^{(N)},\eta} \Rightarrow\tilde R$ as
$N\rightarrow
\infty$. Since $\tilde{R}$ is continuous, $\overline
{R}{}^{(N)}=\overline{M}{}^{(N)}_{\theta
^{(N)},\eta} + \overline{A}{}^{(N)}_{\theta^{(N)},\eta}$ and
$\overline{M}{}^{(N)}
_{\theta^{(N)},\eta} \Rightarrow0$ by
Lemma \ref{cor:1}, it follows that $\overline{R}{}^{(N)} \Rightarrow
\tilde{R}$.
This implies, a.s., $\tilde{R} = \overline{R}$, and thus the second statement
of Proposition \ref{prop:3} follows from the first statement.

The proof of\vspace*{1pt} (\ref{rep-Aconv}) relies on Lemmas
\ref{lem:uni4}--\ref{lem:2} below and the following observations.
Using (\ref{rep-rcomp2}) and the elementary relation $(F^{\eta
^{(N)}_s})^{-1}(N\cdot)=(F^{\overline\eta^{(N)}_s})^{-1}(\cdot)$,
simple algebraic manipulations show
that
%
%
\begin{eqnarray}\label{rep-rcomp3} \overline{A}{}^{(N)}_{\theta
^{(N)},\eta} (t)
\doteq\int_0^t \biggl(\int_0^{\overline{Q}{}^{(N)}(t)+\overline\iota
^{(N)}(t)}h^r((F^{\overline\eta^{(N)}_s})^{-1}(y))\,dy
\biggr)\,ds,\nonumber\\[-8pt]\\[-8pt]
\eqntext{t\in[0,\infty),}
\end{eqnarray}
where, as usual, $\overline\iota
^{(N)}\doteq\iota^{(N)}/N$ and $\iota^{(N)}$ is given by (\ref
{iota}). Next, observe that for all $t\in[0,T]$ and $L\in[0,H^r)$,
%
%
\begin{equation}\label{dcomp1} \bigl|\overline{A}{}^{(N)}_{\theta
^{(N)},\eta
}(t)-\tilde
R(t) \bigr|\leq\overline C{}^{(N)}_1(t,L)+\overline C{}^{(N)}_2(t,L)+\overline C_3(t,L),
\end{equation}
where $\overline C{}^{(N)}_i(t,L), i=1,2$, and $\overline C_3(t,L)$ are defined, for
$t\in[0,\infty)$, by
%
%
\begin{eqnarray} \label{def-c1} \qquad\overline C{}^{(N)}_1(t,L)
&\doteq& \biggl|\int_0^t \biggl(\int_0^{(\overline{Q}{}^{(N)}(s)+\overline\iota
^{(N)}(s))\wedge F^{\overline\eta_s^{(N)}}(L)} h^r((F^{\overline\eta
_s^{(N)}})^{-1}(y))\,dy \biggr)\,ds\\
& &\hspace*{52.4pt}{} - \int_0^t \biggl(\int_0^{\overline{Q}
(s)\wedge F^{\overline\eta_s}(L)}h^r((F^{\overline\eta_s})^{-1}(y))\,dy
\biggr)\,ds \biggr|, \nonumber
\\
%
%
\label{def-c2} \overline C{}^{(N)}_2(t,L) &\doteq&
\biggl|\int
_0^t \biggl( \int_{(\overline{Q}{}^{(N)}(s)+\overline\iota^{(N)}(s))\wedge
F^{\overline\eta_s^{(N)}}
(L)}^{\overline{Q}{}^{(N)}(s)+\overline\iota
^{(N)}(s)}h^r((F^{\overline\eta_s^{(N)}})^{-1}(y))\,dy
\biggr)\,ds \biggr|
\end{eqnarray}
and
%
%
\begin{equation}\label{def-c3} \overline C_3(t,L) \doteq\int
_0^t \biggl( \int
_{\overline{Q}(s)\wedge F^{\overline\eta_s}(L)}^{\overline
{Q}(s)}h^r((F^{\overline\eta_s})^{-1}(y))\,dy
\biggr)\,ds.
\end{equation}

As a precursor to the proof of (\ref{rep-Aconv}) of Proposition \ref
{prop:3}, we first establish
some path properties of the limiting queue measure $\overline{\eta}$
in Lemma\ref{lem:uni4}
and some estimates in Lemma \ref{lem-est}. These two preliminary
results will be used in Lemma \ref{lem:1} to show that for any $L\in
[0,H^r)$, $\lim_{N\rightarrow\infty} \sup_{t\in[0,T]}
|\overline C{}^{(N)}_1(t,L) |=0$ in the case when $h^r$ is
continuous. Next, Lemma \ref{lem:2} extends this to include general
$h^r$ that is locally integrable in $[0,H^r)$. All these results are
then combined to prove Proposition \ref{prop:3}.
\begin{lemma} \label{lem:uni4}
For every $L\in[0,H^r)$, $\overline{\eta}_t$ is continuous at $L$
for almost
every \mbox{$t\geq0$}. Moreover, for $t\in(0,\infty)$ and $L\in[0,H^r)$,
if $\overline{\eta}_t(\{L\})>0$, then $\overline{\eta
}_t(L,L+\varepsilon)>0$ for all
sufficiently small $\varepsilon$.
\end{lemma}
\begin{pf}
It was shown in Corollary \ref{cor:nueta} that $(\overline{\eta
},\overline{E})$
satisfies (\ref{eq-freneg2}) for every bounded Borel measurable
function $f$.
For every $L \in[0,H^r)$, substituting $f = \mathbh{1}_{L}$ in (\ref
{cor:nueta}),
we obtain
%
%
\begin{eqnarray} \label{atL} \overline{\eta}_t (\{L\}) & = &\int
_{[0,H^r )}
\mathbh{1}_{\{L\}} (x+t)
\frac{1 - G^r(x+t)}{1 - G^r(x)} \overline{\eta}_0
(dx)\nonumber\\[-8pt]\\[-8pt]
& &{} + \int
_{[0,t]} \mathbh{1}_{\{L\}} (t-s) \bigl(1
- G^r(t-s)\bigr) \,d \overline{E}(s).\nonumber
\end{eqnarray}
It is easy to see that the right-hand side of the above display is zero
except when $\overline{\eta}_0(\{L-t\})>0$ if $t\leq L$ or when
$\overline{E}(t-L)-\overline{E}
((t-L)-)>0$ if $t>L$. Since the jump times of both $\overline{\eta
}_0$ and
$\overline{E}$ are at most countable, (\ref{atL}) shows that
$\overline{\eta}_t$ is
continuous at $L$ for almost every $t\geq0$.

Next, suppose $\overline{\eta}_t(\{L\})>0$. Then by (\ref{atL}), at
least one
of the following two inequalities must hold:
%
%
\begin{equation}\label{atL1}\int
_{[0,H^r )} \mathbh{1}_{\{L\}} (x+t)
\frac{1 - G^r(x+t)}{1 - G^r(x)} \overline{\eta}_0 (dx)>0
\end{equation}
or
%
%
\begin{equation}\label{atL2}\int_{[0,t]} \mathbh{1}_{\{L\}} (t-s) \bigl(1-
G^r(t-s)\bigr)\, d
\overline{E}(s)>0.
\end{equation}
If (\ref{atL1}) holds, then it must be that $L-t\in[0,H^r)$, $(1 -
G^r(L))/(1 - G^r(L-t))>0$ and $\overline{\eta}_0(\{L-t\})>0$. By Assumption
\ref{ass-jump} and the continuity of $(1 - G^r(\cdot+t))/(1 -
G^r(\cdot))$, it then follows that for all sufficient small
$\varepsilon>0$,
%
%
\begin{equation}\label{dis:lle}\int_{[0,H^r )} \mathbh{1}
_{(L,L+\varepsilon)} (x+t)
\frac{1 - G^r(x+t)}{1 - G^r(x)} \overline{\eta}_0 (dx)>0.
\end{equation}
Substituting
$f=\mathbh{1}_{(L,L+\varepsilon)}$ into (\ref{eq-freneg2}) in Corollary
4.2 shows that $\overline{\eta}_t(L,L+\varepsilon)$ is greater than
or equal
to the left-hand side of (\ref{dis:lle}), and so the lemma is
established in this case. On the other hand, suppose (\ref{atL2})
holds. In this case, $t-L>0$, $1-G^r(t-L)>0$ and $\overline
{E}(t-L)-\overline{E}
((t-L)-)>0$. By Assumption \ref{ass-jump} and the continuity of
$1-G^r(t-\cdot)$, for all sufficiently small $\varepsilon>0$,
$1-G^r(t-\cdot)$ is strictly positive on $(L,L+\varepsilon)$ and
$\overline{E}
((t-L)-)-\overline{E}(t-L-\varepsilon)>0$. Another application of
(\ref
{eq-freneg2}) of Corollary \ref{cor:nueta}, with $f=\mathbh{1}
_{(L,L+\varepsilon)}$, shows that
\[
\overline{\eta}_t(L,L+\varepsilon) \geq\int_0^t \mathbh
{1}_{(L,L+\varepsilon)}
(t-s) \bigl(1
- G^r(t-s)\bigr) \,d \overline{E}(s)>0,
\]
and the proof of the lemma is complete.
\end{pf}
\begin{lemma}
\label{lem-est} Let $T\in[0,\infty)$ and $L\in[0,H^r)$. The
following estimates hold:
\begin{enumerate}
\item For $m \in[0, H^r)$ and every $\ell\in L^1_{\mathrm{loc}}[0, H^r)$ with
support in $[0,m]$, there exists $\tilde L(m,T) < \infty$ such that
%
%
\begin{equation}\label{lem-est-1} \biggl|\int_0^T\langle\ell, \overline
{\eta}
_s\rangle
\,ds \biggr|\leq
{\tilde L(m,T) \int_{[0,H^r)}}|\ell(x)|\,dx.
\end{equation}

\item Suppose $h$ is a measurable function such that $\tilde
C_L^h\doteq\sup_{x\in[0,L]}|h(x)|<\infty$. Then, $\mathbb{P}$-a.s.,
%
%
\begin{equation}
\label{lem-est-2}\sup_N\sup_{s\in[0,T]}\int_0^Lh(x)\overline{\eta}^{(N)}
_s(dx)\leq\tilde C_L^h \sup_N \bigl(\bigl\langle{\mathbf{1}},\overline{\eta}^{(N)}
_0\bigr\rangle
+\overline{E}{}^{(N)}(T) \bigr)<\infty.
\end{equation}
\end{enumerate}
\end{lemma}
\begin{pf}
It was established in Lemma 5.16 of \cite{kasram07} that inequality
(\ref{lem-est-1}) holds with $\overline{\eta}$ replaced by the fluid age
measure $\overline{\nu}$
associated with a many-server queue without abandonments. The proof follows
directly from Proposition 4.15 and the estimate (5.46) of
\cite{kasram07}. Since the dynamic equations (\ref{eqn-prelimit3}) and
(\ref{eq-freneg2}) for $\eta^{(N)}$ and~$\overline{\eta}$,
respectively, are exactly
analogous to the dynamic equations for $\nu^{(N)}$ and~$\overline{\nu}$.
Estimate (5.46) of \cite{kasram07} can be shown to hold for $\overline
{\eta}$
using the same
argument as in~\cite{kasram07}. When combined with Proposition 4.15 of
\cite{kasram07}, this shows that (\ref{lem-est-1}) holds.
Estimate (\ref{lem-est-2}) follows directly from (\ref{def-sn}) and
Assumption \ref{ass-init}.
\end{pf}
\begin{lemma}\label{lem:1} For $T\geq0$ and all but countably many
$L\in[0,H^r)$, given any continuous function $h$ on $[0,\infty)$, as
$N\rightarrow\infty$, for every realization,
%
%
\begin{eqnarray} \label{dis:ext2}\qquad& & \sup_{t\in[0,T]} \biggl|\int
_0^t \biggl(\int_0^{(\overline{Q}{}^{(N)}(s)+\overline\iota^{(N)}(s))\wedge
F^{\overline\eta_s^{(N)}}
(L)} h((F^{\overline\eta_s^{(N)}})^{-1}(y))\,dy \biggr)\,ds
\nonumber\\[-8pt]\\[-8pt]
& &\hspace*{80pt}{}
- \int_0^t \biggl(\int_0^{\overline{Q}(s)\wedge F^{\overline\eta_s}
(L)}h((F^{\overline\eta_s})^{-1}(y))\,dy \biggr)\,ds \biggr|\rightarrow
0.\nonumber
\end{eqnarray}
\end{lemma}
\begin{pf}
Fix $\omega\in\Omega$. To ease the notation, we shall suppress
$\omega$ from the notation.
From the convergence of $\overline{\eta}^{(N)}$ to $\overline{\eta
}$ and $\overline{Q}{}^{(N)}$ to $\overline{Q}
$, it follows that, as $N\rightarrow\infty$, $\overline{\eta
}^{(N)}_s \stackrel{w}{\rightarrow}
\overline{\eta}_s$ and $\overline{Q}{}^{(N)}(s) \rightarrow\overline
{Q}(s)$ for almost every $s\geq
0$. Also, by Lem\-ma~\ref{lem:uni4}, $\overline{\eta}_s$ is continuous
at $L$
for almost every $s\geq0$. Let $s\geq0$ be a time at which $\overline
{\eta}^{(N)}
_s \stackrel{w}{\rightarrow}\overline{\eta}_s$ and $\overline
{Q}^{(N)}(s) \rightarrow\overline{Q}(s)$ as
$N\rightarrow\infty$ and $\overline{\eta}_s$ is continuous at $L$.
Then, as
$N\rightarrow\infty$, $F^{\overline\eta_s^{(N)}}(x)\rightarrow
F^{\overline\eta_s}(x)$ for $x=L$ and all but
a countable number of $x\in[0,H^r)$. Therefore, by Theorem 13.6.3 of
\cite{whi-SPL}, we have $(F^{\overline\eta_s^{(N)}})^{-1}
\rightarrow(F^{\overline\eta_s})^{-1}$ on
$[0,F^{\overline\eta_s}(H^r-))$ in the $M_1$ topology. For $s\in
[0,T]$, we now show
that, as $N\rightarrow\infty$,
%
%
\begin{eqnarray}\label{dis:4}
&&
\int_0^{(\overline{Q}{}^{(N)}
(s)+\overline\iota^{(N)}(s))\wedge F^{\overline\eta
_s^{(N)}}(L)}h((F^{\overline\eta_s^{(N)}}
)^{-1}(y))\,dy\nonumber\\[-8pt]\\[-8pt]
&&\qquad\rightarrow\int_0^{\overline{Q}(s)\wedge F^{\overline
\eta_s}
(L)}h((F^{\overline\eta_s})^{-1}(y))\,dy.\nonumber
\end{eqnarray}
From the inequality $ |\overline\iota^{(N)} |\leq1/N$, we
immediately see that
%
%
\begin{equation} \label{dis:Mconv}\qquad
\bigl(\overline{Q}{}^{(N)}(s)+\overline\iota^{(N)}(s)\bigr)\wedge
F^{\overline\eta_s^{(N)}}(L) \rightarrow\overline{Q}(s)\wedge
F^{\overline\eta_s}(L)\qquad\mbox{as
}N\rightarrow\infty.
\end{equation}
We now consider the
following two cases:

\textit{Case} 1. $\overline{Q}(s)\wedge F^{\overline\eta
_s}(L)<F^{\overline\eta_s}(H^r-)$. In this
case, due
to (\ref{dis:Mconv}), for all sufficiently large $N$, $(\overline{Q}{}^{(N)}
(s)+\overline
\iota^{(N)}(s))\wedge F^{\overline\eta_s^{(N)}}(L) <F^{\overline
\eta_s}(H^r-)$. For each $n\in{\mathbb N}$,
by Theorem
11.5.1 of \cite{whi-SPL} and the continuity of $h$, we obtain for each
$t<F^{\overline\eta_s}(H^r-)$,
\[
\lim_{N \rightarrow\infty}\sup_{u\in
[0,t]} \biggl|\int_0^{u}h((F^{\overline\eta_s^{(N)}})^{-1}(y))\,dy-\int
_0^{u}h((F^{\overline
\eta_s})^{-1}(y))\,dy \biggr| = 0.
\]
By the case assumption, this implies, in particular, that
\[
\lim_{N \rightarrow
\infty} \biggl|\int_0^{\overline{Q}(s)\wedge F^{\overline\eta
_s}(L)}h((F^{\overline\eta_s^{(N)}})^{-1}(y))\,dy -
\int_0^{\overline{Q}(s)\wedge F^{\overline\eta
_s}(L)}h((F^{\overline\eta
_s})^{-1}(y))\,dy \biggr| =
0.
\]
On the other hand, (\ref{dis:Mconv}) and the continuity of $h$ show that
\[
\lim_{N \rightarrow\infty}\int_{(\overline{Q}{}^{(N)}(s)+\overline
\iota^{(N)}(s))\wedge
F^{\overline\eta_s^{(N)}}(L)}^{\overline{Q}(s)\wedge F^{\overline
\eta_s}(L)}h((F^{\overline\eta_s^{(N)}})^{-1}(y))\,dy = 0.
\]
Together, the last two assertions imply (\ref{dis:4}).

\textit{Case} 2. $\overline{Q}(s)\wedge F^{\overline\eta
_s}(L)=F^{\overline\eta_s}(H^r-)$. We
first claim that in this case
%
%
\begin{equation}\label{dis:qlh} \overline{Q}(s)=F^{\overline\eta_s}
(L)=F^{\overline\eta_s}(H^r-).
\end{equation}
Indeed, $F^{\overline\eta_s}(L)\leq F^{\overline\eta_s}(H^r-)$ because
$F^{\overline\eta_s}$ is nondecreasing and $L<H^r$, while $\overline
{Q}(s)\leq\overline{\eta}
_s[0,H^r)=F^{\overline\eta_s}(H^r-)$ by (\ref{fqfreneg}). On the
other hand, the
reverse inequalities $\overline{Q}(s)\geq F^{\overline\eta_s}(H^r-)$
and $F^{\overline\eta_s}(L)\geq
F^{\overline\eta_s}(H^r-)$ hold by the case assumption, and so the
claim follows.
Now, define $\overline{L}\doteq(F^{\overline\eta_s})^{-1}(F^{\overline
\eta_s}(H^r-))$.
Then $\overline L=( F^{\overline\eta_s})^{-1}(F^{\overline\eta
_s}(L))$ by (\ref
{dis:qlh}). Hence, $\overline L\leq L$ and
%
%
\begin{equation}\label{feta-eq} F^{\overline\eta_s}
(\overline L)=F^{\overline\eta_s}(L)=F^{\overline\eta_s}(H^r-).
\end{equation}
This implies $\overline{\eta}_s(\overline
L,H^r)=0$, and from the second assertion of Lemma \ref{lem:uni4}, it
follows that
%
%
\begin{equation}\label{dis:reneg0}
\overline{\eta}_s(\{\overline L\})=0.
\end{equation}
The change of variables formula and (\ref{feta-eq}) then yield
%
%
\begin{eqnarray}
\label{eq-cvf1} \int_0^{\overline{Q}(s)\wedge F^{\overline\eta
_s}(L)}h((F^{\overline\eta
_s})^{-1}(y))\,dy&=&\int_{[0,H^r)}h(x)\overline{\eta}_s(dx)\nonumber\\[-8pt]\\[-8pt]
&=&\int
_{[0,\overline
L]}h(x)\overline{\eta}_s(dx).\nonumber
\end{eqnarray}
Also, by Proposition \ref{lem:uni5} and
another application of the change of variables formula, we have
%
%
\begin{eqnarray}
\label{eq-cvf2}
&&
\int_0^{(\overline{Q}{}^{(N)}(s)+\overline\iota^{(N)}(s))\wedge
F^{\overline\eta_s^{(N)}}(L)}h((F^{\overline\eta
_s^{(N)}})^{-1}(y))\,dy \nonumber\\[-8pt]\\[-8pt]
&&\qquad=\int_{[0,\chi
^{(N)}(s-)]}\mathbh{1}_{[0,L]}(x)h(x)\overline{\eta
}^{(N)}_s(dx).\nonumber
\end{eqnarray}
Expanding the term on the right-hand side of (\ref{eq-cvf2}) and using
the inequality \mbox{$\overline L \leq L$}, we obtain
%
%
\begin{eqnarray}
\label{dis:3}
& &\int_{[0,\chi^{(N)}(s-)]}\mathbh{1}
_{[0,L]}(x)h(x)\overline{\eta}^{(N)}_s(dx) \nonumber\\
&&\qquad=\int_{[0,\overline
L]}\mathbh{1}
_{[0,L]}(x)h(x)\overline{\eta}^{(N)}_s(dx)\nonumber\\[-8pt]\\[-8pt]
&&\qquad\quad{}+\int_{(\chi
^{(N)}(s-)\wedge\overline
L,\chi^{(N)}(s-)]}\mathbh{1}_{[0,L]}(x)h(x)\overline{\eta
}^{(N)}_s(dx) \nonumber\\
&&\qquad\quad{} - \int_{(\chi^{(N)}(s-)\wedge\overline L,\overline L]}\mathbh{1}
_{[0,L]}(x)h(x)\overline{\eta}^{(N)}_s(dx). \nonumber
\end{eqnarray}
By (\ref{eq-cvf1}) and (\ref{eq-cvf2}), the left-hand side and the
first term
on the right-hand side of (\ref{dis:3}), respectively, equal the left-hand
side and right-hand side of (\ref{dis:4}). Therefore, to prove (\ref{dis:4})
it suffices to show that the second and the third terms on the
right-hand side
of (\ref{dis:3}) converge to zero, as $N\rightarrow\infty$. Recall the
constant $\tilde C_L^h$ defined in Lemma \ref{lem-est}. Note that
$\tilde
C_L^h<\infty$ since $h$ is continuous. Therefore, the second term on the
right-hand side of (\ref{dis:3}) is bounded above by $\tilde C_L^h
\overline{\eta}^{(N)}_s(\chi^{(N)}(s-)\wedge\overline L,\chi
^{(N)}(s-)]$. By
(\ref{dis:reneg0}), Portmanteau's theorem and (\ref{feta-eq}), it
follows that
\[
\lim_{N\rightarrow\infty}\overline{\eta}^{(N)}_s\bigl(\chi
^{(N)}(s-)\wedge\overline
L,\chi^{(N)}(s-)\bigr] \leq\lim_{N\rightarrow\infty}\overline{\eta
}^{(N)}_s(\overline
L,H^r)=\overline{\eta}[\overline L,H^r)=0.
\]
On the other hand, the absolute value of the third term on the
right-hand side
of (\ref{dis:3}) is bounded above by $\tilde C_L^h
\overline{\eta}^{(N)}_s(\chi^{(N)}(s-)\wedge\overline L,\overline L]$. We now
argue by contradiction
to show that $\liminf_{N\rightarrow\infty}\chi^{(N)}(s-)\geq\overline
L$ and,
consequently, that $\overline{\eta}^{(N)}_s(\chi^{(N)}(s-)\wedge
\overline L,\overline L]$
converges to zero as $N\rightarrow\infty$. Indeed, suppose this
assertion were false. Then there must exist a subsequence $\{N_k\}
_{k\in{\mathbb N}}$ such that $\lim_{k\rightarrow\infty}\chi^{(N_k)}(s-)=
\overline L-\delta$ for some $\delta>0$. Hence, for $k$ large enough,
$\chi^{(N_k)}(s-)<\overline L-\delta/2$. By Lemma \ref{lem-chi}, we have
$\chi^{(N_k)}(s-)\geq\chi^{(N_k)}(s)$. Hence $\overline{\eta}
_s^{(N_k)}[0,\overline L-\delta/2]\geq\overline{Q}{}^{(N_k)}(s)$ by (\ref{qn}).
Sending $k\rightarrow\infty$ and using the convergence $\overline
{\eta}
_s^{(N_k)}\Rightarrow\overline{\eta}_s$, the fact that $[0,\overline
L-\delta/2]$
is closed and Portmanteau's theorem, we obtain $\overline{\eta
}_s[0,\overline
L-\delta/2]\geq\overline{Q}(s)$. This contradicts the definition of
$\overline L$,
and hence completes the proof of (\ref{dis:4}).

Finally, we deduce (\ref{dis:ext2}) from (\ref{dis:4}) using the
bounded convergence theorem, whose application is justified by the
bounds (\ref{eq-cvf1}), (\ref{eq-cvf2}) and the estimate (\ref{lem-est-2}).
\end{pf}

We now generalize Lemma \ref{lem:1} to allow for a general locally
integrable (not necessarily continuous) function $h^r$ on $[0,H^r)$.
\begin{lemma} \label{lem:2}
Let $L<H^r$, and let $\overline C{}^{(N)}_1(t,L), t\in[0,\infty),
N\in{\mathbb N}$ be defined as in (\ref{def-c1}). Then for every
$T\in
[0,\infty)$, almost surely for $L<H^r$,
%
%
\begin{equation}\label{dis:5}\lim_{N\rightarrow
\infty}\sup_{t\in[0,T]}\overline C{}^{(N)}_1(t,L)=0.
\end{equation}
\end{lemma}
\begin{pf}
Fix $L<H^r$. Since $h^r$ lies in $\mathcal{L}^1_{\mathrm{loc}}[0,H^r)$ and is
nonnegative,
there exists a sequence of nonnegative continuous functions $\{h^r_n\}
_{n\geq
1}$ on $[0,H^r)$ such that $\int_0^L|h^r(x)-h^r_n(x)|\,dx\rightarrow0$ as
$n\rightarrow\infty$ and $h^r_n$ has common compact support in
$[0,H^r)$. For
each $n\in{\mathbb N}$, (\ref{dis:5}) holds with $h^r_n$ in place of $h^r$
due to
Lemma \ref{lem:1}.
Let $l^r_n=|h^r_n-h^r|$ for each $n\geq1$.
Then, in order to prove (\ref{dis:5}), it clearly suffices to show
that the following two limits hold: almost everywhere,
%
%
\begin{equation}\label
{dis:11}\lim_{N\rightarrow\infty}\sup_N\int_0^T \biggl(\int
_0^{(\overline{Q}{}^{(N)}
(s)+\overline\iota^{(N)}(s))\wedge F^{\overline\eta
_s^{(N)}}(L)}l^r_n((F^{\overline\eta_s^{(N)}}
)^{-1}(y))\,dy \biggr)\,ds= 0\hspace*{-28pt}
\end{equation}
and
%
%
\begin{equation}\label{dis:12}\lim_{N\rightarrow\infty
}\int_0^T \biggl(\int_0^{\overline{Q}(s)\wedge F^{\overline\eta
_s}(L)}l^r_n((F^{\overline\eta_s}
)^{-1}(y))\,dy \biggr)\,ds= 0.
\end{equation}

We first consider (\ref{dis:11}). By Proposition \ref{lem:uni5},
applied to
$h=l^r_n$, and the same scaling argument that was used to obtain (\ref
{rep-rcomp3}), for every $N,n\in{\mathbb N}$,
\begin{eqnarray*}
&&
\int_0^T \biggl(\int_0^{(\overline{Q}{}^{(N)}(s)+\overline\iota
^{(N)}(s))\wedge
F^{\overline\eta_s^{(N)}}(L)}l^r_n((F^{\overline\eta
_s^{(N)}})^{-1}(y))\,dy \biggr)\,ds \\
&&\qquad= \int
_0^T \biggl(\int_{[0,\chi^{(N)}(s-)\wedge L]}l^r_n(x)\overline{\eta}^{(N)}
_s(dx) \biggr)\,ds \leq\int_0^T \biggl(\int_{[0,L]}l^r_n(x)\overline{\eta}^{(N)}
_s(dx) \biggr)\,ds.
\end{eqnarray*}
By (\ref{def-waitjn}) and the representation of $\eta^{(N)}$ in (\ref
{def-etan}), we have
\begin{eqnarray*}
&&
\int_0^T \biggl(\int_{[0,L]}l^r_n(x)\overline{\eta}^{(N)}_s(dx) \biggr) \,ds \\
&&\qquad\leq\frac{1}{N} \sum_{j = -\mathcal{E}^{(N)}_0+ 1}^{0} \int_0^T
l^r_n\bigl(w^{(N)}_j(0)+s\bigr)
\mathbh{1}_{\{w^{(N)}_j (0)+s <L\wedge r_j\}} \,ds \\
&&\qquad\quad{} +
\frac{1}{N}\sum_{j = 1}^{E^{(N)}(T)} \int_{\zeta_j^{(N)}}^T
l^r_n\bigl(s-\zeta_j^{(N)}\bigr) \mathbh{1}_{\{s-\zeta_j^{(N)} <L\}} \,ds \\
&&\qquad\leq\sup_N \bigl( \bigl\langle1,\overline{\eta}^{(N)}_0 \bigr\rangle+\overline{E}{}^{(N)}(T)
\bigr)\int_0^Ll^r_n(x) \,dx.
\end{eqnarray*}
Since $\sup_N ( \langle 1,\overline{\eta}{}^{(N)}_0 \rangle +\overline{E}{}^{(N)}(t)
)<\infty$ almost surely, due to Assumption \ref{ass-init}, and
$h_n^r$ converges in $\mathcal{L}^1_{\mathrm{loc}}[0,H^r)$ to $h^r$, we obtain
(\ref{dis:11}).
On the other hand, observe that, by (\ref{lem-est-1}) of Lemma \ref
{lem-est} applied to $l=l^r_n$,
\begin{eqnarray*} \int_0^T \biggl(\int_0^{\overline{Q}(s)\wedge
F^{\overline\eta_s}
(L)}l^r_n((F^{\overline\eta_s})^{-1}(y))\,dy \biggr)\,ds &\leq& \int_0^T
\biggl(\int
_{[0,L]}l^r_n(x)\overline{\eta}_s(dx) \biggr)\,ds \\ &\leq& \tilde
L(L,T)\int
_0^Ll^r_n(x)\,dx.
\end{eqnarray*}
By the convergence of $h_n^r$ to $h^r$ in $\mathcal
{L}^1_{\mathrm{loc}}[0,H^r)$, the last
term on the right-hand side of the above display converges to $0$, as
$n\rightarrow\infty$, and (\ref{dis:12}) follows.
\end{pf}
\begin{pf*}{Proof of Proposition \ref{prop:3}}
Given the
discussion prior
to Lemma \ref{lem:uni4} and, in particular, (\ref{dcomp1}), to
complete the
proof of the proposition, it only remains to show that
%
%
\begin{equation}
\label{dis:6} \lim_{L\rightarrow H^r}\limsup_{N\rightarrow\infty
}\mathbb{E}\Bigl[\sup
_{t\in
[0,T]}\overline C{}^{(N)}_i(t,L) \Bigr]=0,\qquad i=1,2,
\end{equation}
and
%
%
\begin{equation}\label{dis:6'}\lim_{L\rightarrow H^r}\mathbb{E}
[\overline
C_3(T,L) ]=
0.
\end{equation}
For the case $i=1$ in (\ref{dis:6}), this follows from Lemma
\ref{lem:2}
and the dominated convergence theorem, whose application is justified
because, by (\ref{eq-cvf1}), (\ref{eq-cvf2}) and the fact that $\overline
{L} \leq
L$,
\begin{eqnarray*}\mathbb{E}\Bigl[\sup_{t\in[0,T]}\overline
C_1^{(N)}(t,L) \Bigr]
&\leq&\mathbb{E}\biggl[\int_0^T \biggl(\int_{[0,L]}h^r(x)\overline{\eta}^{(N)}
_s(dx) \biggr)\,ds \biggr]\\
& &{} +\mathbb{E}\biggl[\int_0^T \biggl(\int
_{[0,L]}h^r(x)\overline{\eta}_s(dx) \biggr)\,ds \biggr],
\end{eqnarray*}
which is bounded uniformly in $N$ by (\ref{lem-est-2}) and Assumption
\ref{ass-init}.

Now, by Remark \ref{rem-compen}, an application of Lemma 5.8(1) of
\cite{kasram07} (with $\nu$, $h^s$ and $H^s$, resp., replaced by
$\eta$, $h^r$ and $H^r$, resp.), shows that
%
%
\begin{equation}
\label{use-temp}
\lim_{L \rightarrow H^r} \sup_{N} \mathbb{E}\biggl[ \int_0^t \biggl( \int
_{[L,H^r)} h^r(x)
\overline{\eta}^{(N)}_s (dx) \biggr) \,ds \biggr] = 0.
\end{equation}
On the other hand, the definition of $\overline C{}^{(N)}_2(T,L)$ in
(\ref{def-c2}), when combined
with Proposition \ref{lem:uni5} and (\ref{eq-cvf2}), shows that
\[
\sup_N\mathbb{E}\bigl[\overline C{}^{(N)}_2(T,L) \bigr]\leq\sup_N
\mathbb{E}\biggl[\int_0^T \biggl(\int_{[L,H^r)} h^r(x) \overline{\eta}^{(N)}
_s(dx) \biggr) \,ds \biggr].
\]
Taking the limit, as $L\rightarrow H^r$, and invoking (\ref
{use-temp}), it follows that (\ref{dis:6}) holds for $i=2$. Finally,
to show (\ref{dis:6'}), we see that, by the definition of $\overline
C_3(T,L)$ in (\ref{def-c3}) and the change of variables formula,
\begin{eqnarray*} \mathbb{E}[\overline C_3(T,L) ] &= &\mathbb{E}
\biggl[\int_0^t \biggl( \int_{\overline{Q}(s)\wedge F^{\overline\eta
_s}(L)}^{\overline{Q}(s)}h^r((F^{\overline\eta_s}
)^{-1}(y))\,dy \biggr)\,ds \biggr]\\
&\leq& \int_0^t \biggl( \int
_{[L,H^r)}h^r(x)\overline{\eta}_s(dx) \biggr)\,ds.
\end{eqnarray*}
If $h^r$ is bounded, then (\ref{dis:6'}) holds by simply applying the
bounded convergence theorem on the right-hand side of the equality in
the above display. On the other hand, suppose $h^r$ is
lower-semicontinuous on $(L^r,H^r)$ for some $L^r<H^r$. Then, by
Theorem A.3.12 of
\cite{dupellbook} and the
fact that $\mathbb{P}$ a.s., $\overline{\eta}^{(N)}_s \stackrel
{w}{\rightarrow}\overline{\eta}_s$, as $N\rightarrow\infty$,
for a.e. $s\in[0,T]$, this implies that for any
such $s$ and $L>L^r$,
\[
\int_0^t \biggl( \int_{[L,H^r)}h^r(x)\overline{\eta}_s(dx) \biggr)\,ds \leq
\liminf_{N\rightarrow\infty} \int_0^t \biggl( \int
_{[L,H^r)}h^r(x)\overline{\eta}^{(N)}
_s(dx) \biggr)\,ds.
\]
Integrating both sides over $s\in[0,T]$ and taking expectations, an
application of
Fatou's lemma yields
\[
\mathbb{E}[\overline C_3(T,L) ]\leq\liminf_{N\rightarrow\infty
}\mathbb{E}
\biggl[\int_0^t \biggl( \int_{[L,H^r)}h^r(x)\overline{\eta}^{(N)}_s(dx)
\biggr)\,ds \biggr].
\]
Taking the limit as $L\rightarrow H^r$, an application of (\ref{use-temp})
shows that
(\ref{dis:6'}) holds.
\end{pf*}

We now prove the main limit result.
\begin{pf*}{Proof of Theorem \ref{thm:FE}}
Fix $t\in[0,\infty)$ such that $\overline{\nu}{}^{(N)}_t \stackrel
{w}{\rightarrow}\overline{\nu}_t$, $\overline{\eta}^{(N)}_t
\stackrel{w}{\rightarrow}\overline{\eta}_t$, $\overline
{E}{}^{(N)}(t)\rightarrow\overline{E}(t)$, $\overline
{X}{}^{(N)}(t)\rightarrow
\overline{X}(t)$, $\overline{R}{}^{(N)}(t)\rightarrow\overline
{R}(t)$, $\overline{A}{}^{(N)}_{\cdot,\nu
}(t)\stackrel{w}{\rightarrow}\overline{A}_{\cdot,\nu}(t)$,
$\overline{D}{}^{(N)}_{\cdot}(t) \stackrel{w}{\rightarrow}
\overline{A}_{\cdot,\nu}(t)$, $\overline{A}{}^{(N)}_{\cdot,\eta
}(t)\stackrel{w}{\rightarrow}
\overline{A}_{\cdot,\eta}(t)$, $\overline{S}{}^{(N)}_{\cdot
}(t)\stackrel{w}{\rightarrow}\overline{A}_{\cdot
,\eta}(t)$ as $N\rightarrow\infty$. Since $\overline Y{}^{(N)}\rightarrow\overline Y$ a.s., this occurs for $t$ outside a
countable set. By (\ref{dis:13}), this implies that as $N\rightarrow
\infty$,
%
%
\begin{eqnarray}\overline{D}{}^{(N)}_{\varphi}(t)\rightarrow\overline
{A}_{\varphi,\nu
}(t)=\int
_0^t\langle\varphi(\cdot,s)h^s(\cdot,s),\overline{\nu}_s\rangle
\,ds,\nonumber\\[-8pt]\\[-8pt]
\eqntext{\varphi\in\mathcal{C}_b\bigl([0,H^s)\times{\mathbb
R}_+\bigr).}
\end{eqnarray}
An analogous argument also implies that, as $N\rightarrow\infty$,
%
%
\begin{eqnarray}\label{dis:fsn}\overline{S}{}^{(N)}_{\psi
}(t)\rightarrow\overline{A}
_{\psi,\eta
}(t)=\int_0^t\langle\psi(\cdot,s)h^r(\cdot,s),\overline{\eta
}_s\rangle \,ds,\nonumber\\[-8pt]\\[-8pt]
\eqntext{\psi\in\mathcal{C}_b\bigl([0,H^r)\times{\mathbb R}_+\bigr).}
\end{eqnarray}
In particular, when $\varphi=\psi={\mathbf{1}}$, the above two displays
imply that (\ref{cond-radon}) holds.
Also, we immediately obtain that, as $N\rightarrow\infty$, $\langle
\mathbf{1}, \overline{\nu}{}^{(N)}_t\rangle\rightarrow\langle
{\mathbf{1}}, \overline{\nu}_t\rangle$ and
$\langle
{\mathbf{1}}, \overline{\eta}^{(N)}_t\rangle\rightarrow\langle
{\mathbf{1}}, \overline{\eta}
_t\rangle$.
When combining with (\ref{def-xn}), (\ref{comp-prelimit}), (\ref
{def-kn}), (\ref{equivD}), (\ref{def-dn}), (\ref{qn}), (\ref
{rep-R}), this implies that all the equations in Definition \ref
{def-fleqns} are satisfied at time $t$ except (\ref{eq-ftmeas}) and
(\ref{eq-ftreneg}).

It only remains to show that (\ref{eq-ftmeas}) and (\ref{eq-ftreneg})
are also satisfied at time $t$. We shall just prove (\ref
{eq-ftreneg}). The same argument will also show that (\ref{eq-ftmeas})
holds. Dividing (\ref{eqn-prelimit3}) by $N$, we have
\begin{eqnarray*}
\bigl\langle\psi(\cdot, t), \overline{\eta}^{(N)}_{t} \bigr\rangle
& = & \bigl\langle\psi(\cdot, 0), \overline{\eta}^{(N)}_{0} \bigr\rangle+
\int_{0}^t \bigl\langle\psi_x(\cdot,s) + \psi_s(\cdot,s), \overline
{\eta}^{(N)}
_s \bigr\rangle \,ds \\
& &{} - \overline S{}^{(N)}_\psi(t) + \int_{[0,t]} \psi(0,s) \,d\overline
{E}{}^{(N)}(s).\nonumber
\end{eqnarray*}
Since $\overline{\eta}^{(N)}_{0} \stackrel{w}{\rightarrow}\overline
{\eta}_0$ by Assumption \ref
{ass-init}(4), $\overline{\eta}^{(N)}_s \stackrel{w}{\rightarrow
}\overline{\eta}_s$ for a.e. $s \in[0, t]$,
$\overline{\eta}^{(N)}_t \stackrel{w}{\rightarrow}\overline{\eta
}_t$ by our
choice of $t$ and $\psi(\cdot, t)$ and $\psi_x(\cdot, s)+\psi
_s(\cdot, s), s\in[0, t]$, are bounded and continuous,
as $N\rightarrow\infty$, we have
\[
\bigl\langle\psi(\cdot, t), \overline{\eta}^{(N)}_{t} \bigr\rangle
\rightarrow\langle
\psi(\cdot, t), \overline{\eta}_{t} \rangle\quad\mbox{and}\quad
\bigl\langle\psi(\cdot, 0), \overline{\eta}^{(N)}_{0} \bigr\rangle
\rightarrow\langle
\psi(\cdot, 0), \overline{\eta}_{0} \rangle,
\]
and, by the bounded convergence theorem,
\[
\int_{0}^t \bigl\langle\psi_x(\cdot,s) + \psi_s(\cdot,s), \overline
{\eta}^{(N)}
_s \bigr\rangle \,ds \rightarrow\int_{0}^t \langle\psi_x(\cdot,s) + \psi
_s(\cdot,s), \overline{\eta}_s \rangle \,ds.
\]
On the other hand, using an integration-by-parts argument, the facts
that $\overline{E}{}^{(N)}(0)=0$, $\overline{E}{}^{(N)}\rightarrow
\overline{E}$,
$\overline{E}$ is nondecreasing and $\psi_s(0, \cdot)$ is bounded and
continuous on $[0,t]$,
along with the bounded convergence theorem, we see that, as
$N\rightarrow
\infty$,
\[
\int_{[0,t]} \psi(0,s) \,d\overline{E}{}^{(N)}(s) \rightarrow\int
_{[0,t]} \psi(0,s) \,d\overline{E}(s).
\]
Combining the last four displays with (\ref{dis:fsn}), it follows that
(\ref{eq-ftreneg}) holds.
Then it follows that all fluid equations are satisfied for all but
countably many $t$. By right-continuity (with respect to $t$) of each
of the terms in all fluid equations, we conclude that all fluid
equations are a.s. satisfied for all $t\in[0,\infty)$. This completes
the proof of the desired result that $(\overline{X},\overline{\nu
},\overline{\eta})$
satisfies the fluid equations.
\end{pf*}

\subsection{\texorpdfstring{Proof of Theorem \protect\ref{thm:3}}{Proof of Theorem 3.8}}
\label{subs-prf3}

This section is devoted to the proof of Theorem~\ref{thm:3}.
Recall ${\mathcal T}_t^{(N)}(s)$ in (\ref{dis:Tn}) and its fluid
scaled version defined in (\ref{fl-scaling2}).
Observe that the virtual waiting time defined in (\ref{T}) can be
rewritten in terms of the fluid-scaled quantities as
%
%
\begin{equation} \label{T2}\qquad
W^{(N)}(t) \doteq\inf\bigl\{s\geq0\dvtx\overline{D}{}^{(N)}(t+s)
-\overline{D}{}^{(N)}(t)+\overline
{\mathcal T}{}^{(N)}_t(s) > \overline{Q}{}^{(N)}(t) \bigr\}.
\end{equation}
We first show that for each $t\in[0,\infty)$, $\overline
{\mathcal T}{}^{(N)}_t
\Rightarrow\overline{\mathcal T}_t$ as $N \rightarrow\infty$, where
$\overline{\mathcal T}_t$ is defined in (\ref{dis:T}). Notice that
a customer $j$ who arrived into the system before time $t$ and has
not reneged by time $t$ must have a potential waiting time
$w^{(N)}_j(u) > u-t$ for all $u > t$ sufficiently small.
In addition, for that customer to have reneged from the queue
(before entering service) in the period $[t,t+s]$, there must
exist a time $u \in[t,t+s]$ such that the customer is still in queue
(i.e., has not yet entered service) or, equivalently, such that
$w^{(N)}_j(u) < \chi^{(N)}(u-)$,
the waiting time of the head-of-the-line customer just prior to $u$,
and the customer reneges, so that the slope of her
potential waiting time changes
from one to zero.
Therefore, for each $s\in
[0,\infty)$, $\mathcal T_t^{(N)}(s)$ can be alternatively expressed as
\begin{eqnarray*}
\mathcal T_t^{(N)}(s) &=& \sum_{u \in[t,t+s]} \sum_{j=-\mathcal{E}^{(N)}_0+
1}^{E^{(N)}
(u)}
\mathbh{1}_{ \{{dw^{(N)}_j }/{dt}(u-) >0, {dw^{(N)}_j }/{dt}(u
+)=0 \}}\\
&&\hspace*{78.85pt}{}\times
\mathbh{1}_{\{u-t<w^{(N)}_j (u)\leq\chi^{(N)}(u-)\}}.
\end{eqnarray*}
Let
\[
\mathcal T_t^{(N),1}(s) \doteq\sum_{u \in[t,t+s]} \sum_{j=-\mathcal
{E}^{(N)}_0+
1}^{E^{(N)}
(u)}
\mathbh{1}_{ \{{dw^{(N)}_j }/{dt}(u-) >0, {dw^{(N)}_j }/{dt}(u
+)=0 \}}
\mathbh{1}_{\{w^{(N)}_j (u)\leq\chi^{(N)}(u-)\}}
\]
and
\[
\mathcal T_t^{(N),2}(s) \doteq\sum_{u \in[t,t+s]} \sum_{j=-\mathcal
{E}^{(N)}_0+
1}^{E^{(N)}
(u)}
\mathbh{1}_{ \{{dw^{(N)}_j }/{dt}(u-) >0, {dw^{(N)}_j }/{dt}(u
+)=0 \}}
\mathbh{1}_{\{w^{(N)}_j (u)\leq u-t\}}.
\]
It is easy to see that $\mathcal T_t^{(N)}(s)=\mathcal T_t^{(N),1}
(s)-\mathcal T_t^{(N),2}(s)$, $\mathcal T_t^{(N),1}(s)=R^{(N)}(t+s)-R^{(N)}
(t)$, $\mathcal T_t^{(N),2}(s) \leq S^{(N)}(t+s)-S^{(N)}(t)$
and $\mathcal T_t^{(N),2}(s+\delta)-\mathcal T_t^{(N),2}(s) \leq S^{(N)}
(t+s+\delta)-S^{(N)}(t+s)$.
Therefore, an application of Kurtz's criteria in Proposition \ref{Kurtz}
shows
that the relative compactness of
the fluid scaled versions
$\overline{\mathcal T}{}^{(N),1}_t$ and $\overline{\mathcal T}{}^{(N),2}_t$
of
${\mathcal T}_t^{(N),1}$ and ${\mathcal T}_t^{(N),2}$, respectively,
follows from that of $\overline{R}{}^{(N)}$ and $\overline{S}{}^{(N)}$
established in Lemma \ref{lem:rc}.
By a straightforward adaption of the argument used in Proposition \ref
{prop:3} to
show the convergence of $\overline{R}{}^{(N)}$ to $\overline{R}$, we can
conclude that
$\overline{\mathcal T}{}^{(N)}_t\Rightarrow\overline{\mathcal T}_t$
as $N\rightarrow\infty$.

Recall the application of the Skorokhod representation theorem in Theorem
\ref{thm:FE} to assume, without loss of generality, that $\overline Y{}^{(N)}$
converges a.s. to $\overline Y$. Here, we can also assume, in
addition, that
$\overline{\mathcal T}{}^{(N)}_t(s)\rightarrow\overline{\mathcal T}_t$ a.s.,
as $N\rightarrow\infty$. Since $\overline{Q}$ is continuous at $t$
and, by
(\ref{dis:13}), $\overline{A}_{{\mathbf{1}},\nu}=\int_0^\cdot
\langle
h^s,\overline{\nu}_s\rangle \,ds$ is
continuous by the integral representation,
and $\overline{\mathcal T}_t$ has continuous paths by definition,
it follows that, almost surely, $\overline{Q}{}^{(N)}(t)\rightarrow
\overline{Q}(t)$ and for
each $T\in[0,\infty)$, as $N\rightarrow\infty$,
\[
\sup_{s\in[0,T]}\bigl|\overline{D}{}^{(N)}(t+s)-\overline{A}_{{\mathbf
{1}},\nu
}(t+s)\bigr|\rightarrow0 \quad\mbox{and}\quad \sup_{s\in
[0,T]}\bigl|\overline{\mathcal T}{}^{(N)}_t(s)- \overline{\mathcal
T}_t\bigr|\rightarrow
0.
\]
From (\ref{T2}), it is easy to see that
$W^{(N)}(t)\leq(\overline{D}{}^{(N)})^{-1}(\overline
{D}{}^{(N)}(t)+\overline{Q}{}^{(N)}(t))-t$ for each $N$. By the
tightness result established in
Theorem \ref{th-tight}, we know that
$\overline{D}{}^{(N)}(t)+\overline{Q}{}^{(N)}(t)$ is bounded uniformly in
$N$, and due to Lemma 4.10 of
\cite{ramrei03} and the assumption that $\overline{A}_{{\mathbf
{1}},\nu
}$ is
uniformly strictly increasing, we also know that
$(\overline{D}{}^{(N)})^{-1}\rightarrow
(\overline{A}_{{\mathbf{1}},\nu} )^{-1}$ uniformly on compact
sets, as $N\rightarrow\infty$.
Hence,
$W^{(N)}(t)$ is bounded uniformly in $N$. Therefore, there exists a
subsequence, $W^{(N_n)}(t)$, $n \in{\mathbb N}$, that converges to a
limit in
$[0,\infty)$,
which we denote by $W^*$. From (\ref{T2}) and
the right-continuity of $\overline{D}{}^{(N)}, \overline{Q}{}^{(N)}$ and
$\overline{\mathcal T}{}^{(N)}_t$, we then have
$\overline D{}^{(N_n)}(t+\overline W^{(N_n)}(t))-\overline
D{}^{(N_n)}(t)+\overline{\mathcal
T}{}^{(N_n)}_t(\overline W{}^{(N_n)}(t))\geq\overline Q{}^{(N_n)}(t)$.
Sending $n \rightarrow\infty$, we
obtain
%
%
\begin{equation}
\label{ineq-Wstar}
\overline{A}_{{\mathbf{1}},\nu}(t+W^*)- \overline{A}_{{\mathbf
{1}},\nu
}(t)+\overline{\mathcal T}_t(W^*)\geq\overline{Q}(t).
\end{equation}
Together with (\ref{Wbar}),
this shows that $\overline W(t)\leq W^*$. Now, suppose that $\overline
W(t)< W^*$,
and
fix $w$ such that $\overline W(t)<w< W^*$.
Since $\overline{A}_{{\mathbf{1}},\nu}$ is uniformly strictly
increasing and $\overline{\mathcal T}_t$ is nondecreasing, the
inequality $\overline W(t)<w$ implies
that $\overline{A}_{{\mathbf{1}},\nu}(t+w)-\overline{A}_{{\mathbf
{1}},\nu
}(t)+\overline{\mathcal T}_t(w) > \overline{Q}(t)$.
Therefore, for sufficiently large $N$,
we have $\overline{D}{}^{(N)}(t+w)-\overline{D}{}^{(N)}(t)+\overline
{\mathcal T}{}^{(N)}_t(w) > \overline{Q}{}^{(N)}(t)$ and hence $W^{(N)}(t)
\leq w$. In
turn, this implies
that $W^{(N_n)}(t) \leq w$ for sufficiently large $n \in{\mathbb N}$.
Sending $n\rightarrow\infty$ and using the convergence of
$W^{(N_n)}(t)$ to $W^*$, we then obtain $W^*\leq w$.
This contradicts the choice of $w$. Hence $\overline W(t)= W^*$, and
this proves the
desired result.

\begin{appendix}

\section{Explicit construction of the state processes}
\label{ap-markov}

In this section, we construct all state processes and auxiliary
processes described in Section \ref{sec:repdyn} from the initial data
$\{\mathcal{E}^{(N)}_0, X^{(N)}(0),w^{(N)}_j(0), a^{(N)}_j(0),
j=-\mathcal{E}^{(N)}_0+1,\ldots,0\}$, $\{
\alpha_E^{(N)}(t),t\in[0,\infty)\}$, $\{v_j, j\in{\mathbb Z}\}$ and
$\{r_j,
j\in{\mathbb Z}\}$.

Fix $N$ and, for simplicity, we omit the dependence on $N$ in notation.
Let \mbox{$E(0)=0$}. The process $E$ on $[0,\infty)$ can be obtained from
$\alpha_E$ using the relation (\ref{def-ren}). Let $\ell=0$, $\tau
_0=0$, and let $R(\tau_\ell)=D(\tau_\ell)=K(\tau_\ell)=0$,
%
%
\begin{equation}
\label{app-q}
Q(\tau_\ell) \doteq[X(\tau_\ell) - N]^+,
\end{equation}
and for $j>E(\tau_\ell)$, let $w_j(\tau_\ell)=a_j(\tau_\ell)=0$.
Now, for $t \in[\tau_\ell, \infty)$, define
%
%
\begin{equation}
\label{app-chi}
\chi^\ell(t) \doteq
\inf\{ x > 0\dvtx\eta_{\tau_\ell} [0,x] \geq Q(\tau_\ell)\} + t -
\tau_{\ell}.
\end{equation}
Also, for $j=-\mathcal{E}_0+1,\ldots,0,\ldots,E(\tau_\ell)$ and
$t\in[\tau_\ell,\infty)$, let
\begin{eqnarray*}
w_j^\ell(t) &\doteq& \bigl(w_j(\tau_\ell)+t-\tau_\ell\bigr)\wedge r
_j, \\
a_j^\ell(t) &\doteq& \cases{
0, &\quad if $w_j(\tau_\ell)=r_j$
or $w_j(\tau_\ell)\leq\chi^\ell(\tau_\ell)$, \vspace*{2pt}\cr
\bigl(a_j(\tau_\ell)+t-\tau_\ell\bigr)\wedge v_j, &\quad if $\chi^\ell(\tau
_\ell
)<w_j(\tau_\ell)<r_j$,}
\\
\eta^\ell_t &\doteq& \sum_{j=-\mathcal{E}_0+1}^{E(\tau
_\ell)}\delta_{w_j(t)}\mathbh{1}_{ \{{d w_j}/{dt}(t+) >0 \}
}, \\
\nu^\ell_t &\doteq& \sum_{j=-\mathcal{E}_0+1}^{E(\tau_\ell
)}\delta_{a_j(t)}\mathbh{1}_{ \{{d a_j}/{dt}(t+) >0 \}}, \\
R^\ell(t) &\doteq& \sum_{j=-\mathcal{E}_0+1}^{E(\tau_\ell)}\sum
_{s\in[0,t]}\mathbh{1}_{ \{w_j(s)\leq\chi^l(s-), {d
w_j}/{dt}(s-) >0, {d w_j }/{dt}(s+)=0 \}}, \\
D^\ell(t)
&\doteq& \sum_{j=-\mathcal{E}_0+1}^{E(\tau_\ell)}\sum_{s\in
[0,t]}\mathbh{1}_{ \{{d a_j}/{dt}(s-) >0, {d a_j
}/{dt}(s+)=0 \}}.
\end{eqnarray*}
Next, define
\[
\tau_{\ell+ 1} \doteq\inf\bigl\{t > 0\dvtx\bigl(D^{\ell}(t) -
D(\tau_{\ell})\bigr) \wedge
\bigl(R^{\ell}(t) - R(\tau_{\ell})\bigr) \wedge\bigl(E(t) - E(\tau_{\ell})\bigr) > 0
\bigr\}.
\]
For $t \in[\tau_{\ell},\tau_{\ell+1})$, let $Y(t) = Y^\ell(t)$ for $Y =
w_j, a_j, j \in-\mathcal{E}_0+1, \ldots, E(\tau_\ell)$, $R, D, \eta,
\nu$ and $\chi$ and set $Y(t) = Y(\tau_\ell)$ for $Y = X, Q, w_j, a_j,
j>E(\tau_\ell)$. Moreover, define
\begin{eqnarray*}
X(\tau_{\ell+1}) & \doteq& X(\tau_{\ell}) + E(\tau_{\ell+1}) -
E(\tau_{\ell})
- D(\tau_{\ell+1}) + D(\tau_{\ell}) \\
& &{} - R(\tau_{\ell+1}) + R(\tau_{\ell}), \\
\eta_{\tau_{\ell+1}} & \doteq& \eta^{\ell}_{\tau_{\ell+1}} +
\bigl(E(\tau_{\ell+1}) - E(\tau_{\ell})\bigr) \delta_0,
\end{eqnarray*}
and, if $E(\tau_{\ell+ 1}) > E(\tau_{\ell})$, then $E(\tau_{\ell+
1}) = E(\tau_{\ell})+1$, and then let
$w_j (\tau_{\ell+1}) \doteq0$ for $j \in\{ E(\tau_{\ell})+1,
\ldots,
E(\tau_{\ell+1})\}$.
In this case, $Q(\tau_{\ell+1})$ and $\chi(\tau_{\ell+1})$ can be
defined via equations (\ref{app-q}) and (\ref{app-chi}), but with
$\ell$
replaced by $\ell+ 1$, and
the procedure can be reiterated.
Now, $\max\{\ell\dvtx\tau_\ell\leq t\}$ is
bounded by $\mathcal{E}_0 + E(t)$, and is therefore
a.s. finite. Therefore, $\tau_{\ell} \rightarrow\infty$ as $\ell
\rightarrow
\infty$,
and so the above procedure constructs the above processes on
$[0,\infty)$. $K$ and $S$ can then be defined, respectively,
via equations (\ref{def-kn}) and (\ref{def-sn}).

For each $j\geq-\mathcal{E}^{(N)}_0$, by the construction, we have
\begin{eqnarray*}w_j(t)&=&\sum_{E(\ell)\geq j}\mathbh{1}_{[\tau
_\ell
,\tau_\ell+1)}(t)\bigl(w_j(\tau_\ell)+t-\tau_\ell\bigr)\wedge r_j \\
&=& \cases{
t\wedge r_j, &\quad if $j=-\mathcal{E}^{(N)}_0,\ldots,
0$,\cr
(t-\zeta_j)\wedge r_j, &\quad otherwise,}
\end{eqnarray*}
where $\zeta_j=\inf\{t>0\dvtx E(t)=j\}$. Hence the process $w_j$ defined
above is
indeed the potential waiting time process of customer $j$. It is also
not to
hard to see that the process $a_j$ defined above is the age process of
customer $j$ and satisfies (\ref{adif}). We next show that the process
$\chi$ constructed above satisfies (\ref{def-chi}). It is easy to see
that $\chi(0)=\chi^0(0)$ by (\ref{app-chi}) with $t=0$ and $\ell=
0$. The $\chi(0)$ satisfies (\ref{def-chi}) for $t=0$. When $t\in
[\tau_0, \tau_1)$, $Q(t)=Q(0)$, $\eta_t=\eta^0_t$ and $\chi
(t)=\chi^0(t)$. Then we have
\[
\chi^0(t)=\inf\{ x > 0\dvtx\eta_{\tau_0} [0,x] \geq Q(\tau_0)\} + t
- \tau_{0}=\inf\{ x > 0\dvtx\eta_{t} [0,x] \geq Q(t)\}.
\]
Hence $\chi$ satisfies (\ref{def-chi}) on the interval $[\tau_0,\tau
_1)$. By the standard induction argument, we can see that $\chi$
satisfies (\ref{def-chi}) for all $t\geq0$.

For each $t\geq0$, by the construction, we have
\begin{eqnarray*}\eta_t &=&\sum_{\ell= 0}^\infty\mathbh{1}_{[\tau
_\ell
,\tau_\ell+1)}(t)\sum_{j=-\mathcal{E}_0+1}^{E(\tau_\ell)}\delta
_{w_j(t)}\mathbh{1}_{ \{{d w_j}/{dt}(t+) >0 \}} \\
&=& \sum
_{\ell= 0}^\infty\mathbh{1}_{[\tau_\ell,\tau_\ell+1)}(t)\sum
_{j=-\mathcal{E}_0+1}^{E(t)}\delta_{w_j(t)}\mathbh{1}_{ \{{d
w_j}/{dt}(t+) >0 \}} \\
&=& \sum_{j=-\mathcal
{E}_0+1}^{E(t)}\delta_{w_j(t)}\mathbh{1}_{ \{{d w_j}/{dt}(t+)
>0 \}} .
\end{eqnarray*}
This shows that the $\eta$ constructed satisfies (\ref{def-etan}). A
similar argument shows that the processes $\nu$, $D$ and $R$
constructed satisfy (\ref{def-nun}), (\ref{def-depart}) and (\ref
{def-crp}), respectively. Finally, $K$ and $S$ satisfy (\ref{def-kn})
and (\ref{def-sn}) by construction.

Recall that, for $t\in[0,\infty)$, $\tilde{\mathcal{F}}_t$ is the
$\sigma$-algebra generated by
\[
\bigl(\mathcal{E}_0,X(0),\alpha_E(s), w_j(s), a_j(s), j \in\{-\mathcal
{E}_0+1, \ldots, 0\}\cup{\mathbb N}, s \in[0,t]\bigr\}
\]
and $\{\mathcal{F}_t\}$ is the associated completed,
right-continuous filtration.
\begin{lemma} \label{app:adapt}
The processes $w_j, a_j, j \geq-\mathcal{E}_0+1$ and $E, R, D,
\eta
, \nu, \chi, X, Q$, $K, S$ are c\`{a}dl\`{a}g and $\{\mathcal
{F}_t\}
$-adapted.
\end{lemma}
\begin{pf}
The c\`{a}dl\`{a}g property of those processes follows from the
construction. Now we show that all the processes are $\{\mathcal{F}_t\}
$-adapted. Indeed, it follows immediately from (\ref{def-ren}),
(\ref{def-etan}), (\ref{def-nun}), (\ref{def-depart}) and (\ref{def-cvrp})
that $E, \eta, \nu, D$ and $S$ are $\mathcal{F}_t$-adapted. We next
show that
$\chi$ is $\mathcal{F}_t$-adapted. By equations (\ref{def-qnt}) and
(\ref{def-chi}) evaluated at time~$0$, it follows that $\chi(0)$ is a function
of $X(0)$ and $\eta_0$ and hence $\mathcal{F}_0$-adapted. Now, let $t>
0$. For
each $\ell\geq0$, by the induction argument, $\chi^\ell(t)$ is
$\mathcal{F}_t$-adapted, and $\tau_\ell$ is an $\mathcal
{F}_t$-stopping time. Since
$\chi_t=\chi^\ell_t$ if $t\in[\tau_\ell,\tau_{\ell+1})$, $\chi
$ is
$\mathcal{F}_t$-adapted.~Equations (\ref{def-crp}) and (\ref{def-dn})
show that $R$
and $X$ are $\mathcal{F}_t$-adapted,
and it follows from (\ref{def-qnt}) and (\ref{def-kn}) that $Q$ and
$K$ are
$\mathcal{F}_t$-adapted.
\end{pf}

The next lemma establishes some basic properties of $\chi(t)$, the
waiting time of the head-of-the-line customer at time $t$, defined in
(\ref{def-chi}).
\begin{lemma} \label{lem-chi} $\chi$ is piecewise linear with
downward jumps
that occur when the head-of-the-line customer either enters service
(due to
a departure from service) or reneges from the queue. Hence, $\chi
(t-)\geq\chi(t)$ for every $t\in(0,\infty)$. Moreover, for every
$t>0$, there exists $\varepsilon_t(\omega)\in(0,t)$ such that for
all $\tilde t\in(t-\varepsilon_t(\omega),t)$, $\chi(t-)-\chi
(\tilde t-)=t-\tilde t>0$.
\end{lemma}
\begin{pf}
By the construction, $\chi_t=\chi^\ell_t$ if $t\in
[\tau_\ell,\tau_{\ell+1})$. Since $\chi^\ell$ is linear on
$[\tau_\ell,\tau_{\ell+1})$, $\chi$ is piecewise linear. Also
$\chi$
can only jump at $\tau_{\ell+1}$, $\ell\geq0$. Based on the definition
of $\tau_{\ell+1}$, it is not hard to see that $\chi$ can only have a
downward jump at $\tau_{\ell+1}$ when the head-of-the-line customer
either enters service [$D^{\ell}(\tau_{\ell+1}) - D(\tau_{\ell
})>0$] or
reneges from the queue [$R^{\ell}(\tau_{\ell+1}) - R(\tau_{\ell})>0$].
Then we have $\chi(t-)\geq\chi(t)$ for every $t\in(0,\infty)$. The
last statement of the lemma follows from the fact that $\chi$ is
c\`{a}dl\`{a}g and piecewise linear.
\end{pf}\vspace*{-14pt}


\section{Strong Markov property} \label{sec:SMF}

In this section we show that the state descriptor $V^{(N)} =
(\alpha_E^{(N)}, X^{(N)}, \nu^{(N)},\break \eta^{(N)})$ is a strong Markov
process with
respect to the filtration $\{\mathcal{F}_t^{(N)}, t\geq0\}$ defined in
Section \ref{subsub-filt}. To ease the notation, we shall suppress the
superscript $(N)$ from the notation.

Let $\mathcal{M}_D[0,H^s)$ and $\mathcal{M}_D[0,H^r)$ be the subsets
of $\mathcal{M}_F[0,H^s)$ and $\mathcal{M}_F[0,H^r)$,
respectively, such that each measure in $\mathcal{M}_D[0,H^s)$ and
$\mathcal{M}_D[0,H^r)$ takes the
form $\sum_{i=1}^k \delta_{x_i}$. Define
%
%
\begin{equation}
\mathcal{V}\doteq\left\{
\matrix{(\alpha, x, \mu, \pi) \in{\mathbb R}_+\times{\mathbb Z}
_+\times\mathcal{M}_D[0,H^s)\times\mathcal{M}_D[0,H^r)\mbox{:}\cr
x\leq\langle{\mathbf{1}}, \mu\rangle+\langle{\mathbf{1}}, \pi
\rangle,
\langle{\mathbf{1}},
\mu\rangle\leq N}
\right\},
\end{equation}
where ${\mathbb R}_+$ is endowed with the Euclidean topology $d$,
${\mathbb Z}_+$ is
endowed with the discrete topology $\rho$ and $\mathcal{M}_D[0,H^s)$
and $\mathcal{M}_D[0,H^r)$
are endowed with the weak topology, respectively. The space
$\mathcal{V}$ is a closed subset of ${\mathbb R}_+\times{\mathbb
Z}_+\times\mathcal{M}_F[0,H^s)
\times\mathcal{M}_F[0,H^r)$ and is endowed with the usual product
topology. Since
${\mathbb R}_+\times{\mathbb Z}_+\times\mathcal{M}_F[0,H^s)\times
\mathcal{M}_F[0,H^r)$ is a Polish space, then
the closed subset $\mathcal{V}$ is also a Polish space. Now, denote
\[
V(t)\doteq(\alpha_E(t), X(t), \nu_t,\eta_t),\qquad t\geq0.
\]
It is obvious that $V$ is a $\mathcal{V}$-valued process adapted
to the filtration $\{\mathcal{F}_t^V, t\geq0\}$, the natural
filtration generated by $V$.

For each $y, z\in\mathcal{V}$ and $t\geq0$, let
%
%
\begin{equation}\label{dis:semiG}
P_t(y,z)=\mathbb{P}\bigl(V(t)=z|V(0)=y\bigr).
\end{equation}
For any measurable function
$\psi$ defined on
$\mathcal{V}$ and $t\geq0$, define the function $P_t\psi$ on
$\mathcal{V}$ as
%
%
\begin{equation}
\label{def:Ppsi} P_t\psi(y)=\mathbb{E}[\psi(V(t))|V(0)=y],\qquad
y\in\mathcal{V}.
\end{equation}

\begin{lemma} \label{lem:Mark}
The state descriptor $V$ is strong Markov with respect to $\{
\mathcal{F}_t$, $t\geq0\}$, and hence is strong Markov with respect to
$\{
\mathcal{F}_t^V, t\geq0\}$. Moreover, $\{P_t, t\geq0\}$ in
(\ref
{dis:semiG}) is the Markov semigroup of $V$.
\end{lemma}
\begin{pf}
To establish the strong Markov property, we shall identify $V$ as
a, so-called, piecewise deterministic Markov process (cf. \cite
{jacobsen}). From the explicit pathwise construction of $V$ in
Appendix \ref{ap-markov}, it follows that $V$ is a piecewise deterministic
process with jump times $\{\tau_1,\tau_2,\ldots\}$. Each jump time
is either the arrival time of a new customers or the time of a service
completion or the time to the end of a patience time. Note that, due to
the nonidling condition, the time of entry into service of a customer
must coincide with either the arrival time of that customer or the time
of service completion of another customer. Let $\tau_0=0$. For each
integer $n\geq0$, let $P_n=V(\tau_n)$. Then $\{(\tau_n,P_n),
n\geq0\}$ forms a marked point process. For each $n\geq0$, $V$
evolves in a deterministic fashion on $[\tau_n,\tau_{n+1})$. For each
$t\geq0$ and $y\in\mathcal{V}$ with $y=(\alpha, x, \sum_{i=1}^k
\delta_{u_i}, \sum_{j=1}^l \delta_{z_i})$ and $k\leq N$, define
{\renewcommand{\theequation}{$*$}
\begin{equation}\label{zv1}
\phi_t(y)\doteq \Biggl(\alpha+t, x, \sum_{i=1}^k \delta_{u_i+t}, \sum
_{j=1}^l \delta_{z_i+t} \Biggr).
\end{equation}}

\noindent It is easy to see that
\[
\phi_{t+s}(y)=\phi_s(\phi_t(y)),\qquad \phi_0(y)=y,
\]
and the map $t\mapsto\phi_t(y)$ is continuous in the interval
$[0,\infty)$. For each $t\geq0$, let
\[
\langle t\rangle=\max\{n\geq1\dvtx\tau_n \leq t\}
\]
with the convention that $\max\varnothing= 0$. We can see that
%
%
\setcounter{equation}{3}
\begin{equation}
V(t)=\phi_{t-\tau_{\langle t\rangle}}(V_{\tau_{\langle t\rangle
}}).
\end{equation}

The jump dynamics are captured by $\{r_t(y,C), t\geq0, y\in
\mathcal{V}, C\subset\mathcal{V}\}$. For each $t\geq0, y\in
\mathcal{V}, C\subset\mathcal{V}$, $r_t(y,C)$ is the conditional
probability that a jump leads to a state in\vspace*{1pt} $C$, given that the jump
occurs at time $t$ from state $y$. Let $y=(\alpha, x, \sum_{i=1}^k
\delta_{u_i}, \sum_{j=1}^l \delta_{z_i})$. Recall that there are only
three types of jump times for the process $V$. Given that $V$
jumps at time $t$ from state $y$, if we know which type the jump time
$t$ is, then we know to which state the process $V$ jumps to. For
example, suppose that the number $k$ in the expression of $y$ is less
than $N$, then, at state $y$, there is at least one idle server. If the
jump\vspace*{1pt} is due to the new arrival, then the process $V$ will jump\vspace*{1pt} to
state $(0, x+1, \sum_{i=1}^k \delta_{u_i}+\delta_0, \sum_{j=1}^l
\delta_{z_i}+\delta_0)$. Let $p_1, p_2, p_3$, respectively, be the
conditional probability that the jump at time $t$ is due to the arrival
of a new customer, service completion of a customer in service, the end
of patience time for some customer in the system, respectively, given
that the jump occurs at time $t$ from state $y$. Then the probability
measure $r_t(y,\cdot)$ can be easily written from $y$ and $p_i$,
$i=1,2,3$.

The jump time dynamics are captured by the survivor functions
$\{\overline H_{s,y}(t)\dvtx\break 0\leq s\leq t, y\in\mathcal{V}\}$, where
$\overline H_{s,y}(t)$ is the conditional probability that the time for
the next jump is more than time $t$ given the state being at $y$ at
time $s$, in other words, for $y=(\alpha, x, \sum_{i=1}^k
\delta_{u_i}, \sum_{j=1}^l \delta_{z_i})$,
{\renewcommand{\theequation}{$**$}
\begin{eqnarray}\label{zv2}
\overline
H_{s,y}(t)&=&\frac{1- F(\alpha+t-s)}{1-
F(\alpha)}\prod_{i=1}^k\frac{1-G^s(u_i+t-s)}{1-G^s(u_i)}\nonumber\\[-8pt]\\[-8pt]
&&\hspace*{0pt}{}\times\prod
_{j=1}^l\frac{1-G^r(z_j+t-s)}{1-G^r(z_j)}.\nonumber
\end{eqnarray}}

\noindent It is easy to see that $\overline H_{s,y}(t)$ satisfies
\[
\overline
H_{s,y}(u)=\overline H_{s,y}(t)\overline H_{t,\phi_{t-s}(y)}(u),\qquad
s\leq t\leq u.
\]
Then by Theorem 7.3.2 of \cite{jacobsen}, $V$ is
a piecewise deterministic Markov process constructed from
$\{(\tau_n,P_n), n\geq0\}$ using functions $\phi_t$ for the
deterministic part, survivor functions $\overline H_{s,y}$ for jump
time distributions and transition probabilities $r_t$ for the jumps.
Thus it follows from Theorem 7.5.1 of \cite{jacobsen} that $V$ is
a strong Markov process. The second part of the lemma follows directly
from the definition of the $\{P_t, t\geq0\}$ in
(\ref{dis:semiG}).
\end{pf}
\begin{remark}
For future purposes, we note that
the results of this paper including, in particular, the strong Markov
property established above, continue to be valid if the state component
$\alpha_E^{(N)}$ introduced in Section \ref{subs-modyn} is, instead, defined as follows:
\[
\alpha_E^{(N)}(s) \doteq \cases{
s, &\quad if $E^{(N)}(s) = 0$, \cr
\inf\bigl\{u>s\dvtx E^{(N)}(u)>E^{(N)}(s)\bigr\}-s, &\quad if $E^{(N)}(s) >
0$.}
\]
Observe that\vspace*{-1pt} when $E^{(N)}(s) > 0$, $\alpha_E^{(N)}(s)$ represents
the time from $s$ until the next arrival, and if $E^{(N)}$ is a
renewal process, then $\alpha_E^{(N)}$ is simply the forward recurrence time process.
A minor variation of the proof of Lemma \ref{lem:Mark} given above shows that
the strong Markov property holds in this case as well.  First, the
definition of $\phi_t(y)$ should be modified by replacing $\alpha + t$ by
$\alpha - t$ in (\ref{zv1}). With $V$, $r_t$ and
$p_1$, $p_2$, $p_3$  defined as before,
in this case,
the probability measure $r_t(y,\cdot)$ can be easily determined from $y$,
the distribution of the remaining time from $t$ to the next arrival
and $p_i$, $i=1,2,3$. Note that if  $\alpha > 0$ at time $t$,
then $p_1=0$. On the other hand, if $\alpha = 0$ at time $t$,
then $V$ jumps at time $t$ due to the arrival of a
new customer, and, hence, $p_1=1$.
Moreover, given that $V$
jumps at time $t$ from state~$y$, if the type of the jump at time
$t$ is known, then it is possible to determine the state to which
the process $V$ jumps.
For example,
suppose that the number $k$ in the expression for $y$ is less than $N$.
Then, at state $y$, there is at least one idle server. If the jump is
due to a new arrival, then the state $V$ will jump to the
region $\{c\in [0,\infty)\dvtx (c, x+1, \sum_{i=1}^k
\delta_{u_i}+\delta_0, \sum_{j=1}^l \delta_{z_i}+\delta_0)\}$
according\vspace*{1pt} to the distribution of the time to the next arrival
(which is determined  by the current state $\alpha$ due to the assumption that
$\alpha_E^{(N)}$ is Markov with respect to its own filtration).
Once again, the jump time dynamics are  captured by the
survivor functions, with the only difference that now
the  ratio $(1-F(\alpha+t-s))/(1-F(\alpha))$
on the right-hand side of (\ref{zv2}) should be replaced by
$\mathbh{1}_{\{\alpha\geq t-s\}}$.  The rest of the proof then follows
as before.
\end{remark}
\end{appendix}


%
\printaddresses

\mbox{}

\end{document}